\documentclass{cedram-aif}

% ...

\usepackage{amsfonts,latexsym,enumerate}
\usepackage{amsmath}
\usepackage{amscd}
\usepackage{float,amsmath,amssymb,mathrsfs,bm,multirow,graphics}
\usepackage[dvips]{graphicx}
\usepackage[percent]{overpic}

% ...

\newcommand{\R}{{\Bbb R}}

\newcommand{\C}{{\Bbb C}}

\newcommand{\D}{{\Bbb D}}

\newcommand{\diag}{\text{\upshape diag}}
\newcommand{\re}{\text{\upshape Re\,}}
\newcommand{\im}{\text{\upshape Im\,}}

\newcommand{\Arg}{\text{\upshape Arg\,}}

\newcommand{\proofend}{\hfill$\Box$\bigskip}

\newtheorem{claim}{Claim}
\newtheorem{figuretext}{Figure}

\newtheorem{assumptions}{Assumptions}

\title[Long-time asymptotics for the Degasperis--Procesi equation]
{Long-time asymptotics for the Degasperis--Procesi equation on the half-line}

\alttitle{Comportement asymptotique en temps grand de l'\'equation de Degasperis--Procesi sur la demi-droite}

\author{\firstname{Anne} \lastname{Boutet de Monvel}}
\address{Institut de Math\'ematiques de Jussieu-PRG \\
Universit\'e Paris Diderot \\
75205 Paris Cedex 13 (France)}
\email{anne.boutet-de-monvel@imj-prg.fr}

\author{\firstname{Jonatan} \lastname{Lenells}}
\address{Department of Mathematics \\
KTH Royal Institute of Technology \\
10044 Stockholm (Sweden)}
\email{jlenells@kth.se}

\author{\firstname{Dmitry} \lastname{Shepelsky}}
\address{Mathematical Division \\
Institute for Low Temperature Physics \\
61103 Kharkiv (Ukraine) and \\
School of Mathematics and Computer Sciences \\
V. N. Karazin Kharkiv National University \\
61022 Kharkiv (Ukraine)}
\email{shepelsky@yahoo.com}

\dedicatory{Dedicated to the memory of Louis Boutet de Monvel}

\thanks{The authors are grateful to the referee for many helpful suggestions. J.L. acknowledges support from the European Research Council, Consolidator Grant No. 682537, the Swedish Research Council, Grant No. 2015-05430, the G\"oran Gustafsson Foundation, Sweden, and the EPSRC, UK}

  
\altkeywords{\'Equation de Degasperis--Procesi, asymptotique en temps grand, probl\`eme de Riemann--Hilbert, probl\`eme aux limites.}


\begin{document}
%% Abstract 
\begin{abstract}
We analyze the long-time asymptotics for the Degasperis--Procesi equation on the half-line. By applying nonlinear steepest descent techniques to an associated $3 \times 3$-matrix valued Riemann--Hilbert problem, we find an explicit formula for the leading order asymptotics of the solution in the similarity region in terms of the initial and boundary values. 
\end{abstract}

%% French abstract
\begin{altabstract}
Nous \'etudions le comportement asymptotique en temps grand de l'\'equation de Degasperis--Procesi sur la demi-droite. L'appli\-cation de techniques de descente de plus grande pente non lin\'eaire \`a un probl\`eme de Riemann--Hilbert matriciel $3\times 3$ associ\'e nous permet d'obtenir une formule explicite, en termes des donn\'ees initiale et au bord, pour le terme dominant de l'asymptotique de la solution dans la r\'egion de similarit\'e.
\end{altabstract}

\maketitle

\vspace{-1cm}

\setcounter{tocdepth}{1}
\tableofcontents

\section{Introduction}
The nonlinear steepest descent method introduced in \cite{DZ1993} provides a powerful technique for determining asymptotics of solutions of nonlinear integrable PDEs. By appropriately deforming the contour of the associated Riemann--Hilbert (RH) problem, the long-time behavior of the solution can be determined by adding up the contributions from the individual critical points. In this way the asymptotics associated with the modified KdV \cite{DZ1993}, the nonlinear Schr\"odinger \cite{DIZ1993}, and several other integrable equations posed on the real line have been rigorously established, see \cite{BKST2009, DVZ1994, K1993}. More recently, a number of works treating periodic problems \cite{KT2012} as well as initial-boundary value problems \cite{BFS2004, BS2009} have also appeared. 

In this paper we use the method of nonlinear steepest descent to analyze long-time asymptotics for the Degasperis--Procesi (DP) equation 
\begin{equation}\label{DP}
  u_t - u_{txx} + 3\kappa u_x+4uu_x - 3u_xu_{xx} - uu_{xxx}=0,\qquad \kappa > 0,
\end{equation}
posed in the domain
\begin{equation}\label{halflinedomain}
   \Omega = \left\{(x,t) \in \R^2\, |\, 0 \leq x < \infty, \; 0 \leq t < \infty \right\}.
\end{equation}
Our main result (see Theorem \ref{mainth} below) gives an explicit formula for the leading order asymptotics of $u(x,t)$ in the similarity region $0 < \frac{x}{t} < 3$ in terms of the initial and boundary values. 
In this region it has the form of slowly decaying oscillations, whereas in the complementary region $\frac{x}{t} > 3$ it is dominated by solitons, if any, see \cite{BKST2009,BS2009,BS2013}.

Equation (\ref{DP}) was discovered in \cite{DP1999} using methods of asymptotic integrability. A Lax pair and a bi-Hamiltonian structure were derived in \cite{DHH2002}.
An interesting aspect of (\ref{DP}) is the existence of peaked solutions \cite{DHH2002} as well as weak solutions with a very low degree of regularity \cite{CK2006}. The latter class  includes a class of discontinuous generalizations of the peakons called shock-peakons \cite{Lu2007}. The asymptotic behavior of the solution of (\ref{DP}) on the line was determined in \cite{BS2013}.  
In \cite{L2013} the solution of the initial-boundary value problem of (\ref{DP}) on the half-line was expressed in terms of the solution of a $3 \times 3$-matrix RH problem. 

Compared with most other applications of the nonlinear steepest descent approach, the asymptotic analysis of (\ref{DP}) presents a number of additional difficulties:

\begin{enumerate}[$(a)$]
\item The RH problem associated with (\ref{DP}) involves $3\times 3$ matrices instead of $2 \times 2$ matrices. This implies that the standard uniqueness results for $L^2$-RH problems (such as Theorem 7.18 of \cite{D1999}) do not apply. However, it turns out that in an appropriate function space, which we denote by $\dot{L}^3$, uniqueness holds also for $3\times 3$-matrix valued RH problems, see \cite{LCarleson}. 
%Thus, by developing the nonlinear steepest descent approach in the $\dot{L}^3$-setting, rigorous asymptotic formulas can still be obtained. 

\vspace{.1cm}
\item The $t$-part of the Lax pair associated with (\ref{DP}) has singularities at the points $K_j = e^{\frac{\pi i j}{3} - \frac{\pi i}{6}}$, $j = 1, \dots, 6$. In \cite{L2013} this difficulty was overcome by utilizing two different sets of eigenfunctions which were solutions of two different Lax pairs (a similar idea was used already in \cite{BS2013} to recover $u(x,t)$ for the problem on the line). Here we adopt a similar approach; however, in order to obtain a RH problem suitable for the asymptotic analysis of (\ref{DP}), we use a modification of the RH problem in \cite{L2013}. The modified problem has the advantage that, after the appropriate contour deformations prompted by the nonlinear steepest descent method have been performed, the RH problem involves only one set of eigenfunctions near each of the twelve critical points. This leads to a jump matrix near each critical point of an appropriate form.

\vspace{.1cm}
\item The Lax pair associated with (\ref{DP}) has singularities at the sixth roots of unity $\varkappa_j = e^{\frac{\pi i(j-1)}{3}}$, $j = 1, \dots, 6$. In \cite{BS2013, L2013} this difficulty was overcome by considering a regular RH problem for an associated row vector. Here, rather than trying to develop a nonlinear steepest descent approach for row vector RH problems, we carry out the steepest descent analysis using a regular $3 \times 3$-matrix valued solution which, in general, is different from the original solution. However, by uniqueness for the row vector RH problem, the row vectors associated with these two solutions coincide. 

\vspace{.1cm}
\item On the half-line, the jump contour for the RH problem associated with (\ref{DP}) involves nontransversal intersection points, see Figure \ref{Dns.pdf}. This implies that the standard theory of $L^p$-RH problems does not apply. We circumvent this difficulty by employing the theory of $L^p$-RH problems developed in \cite{LCarleson} for general Carleson jump contours. 
\end{enumerate}

Our analysis determines the asymptotic behavior of $u(x,t)$ provided that all boundary values $\{\partial_x^ju(0, t)\}_0^2$ are known. However, for a well-posed problem, only a subset of the initial and boundary values can be independently prescribed. If all boundary values are not known, our asymptotic formula (see Theorem \ref{mainth}) still provides some information on the solution, but since the function $r(k)$ is unknown, the precise form of the asymptotics remains undetermined. In general, the computation of the unknown boundary values (i.e. the construction of the generalized Dirichlet-to-Neumann map) involves the solution of a nonlinear Volterra integral equation. We do not consider the construction of the Dirichlet-to-Neumann map in this paper. We also do not consider the existence of so-called linearizable boundary conditions for which the unknown boundary values can be eliminated thanks to additional symmetries. 

In Section \ref{prelsec}, we give a short review of the RH approach for (\ref{DP}) on the half-line.
In Section \ref{RHsec}, we formulate a RH problem suitable for determining the long-time asymptotics. 
In Section \ref{steepsec}, we prove a nonlinear steepest descent theorem appropriate for analyzing the asymptotics in the similarity region. 
In Section \ref{similaritysec}, we prove our main theorem.

\section{Preliminaries}\label{prelsec}
We consider initial-boundary value problems for (\ref{DP}) for which the initial and boundary values 
\begin{subequations}\label{boundaryconditions}
\begin{align}
& u_0(x) = u(x,0), \qquad x \geq 0,
	\\ 
& g_0(t) = u(0,t), \quad g_1(t) = u_x(0,t), \quad g_2(t) = u_{xx}(0,t), \qquad t \geq 0,
\end{align}
\end{subequations}
satisfy the three conditions
\begin{subequations}\label{qassumptions}
\begin{align} 
 & u_0(x) - u_{0xx}(x) + \kappa > 0, \qquad x \geq 0, 
	\\
&  g_0(t) - g_2(t) + \kappa > 0, \qquad t \geq 0,
	\\
&  g_0(t) \leq 0,  \qquad t \geq 0.
\end{align}
\end{subequations}
The assumptions in (\ref{qassumptions}) imply the following positivity condition which is needed for the spectral analysis:
\begin{align}\label{qpositive}
  u(x,t) - u_{xx}(x,t) + \kappa > 0, \qquad (x,t) \in \Omega.
\end{align}
In view of (\ref{qpositive}), we may define $q(x,t)$ by
\begin{equation}\label{qdef}
  q(x,t) = \big(u(x,t) - u_{xx}(x,t) + \kappa\big)^{\frac{1}{3}}, \qquad (x,t) \in \Omega.
\end{equation}

We next give a short review of the RH approach for (\ref{DP}) on the half-line; see \cite{L2013} for further details. 
We suppose that $\{g_j\}_0^2$ belong to the Schwartz class $\mathcal{S}(\R_+)$ and that there exists a unique smooth solution $u(x,t)$ of (\ref{DP}) in $\Omega$ such that (\ref{boundaryconditions}) and (\ref{qassumptions}) are satisfied and $u(\cdot, t) \in \mathcal{S}(\R_+)$ for each $t \geq 0$.
For simplicity, we henceforth assume that $\kappa = 1$.

\subsection{Lax pairs}
Equation (\ref{DP}) admits the Lax pair \cite{BS2013, CIL2010}
\begin{equation}\label{psilax}
\begin{cases}
  \psi_x(x,t,k) = L(x,t,k) \psi(x,t,k), \\
  \psi_t(x,t,k) = Z(x,t,k) \psi(x,t,k),
\end{cases}
\end{equation}
where $k \in \hat{\C} = \C \cup \{\infty\}$ is the spectral parameter, $\psi(x,t, k)$ is a $3\times 3$-matrix valued eigenfunction, the $3\times 3$-matrix valued functions $L$ and $Z$ are defined by
$$L(x,t, k) = \begin{pmatrix} 0 & 1 & 0 \\ 0 & 0 & 1 \\ \lambda q^3 & 1 & 0 \end{pmatrix},
\quad
Z(x,t, k) = \begin{pmatrix} u_x - \frac{2}{3\lambda} & - u & \frac{1}{\lambda} \\ u + 1 & \frac{1}{3\lambda} & -u \\
u_x - \lambda u q^3 & 1 & - u_x + \frac{1}{3\lambda} \end{pmatrix},$$
and $\lambda = \lambda(k)$  is defined by
$$\lambda = \frac{1}{3\sqrt{3}}\biggl(k^3 + \frac{1}{k^3}\biggr).$$
Let $\omega = e^{\frac{2\pi i}{3}}$. Define $\{l_j\}_1^3$ and $\{z_j\}_1^3$ by 
\begin{align}\label{lmexpressions}
&l_j(k) = \frac{1}{\sqrt{3}}\left(\omega^j k + \frac{1}{\omega^j k}\right), \quad
z_j(k) = \sqrt{3}\left(\frac{(\omega^j k)^2 + (\omega^j k)^{-2}}{k^3 + k^{-3}}\right), \quad \, k \in \C.
\end{align}
Let
\begin{align}\label{Pdef}
P(k) = \begin{pmatrix} 1 & 1 & 1 \\ l_1(k) & l_2(k) & l_3(k) \\ l_1^2(k) & l_2^2(k) & l_3^2(k) \end{pmatrix}, \qquad k \in \C,
\end{align}
and define $\{V_j(x,t,k), \tilde{V}_j(x,t,k)\}_1^2$ by
\begin{align*}
& V_1 = P^{-1}\begin{pmatrix} 0 & 0 & 0 \\ 0 & 0 & 0 \\ \lambda (q^3 - 1) & 0 & 0 \end{pmatrix}P,
	\\
& V_2 = P^{-1}\begin{pmatrix} u_x & - u & 0 \\ u & 0 & -u \\ u_x - \lambda u q^3 & 0 & -u_x \end{pmatrix}P,
 	\\
& \tilde{V}_1 = P^{-1}\begin{pmatrix} \frac{q_x}{q} & 0 & 0 \\ 0 & 0 & 0 \\ 0 & \frac{1}{q} - q & - \frac{q_x}{q} \end{pmatrix}P,
	\\
& \tilde{V}_2 = P^{-1} \left[\begin{pmatrix} -\frac{uq_x}{q} & 0 & 0 \\ \frac{u + 1}{q} - 1 & 0 & 0 \\ \frac{u_x}{q^2} & \frac{1}{q} - 1 + uq & \frac{uq_x}{q} \end{pmatrix} + \frac{q^2 - 1}{\lambda}\begin{pmatrix} 0 & 0 & 1 \\ 0 & 0 & 0 \\ 0 & 0 & 0 \end{pmatrix}\right] P.
\end{align*}
Let $\mathcal{L} = \diag(l_1 , l_2 , l_3 )$ and $\mathcal{Z}  = \diag(z_1 , z_2 , z_3 )$. The eigenfunctions $\Psi$ and $\tilde{\Psi}$ introduced by
\begin{subequations}\label{PhitildePhidef}
\begin{align}\label{Phidef}
 & \psi(x,t,k) = P(k)\Psi(x,t,k) e^{\mathcal{L}(k) x + \mathcal{Z}(k) t},
  	\\\label{tildePhidef}
&  \psi(x,t,k) = \mathcal{D}(x,t) P(k) \tilde{\Psi}(x,t,k) e^{\mathcal{L}(k)y(x,t) + \mathcal{Z}(k) t},	
\end{align} 
\end{subequations} 
where 
\begin{align}\label{ydef}
&  y(x,t) = \int_{(0,0)}^{(x,t)} q(x',t')\left(dx' - u(x', t')dt'\right), 
  	\\\nonumber
&  \mathcal{D}(x,t) = \begin{pmatrix} \frac{1}{q(x,t)} & 0 & 0 \\ 0 & 1 & 0 \\ 0 & 0 & q(x,t) \end{pmatrix},
\end{align}  
satisfy the Lax pair equations
\begin{subequations}\label{laxpairs}
\begin{align}\label{Philax}
& \begin{cases}
  \Psi_x - [\mathcal{L}, \Psi] = V_1 \Psi, \\
  \Psi_t - [\mathcal{Z}, \Psi] = V_2 \Psi,
\end{cases}
\end{align}
and 
\begin{align}\label{tildePhilax} 
& \begin{cases}
  \tilde{\Psi}_x - [q\mathcal{L}, \tilde{\Psi}] = \tilde{V}_1 \tilde{\Psi}, \\
  \tilde{\Psi}_t - [\mathcal{Z} - uq\mathcal{L}, \tilde{\Psi}] = \tilde{V}_2 \tilde{\Psi},
\end{cases}
\end{align}
\end{subequations}
respectively. 

\subsection{Analytic eigenfunctions}
\begin{figure}
\begin{center}
 \begin{overpic}[width=.65\textwidth]{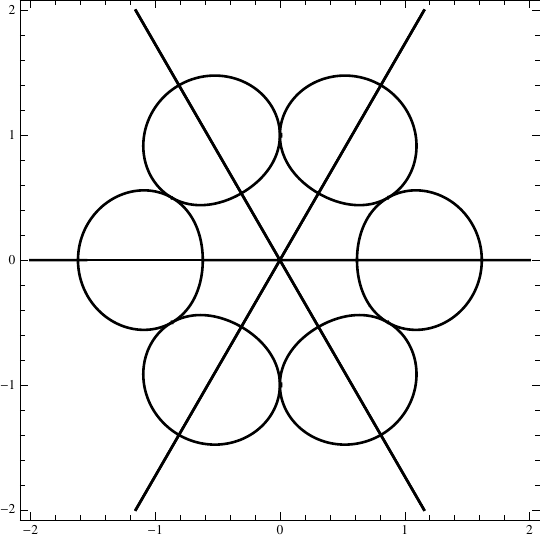}
      \put(86,72){$D_1$}
      \put(59,55){$D_1$}
      \put(50,90){$D_2$}
      \put(50,61){$D_2$}
     \put(14,72){$D_3$}
      \put(40,55){$D_3$}
     \put(14,28){$D_4$}
      \put(40,45){$D_4$}
      \put(50,10){$D_5$}
      \put(50,39){$D_5$}
      \put(86,28){$D_6$}
      \put(59,45){$D_6$}
      \put(74,56){$D_7$}
      \put(67,68){$D_8$}
      \put(56,74){$D_9$}
      \put(42,74){$D_{10}$}
      \put(32,68){$D_{11}$}
      \put(25,56){$D_{12}$}
      \put(25,44){$D_{13}$}
      \put(32,32){$D_{14}$}
      \put(42,26){$D_{15}$}
      \put(56,26){$D_{16}$}
      \put(67,32){$D_{17}$}
      \put(74,44){$D_{18}$}
 \end{overpic}
    \qquad \qquad
     \begin{figuretext}\label{Dns.pdf}
       The sets $D_n$, $n = 1, \dots, 18$, which decompose the complex $k$-plane. 
         \end{figuretext}
     \end{center}
\end{figure}

Let $\gamma_j$, $j = 1,2,3$, denote contours in the $(x, t)$-plane connecting $(x_j, t_j)$ with $(x,t)$, where $(x_1, t_1) = (0, \infty)$, $(x_2, t_2) = (0, 0)$, and $(x_3, t_3) = (\infty, t)$. The contours can be chosen to consist of straight line segments parallel to the $x$- or $t$-axis. For a diagonal matrix $D$, let $\hat{D}$ denote the operator which acts on a matrix $A$ by $\hat{D}A = [D, A]$, i.e. $e^{\hat{D}}A = e^D A e^{-D}$. 
We define solutions $\{\Psi_n(x,t,k)\}_1^{18}$ and $\{\tilde{\Psi}_n(x,t,k)\}_1^{18}$ of the Lax pairs (\ref{Philax}) and (\ref{tildePhilax}) respectively, by the solutions of the integral equations
\begin{subequations}\label{Fredholm}
\begin{align}\label{PhiFredholm}
& (\Psi_n)_{ij}(x,t,k) = \delta_{ij} + \int_{\gamma_{ij}^n} \left(e^{\hat{\mathcal{L}}(k) x + \hat{\mathcal{Z}}(k) t} W_n(x',t',k)\right)_{ij}, 
	\\ \label{tildePhiFredholm}
& (\tilde{\Psi}_n)_{ij}(x,t,k) = \delta_{ij} + \int_{\gamma_{ij}^n} \left(e^{\hat{\mathcal{L}}(k) y(x,t) + \hat{\mathcal{Z}}(k) t} \tilde{W}_n(x',t',k)\right)_{ij},  
\end{align}
\end{subequations}
where $k \in D_n$, $i,j = 1, 2,3$, $n = 1, \dots, 18$, and the contours $\gamma^n_{ij}$ are given by
 \begin{align} \label{gammaijnudef}
 \gamma_{ij}^n =  \begin{cases}
 \gamma_1,  \qquad \re l_i(k) < \re l_j(k), \quad \re z_i(k) \geq \re z_j(k),
	\\
\gamma_2,  \qquad \re l_i(k) < \re l_j(k),\quad \re z_i(k) < \re z_j(k),
	\\
\gamma_3,  \qquad \re l_i(k) \geq \re l_j(k),
	\\
\end{cases} \ \ \text{for} \ k \in D_n,
\end{align}
the closed one-forms $W_n(x,t,k)$ and $\tilde{W}_n(x,t,k)$ are defined by
\begin{align*}  
  W_n = e^{-\hat{\mathcal{L}} x - \hat{\mathcal{Z}} t}(V_1 dx + V_2 dt) \Psi_n,
  \qquad   \tilde{W}_n = e^{-\hat{\mathcal{L}} y- \hat{\mathcal{Z}} t}(\tilde{V}_1 dx + \tilde{V}_2 dt) \tilde{\Psi}_n,
\end{align*}  
and the open sets $\{D_n\}_1^{18}$ are displayed in Figure \ref{Dns.pdf}. 
Precise definitions of all the sets $D_n$ can be found in \cite{L2013}; here we only give the definitions of the $D_n$ relevant near the positive real axis and near $K_1$: 
\begin{align*}
&D_1 = \{k \in \hat{\C}\,|\, \text{Re}\,l_1 < \text{Re}\,l_2 < \text{Re}\,l_3 \text{  and  } 
\text{Re}\,z_1 < \text{Re}\,z_2 < \text{Re}\,z_3 \},
	\\
&D_6 = \{k \in \hat{\C}\,|\, \text{Re}\,l_2 < \text{Re}\,l_1 < \text{Re}\,l_3 \text{  and  } 
\text{Re}\,z_2 < \text{Re}\,z_1 < \text{Re}\,z_3 \},
	\\
&D_7 = \{k \in \hat{\C}\,|\, \text{Re}\,l_1 < \text{Re}\,l_2 < \text{Re}\,l_3 \text{  and  } 
\text{Re}\,z_2 < \text{Re}\,z_1 < \text{Re}\,z_3 \},
	\\
&D_8 = \{k \in \hat{\C}\,|\, \text{Re}\,l_1 < \text{Re}\,l_2 < \text{Re}\,l_3 \text{  and  } 
\text{Re}\,z_1 < \text{Re}\,z_3 < \text{Re}\,z_2 \},	
	\\
&D_{18} = \{k \in \hat{\C}\,|\, \text{Re}\,l_2 < \text{Re}\,l_1 < \text{Re}\,l_3 \text{  and  } 
\text{Re}\,z_1 < \text{Re}\,z_2 < \text{Re}\,z_3 \}.
\end{align*}

\begin{figure}
\begin{center}
 \begin{overpic}[width=.65\textwidth]{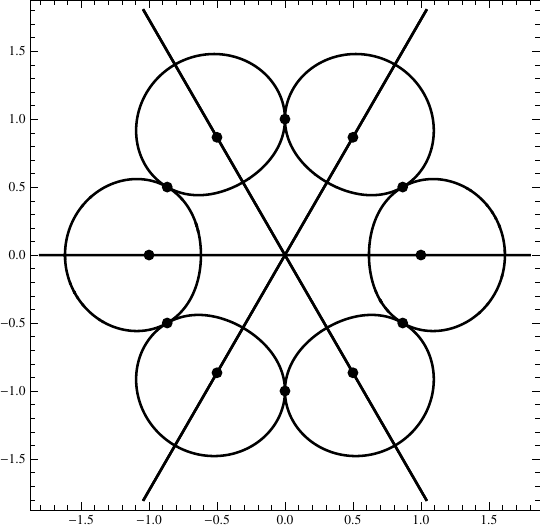}
      \put(75,59){\small $K_1$}
      \put(54.5,74.5){\small $K_2$}
      \put(31,64.5){\small $K_3$}
      \put(26,40){\small $K_4$}
      \put(46,24){\small $K_5$}
      \put(70,34){\small $K_6$}
            \put(75,46.5){ $\varkappa_1$}
      \put(65,70){ $\varkappa_2$}
      \put(40,73){ $\varkappa_3$}
      \put(24,53){ $\varkappa_4$}
      \put(33,29){ $\varkappa_5$}
      \put(59,26){ $\varkappa_6$}
 \end{overpic}
    \qquad \qquad
     \begin{figuretext}\label{Kjs.pdf}
       The points $K_j = e^{\frac{\pi i j}{3} - \frac{\pi i}{6}}$, $j = 1, \dots, 6$, where $\lambda = 0$, and the points $\varkappa_j = e^{\frac{\pi i(j-1)}{3}}$, $j = 1, \dots, 6$, where $P^{-1}(k)$ has poles.
            \end{figuretext}
     \end{center}
\end{figure}

Let $K_j = e^{\frac{\pi i j}{3} - \frac{\pi i}{6}}$, $j = 1, \dots, 6$, denote the points where $\lambda = 0$ and let $\varkappa_j = e^{\frac{\pi i(j-1)}{3}}$, $j = 1, \dots, 6$, denote the sixth roots of unity, see Figure \ref{Kjs.pdf}.
Away from the sets $\{\infty, 0\} \cup \{\varkappa_j\}_1^6 \cup \{k_j\}$ and $\{\varkappa_j, K_j\}_1^6 \cup \{k_j\}$, respectively, $\Psi_n$ and $\tilde{\Psi}_n$ are bounded and analytic functions of $k \in D_n$ with continuous extensions to $\bar{D}_n$. Here $\{k_j\}$ denotes a possibly empty set of singularities at which the Fredholm determinant of integral equations (\ref{Fredholm}) vanishes; for simplicity, we henceforth assume that the set $\{k_j\}$ is empty (solitonless case). 
For those $n$ for which the indicated limiting points lie on the boundary of the corresponding $D_n$,
\begin{align*}
& \Psi_n(x,t,k) = I + O(k - K_j) \quad \text{as} \quad k \to K_j, \ k \in D_n, \ j = 1, \dots, 6,
	\\
&  \tilde{\Psi}_n(x,t,k) = I + O(1/k) \quad \text{as} \quad k \to \infty, \ k \in D_n, 
 	\\
&  \tilde{\Psi}_n(x,t,k) = I + O(k) \quad \text{as} \quad k \to 0, \ k \in D_n,
\end{align*}
where $I$ denotes the identity matrix.

We define spectral functions $\{S_n(k)\}_1^{18}$ and $\{\tilde{S}_n(k)\}_1^{18}$ by
\begin{align}\label{Sndef}
S_n(k) = \Psi_n(0,0,k), \qquad \tilde{S}_n(k) = \tilde{\Psi}_n(0,0,k), \qquad k \in D_n.
\end{align}

\subsection{Symmetries}
Define sectionally analytic functions $S_*(k)$ and $\tilde{S}_*(k)$ for $k \in \C$ by setting $S_*(k) = S_n(k)$ and $\tilde{S}_*(k) = \tilde{S}_n(k)$ for $k \in D_n$. If $f$ denotes one of the $3 \times 3$-matrix valued functions $\mathcal{L}$, $\mathcal{Z}$, $M$, $S_*$, or $\tilde{S}_*$, then $f$ obeys the symmetries
\begin{subequations}\label{symmetries}
\begin{align}\label{symmetriesa}
&  f(k) = \mathcal{A} f(\omega k)\mathcal{A}^{-1}, \qquad k \in \C, 
	\\ \label{symmetriesb}
& f(k) = \mathcal{B} f(1/k)\mathcal{B}, \qquad k \in \C,	
	\\ \label{symmetriesc}
& f(k) = \mathcal{B} \overline{f(\overline{k})}\mathcal{B},	\qquad k \in \C,
\end{align}
\end{subequations}
where $\mathcal{A}$, $\mathcal{B}$ are defined by
\begin{align}\label{ABdef}
\mathcal{A} =  \begin{pmatrix}
0 & 0 & 1 \\
 1 & 0 & 0 \\
 0 & 1 & 0\end{pmatrix}, \qquad
  \mathcal{B} =  \begin{pmatrix}
0 & 1 & 0 \\
1 & 0 & 0 \\
0 & 0 & 1 \end{pmatrix}.
\end{align}

\begin{figure}
\begin{center}
 \begin{overpic}[width=.65\textwidth]{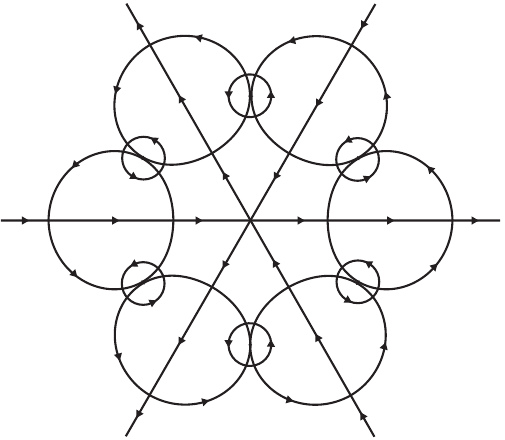}
      \put(86,64){\small $E_1$}
      \put(58,48){\small $E_1$}
      \put(48,82){\small $E_2$}
      \put(48,53){\small $E_2$}
     \put(12,64){\small $E_3$}
      \put(38,48){\small $E_3$}
     \put(12,22){\small $E_4$}
      \put(38,37){\small $E_4$}
      \put(48,2){\small $E_5$}
      \put(48,31){\small $E_5$}
      \put(86,22){\small $E_6$}
      \put(58,37){\small $E_6$}
      \put(75,48){\small $E_7$}
      \put(67,65){\small $E_8$}
      \put(57,70){\small $E_9$}
      \put(38,70){\small $E_{10}$}
      \put(28,65){\small $E_{11}$}
      \put(19,48){\small $E_{12}$}
      \put(19,37){\small $E_{13}$}
      \put(28,20){\small $E_{14}$}
      \put(38,14){\small $E_{15}$}
      \put(57,14){\small $E_{16}$}
      \put(67,20){\small $E_{17}$}
      \put(75,37){\small $E_{18}$}
 \end{overpic}
    \qquad \qquad
     \begin{figuretext}\label{Econtour.pdf}
       The sets $E_n$ which decompose the complex $k$-plane. 
         \end{figuretext}
     \end{center}
\end{figure}

\begin{figure}
\begin{center}
 \begin{overpic}[width=.42\textwidth]{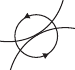}
      \put(85,66){\small $E_1$}
      \put(7,22){\small $E_1$}
      \put(66,6){\small $E_7$}
      \put(27,84){\small $E_8$}
      \put(24.2,34){\tiny $E_{19}$}
      \put(65.5,56.7){\tiny $E_{19}$}
      \put(51,31){\small $E_{25}$}
      \put(40,58){\small $E_{26}$}
 \end{overpic}
    \qquad \qquad
     \begin{figuretext}\label{EnearK1.pdf}
       The sets $E_n$ for $k$  near $K_1$. 
         \end{figuretext}
     \end{center}
\end{figure}

\section{A Riemann--Hilbert problem}\label{RHsec}
We use the eigenfunctions $\Psi_n$ and $\tilde{\Psi}_n$ to define a RH problem suitable for analyzing the long-time asymptotics. 

Choose a small radius $r>0$ and let $B_j$ denote the open disk of radius $r$ centered at $K_j$, $j = 1, \dots, 6$. Let $B = \cup_{j=1}^6 B_j$ and define open sets $\{E_n\}_1^{36}$ by, see Figures \ref{Econtour.pdf} and \ref{EnearK1.pdf}:
\begin{align}\label{EnDnU}
E_n = D_n \setminus \bar{B}, \qquad  E_{n + 18} = D_n \cap B, \qquad n = 1, \dots, 18.
\end{align}
The eigenfunctions $\{\tilde{\Psi}_n\}_{1}^{18}$ are well-behaved near $k = \infty$ and $k = 0$ while the eigenfunctions $\{\Psi_n\}_1^{18}$ are well-behaved near the $K_j$'s. We formulate a RH problem relative to the contour shown in Figure \ref{Econtour.pdf} (see also Figure \ref{EnearK1.pdf}) by using $\tilde{\Psi}_n$ and $\Psi_n$  for $k$ in $E_n$ and $E_{n +18}$, respectively. 

Let $y = y(x,t)$ be the function defined in (\ref{ydef}). The map $\mathcal{G}:(x,t) \mapsto (y,t)$ is a bijection from $\Omega = \{x \geq 0, t \geq 0\}$ onto $\mathcal{G}(\Omega) \subset \R^2$. Thus, for each $(y,t) \in \mathcal{G}(\Omega)$, we may define a sectionally meromorphic function $M(y,t,k)$ by
\begin{align}\label{Mndef}
M(y, t, k) = \begin{cases}
\tilde{\Psi}_n(x,t,k), & k \in E_{n},
	\\
P(k)^{-1}\mathcal{D}(x,t)^{-1} P(k) \Psi_n(x, t, k)e^{(x-y + \nu_0)\mathcal{L}(k)}, & k \in E_{n + 18},  
\end{cases}
\end{align}
where $n = 1, \dots, 18$ and the constant $\nu_0 \in \R$ is defined by
$$\nu_0 = \lim_{x \to\infty} (y-x) = \int_0^\infty(q(x,0) - 1) dx.$$ 
Let $M_n$ denote the restriction of $M$ to $E_n$.
The definition (\ref{Mndef}) and the relations (\ref{PhitildePhidef}) imply that $M$ satisfies the jump condition
\begin{align}\label{jumpcond}
M_n = M_m J_{m, n}, \qquad k \in \bar{E}_n \cap \bar{E}_m,
\end{align}
where  
\begin{align}\label{Jmndef}
\begin{cases}
J_{m, n}(y,t,k) = e^{y \hat{\mathcal{L}}  + t \hat{\mathcal{Z}}} \big(\tilde{S}_m^{-1}(k)\tilde{S}_n(k)\big),
	\\
J_{n, n+18}(y,t,k) = e^{y \hat{\mathcal{L}}  + t \hat{\mathcal{Z}}} C_n(k),
	\\
J_{m + 18,n + 18}(y,t,k) = e^{y \hat{\mathcal{L}}  + t \hat{\mathcal{Z}}}  e^{-\nu_0\hat{\mathcal{L}}}\big(S_m^{-1}(k)S_n(k)\big),
\end{cases}
\ n,m = 1, \dots, 18.
\end{align}
The functions $\{C_n(k)\}_1^{18}$ are defined as follows. By (\ref{PhitildePhidef}) the functions
\begin{align*}
  \psi_n & = DP\tilde{\Psi}_ne^{y\mathcal{L} + t\mathcal{Z}} = DPM_n e^{y\mathcal{L} + t\mathcal{Z}},
  	\\
  \psi_{n+18} & = P\Psi_n e^{x\mathcal{L} + t\mathcal{Z}} = DP M_{n+18} e^{y\mathcal{L} + t\mathcal{Z} - \nu_0 \mathcal{L}}
  \end{align*}
solve the same differential equations (\ref{psilax}), hence $\psi_{n+18} = \psi_n \tilde{C}_n(k)$ with  $\tilde{C}_n(k)$ independent of $(x,t)$. Thus, $M_{n+18} = M_n e^{y\hat{\mathcal{L}} + t\hat{\mathcal{Z}}}C_n(k)$ with $C_n(k) = \tilde{C}_n(k) e^{\nu_0\mathcal{L}}$, and $C_n(k) = \psi_n^{-1} \psi_{n+18} e^{\nu_0\mathcal{L}}$ satisfies
\begin{align}\label{Cndef}
C_n(k) = e^{-y\mathcal{L} - t\mathcal{Z}} \tilde{\Psi}_n(x,t,k)^{-1} P(k)^{-1} \mathcal{D}(x,t)^{-1} P(k) \Psi_n(x,t,k) e^{x \mathcal{L} + t \mathcal{Z}} e^{\nu_0 \mathcal{L}}.
\end{align}

\begin{prop}\label{Mnprop}
Let $E = \cup_{n=1}^{36} E_n$ and let $(y,t) \in \mathcal{G}(\Omega)$. Except for possible singularities at the points $\{\varkappa_j\}_1^6$, $M(y,t,k)$ is a bounded and analytic function of $k \in E$. Moreover, 
$$M(x,t,k) = I + O(1/k) \quad \text{uniformly as} \quad k \to \infty, \quad k \in \C.$$
\end{prop}
\begin{proof}
Since $P(k)^{-1}\mathcal{D}(x,t)^{-1} P(k)$ is an analytic function of $k \in \hat{\C}$ except for poles at the points $\varkappa_j$, the result  follows immediately from the properties of the functions $\Psi_n$ and $\tilde{\Psi}_n$. 
\end{proof}

The singularity structure of $M_n$ at the $\varkappa_j$'s implies that the function $N$ defined by
\begin{equation}\label{nundef}
N(y, t, k) = (1,1,1) M(y,t,k), \qquad k \in E,
\end{equation}
is nonsingular at the $\varkappa_j$'s, see \cite{BS2013, L2013}. Together with Proposition \ref{Mnprop} and Lemma \ref{EpCnlemma}, this implies that, for each $(y,t) \in \mathcal{G}(\Omega)$, 
$$N(y, t, \cdot) \in (1,1,1) + \dot{E}^3(E) \cap E^\infty(E),$$
where the function spaces $\dot{E}^3(E)$ and $E^\infty(E)$ are defined in Appendix \ref{RHapp}. Thus the following result holds.

\begin{prop}\label{Nprop}
For each $(y,t) \in \mathcal{G}(\Omega)$, the function $N(y,t, \cdot)$ is a row vector solution of the $L^3$-RH problem
\begin{align}\label{RHnu}
\begin{cases} N(y,t, \cdot) \in (1,1,1) + \dot{E}^3(E), 
	\\
N_n(y,t,k) =  N_m(y,t,k) J_{m,n}(y,t,k) \quad \text{for a.e.} \ k \in \bar{E}_n \cap \bar{E}_m,
\end{cases} 
\end{align}
where $n,m= 1, \dots, 36$.
\end{prop}

\begin{rema}
Although the RH problem formulated in (\ref{RHnu}) differs from the one used in \cite{L2013}, both problems rely on the same idea of using the $\Psi_n$'s near the $K_j$'s and the $\tilde{\Psi}_n$'s near $\{0, \infty\}$. The RH problem formulated in (\ref{RHnu}) is better adapted for our present purposes because it  uses only one set of eigenfunctions, namely the $\tilde{\Psi}_n$'s, near the three lines $\R$, $\omega \R$, and $\omega^2 \R$.
\end{rema}

\subsection{Jump matrix}
Define functions $\{r(k), h(k), \check{r}(k), \check{h}(k), r_1(k), \check{r}_1(k)\}$ by
\begin{align}\nonumber
& r(k) = (\tilde{S}_{18}(k)^{-1}\tilde{S}_7(k))_{21}, \qquad k \in \bar{E}_7 \cap \bar{E}_{18},
	\\\nonumber
& h(k) = (\tilde{S}_{6}(k)^{-1}\tilde{S}_{18}(k))_{21}, \qquad k \in \bar{E}_{18},
	\\\nonumber
& \check{r}(k) = ( \tilde{S}_{13}(k)^{-1}\tilde{S}_{12}(k))_{21}, \qquad k \in \bar{E}_{12} \cap \bar{E}_{13},
	\\\nonumber
& \check{h}(k) = (\tilde{S}_{4}(k)^{-1}\tilde{S}_{13}(k))_{21}, \qquad k \in \bar{E}_{13},
	\\\nonumber
& r_1(k) % = \frac{\tilde{s}_{21}(k)}{\tilde{s}_{22}(k)}
= (\tilde{S}_{6}(k)^{-1}\tilde{S}_1(k))_{21}, \qquad k \in \R_+,
	\\ \label{rhdef}
& \check{r}_1(k)  % = \frac{m_{12}(\tilde{s}(k))}{m_{11}(\tilde{s}(k))}
= (\tilde{S}_{4}(k)^{-1}\tilde{S}_3(k))_{21}, \qquad k \in \R_-.
\end{align}
The domains of definition of the functions $h(k)$, $\check{h}(k)$, $r_1(k)$, and $\check{r}_1(k)$ in (\ref{rhdef}) can be understood as follows.
For $j = 1, \dots, 6$, the function $\tilde{S}_j(k)$ is defined in terms of the initial data alone. This means that $\tilde{S}_j(k)$ has an analytic continuation to the sector $\arg k \in ((j-1)\pi/3, j\pi/3)$ for each $j = 1, \dots, 6$. 
It follows that $h(k)$  and $\check{h}(k)$  are well-defined for $k \in \bar{E}_{18}$ and $k \in \bar{E}_{13}$, respectively. Similarly, $r_1(k)$ and $\check{r}_1(k)$ are well-defined for $k \in \R_+$ and $k \in \R_-$, respectively.

The approach of Section 5 of \cite{L2013} shows that with the contour oriented as in Figure \ref{Econtour.pdf} the jump matrix for the RH problem (\ref{RHnu}) is given for $k$ near $\R$ by
\begin{align}\label{Jformula1}
J = \begin{cases} 
J_{6,1} = 
\begin{pmatrix}
 1 & - \overline{r_1(\bar{k})} e^{-t\Phi} & 0 \\
r_1(k)e^{t\Phi} & 1 - |r_1(k)|^2 & 0 \\
 0 & 0 & 1 
\end{pmatrix}, \quad& k \in \bar{E}_1 \cap \bar{E}_6, 
	\\
J_{18,7} = \begin{pmatrix}
 1 & - \overline{r(\bar{k})} e^{-t\Phi} & 0 \\
 r(k)e^{t\Phi} & 1 - |r(k)|^2 & 0 \\
 0 & 0 & 1 
\end{pmatrix}, & k \in \bar{E}_7 \cap \bar{E}_{18}, 
	\\
J_{1,7} 
% = \mathcal{B} \overline{J_{6,18}(\bar{k})} \mathcal{B}
= \begin{pmatrix}
 1 & \overline{h(\bar{k})} e^{-t\Phi} & 0 \\
0 & 1 & 0 \\
 0 & 0 & 1 
\end{pmatrix}, & k \in \bar{E}_1 \cap \bar{E}_7, \\
J_{6,18} =  \begin{pmatrix}
 1 & 0 & 0 \\
h(k) e^{t\Phi} & 1 & 0 \\
 0 & 0 & 1 
\end{pmatrix}, & k \in \bar{E}_6 \cap \bar{E}_{18}, \\
J_{4,3} =  \begin{pmatrix}
 1 & -\overline{\check{r}_1(\bar{k})}e^{-t\Phi} & 0 \\
\check{r}_1(k)e^{t\Phi} & 1 - |\check{r}_1(k)|^2 & 0 \\
 0 & 0 & 1 
\end{pmatrix}, & k \in \bar{E}_3 \cap \bar{E}_4, \\
J_{13,12} =  \begin{pmatrix}
 1 & -\overline{\check{r}(\bar{k})} e^{-t\Phi} & 0 \\
\check{r}(k) e^{t\Phi} & 1 - |\check{r}(k)|^2 & 0 \\
 0 & 0 & 1 
\end{pmatrix}, & k \in \bar{E}_{12} \cap \bar{E}_{13}, \\
J_{3,12} 
% = \mathcal{B} \overline{J_{4,13}(\bar{k})} \mathcal{B}
= \begin{pmatrix}
 1 & \overline{\check{h}(\bar{k})} e^{-t\Phi}  & 0 \\
0 & 1 & 0 \\
 0 & 0 & 1 
\end{pmatrix}, & k \in \bar{E}_3 \cap \bar{E}_{12}, \\
J_{4,13} = \begin{pmatrix}
 1 & 0 & 0 \\
\check{h}(k) e^{t\Phi} & 1 & 0 \\
 0 & 0 & 1 
\end{pmatrix}, & k \in \bar{E}_4 \cap \bar{E}_{13},
\end{cases}
\end{align}
where $\Phi := \Phi(\zeta, k)$ with $\zeta = y/t$ and
\begin{align}\label{Phizetakdef}
& \Phi(\zeta, k) = (l_2(k) - l_1(k))\zeta + (z_2(k) - z_1(k)). 
\end{align}
From (\ref{symmetriesb}) and (\ref{symmetriesc}), we infer the symmetries 
$$r(k^{-1}) = \overline{r(\bar{k})}, \qquad r_1(k^{-1}) = \overline{r_1(\bar{k})}, \qquad h(k^{-1}) = \overline{h(\bar{k})}.$$ 
In particular, $|r(k^{-1})| = |r(k)|$ for $k \in \bar{E}_7 \cap \bar{E}_{18}$ and $|r_1(k^{-1})| = |r_1(k)|$ for $k \in \R_+$. The functions $\check{r}$, $\check{r}_1$, and $\check{h}$ satisfy analogous symmetries.
The matrices $J_{1,2}$, $J_{5,6}$ and $J_{3,2}$, $J_{4,5}$ are not listed in (\ref{Jformula1}) because they can be recovered by symmetries from $J_{3,4}$ and $J_{1,6}$, respectively.

Proceeding as in Section 5 of \cite{L2013} we also obtain the identities 
\begin{align*}
& r_1(k) = r(k) + h(k), \qquad k \in \bar{E}_7 \cap \bar{E}_{18}, 
	\\
& \check{r}_1(k) = \check{r}(k) + \check{h}(k), \qquad k \in \bar{E}_{12} \cap \bar{E}_{13}.
\end{align*}
The first of these identities ensures that appropriate cyclic products of the relevant jump matrices equal the identity matrix at the intersection points where the sets $E_1, E_6, E_7, E_{18}$ meet. Similarly, the second identity ensures that appropriate cyclic products of the jump matrices equal the identity matrix where the sets $E_3, E_4, E_{12},E_{13}$ meet.

In a similar way, we find that the jump matrix $J$ for $k$ near $K_1$ is given by (see Figure \ref{EnearK1.pdf})
\begin{align}\label{Jformula2}
J = \begin{cases} 
J_{19, 25} 
= e^{y \hat{\mathcal{L}}  + t \hat{\mathcal{Z}}}   \begin{pmatrix}
1 & f_1(k) & 0 \\
0 & 1 & 0 \\
0 & 0 & 1 
\end{pmatrix},\quad &  k \in \bar{E}_{19} \cap \bar{E}_{25}, 
	\\
J_{19, 26} =  e^{y \hat{\mathcal{L}}  + t \hat{\mathcal{Z}}} \begin{pmatrix}
1 & 0 & 0 \\
0 & 1 & f_2(k) \\
0 & 0 & 1 
\end{pmatrix}, & k \in \bar{E}_{19} \cap \bar{E}_{26},
\end{cases}
\end{align}
where the functions $f_1$ and $f_2$ are bounded and continuous on the given subcontours.

We finally need the form of the jump matrix $J$  for $k$ on the circles where the $E_n$'s and $E_{n+18}$'s meet. 

\begin{lemma}\label{Jlemma}
With the contour oriented as in Figure \ref{Econtour.pdf}, we have
\begin{align}\label{Jformula3}
J = \begin{cases} 
J_{1, 19} = 
e^{y \hat{\mathcal{L}}  + t \hat{\mathcal{Z}}}\begin{pmatrix}
1 & g_1(k) & g_2(k) \\ 0 & 1 & g_3(k) \\ 0 & 0 & 1
\end{pmatrix}, \quad& k \in \bar{E}_1 \cap \bar{E}_{19}, 
	\\
J_{7, 25} = e^{y \hat{\mathcal{L}}  + t \hat{\mathcal{Z}}} \begin{pmatrix}
1 & 0 & g_4(k) \\ 0 & 1 & g_5(k) \\ 0 & 0 & 1
\end{pmatrix}, & k \in \bar{E}_7 \cap \bar{E}_{25}, 
	\\
J_{8, 26} 
= e^{y \hat{\mathcal{L}}  + t \hat{\mathcal{Z}}} \begin{pmatrix}
1 & g_6(k) & g_7(k) \\ 0 & 1 & 0 \\ 0 & 0 & 1
\end{pmatrix}, & k \in \bar{E}_8 \cap \bar{E}_{26}, 
\end{cases}
\end{align}
where the functions $\{g_j(k)\}_1^7$ are bounded and continuous on the given subcontours. 
\end{lemma}
\begin{proof}
In view of (\ref{Jmndef}), it is enough to show that 
\begin{align}\nonumber
&C_1(k) = \begin{pmatrix} 1 & * & * \\ 0 & 1 & * \\ 0 & 0 & 1 \end{pmatrix}, \quad k \in \bar{D}_1,
	\\\nonumber
& C_7(k) = \begin{pmatrix} 1 & 0 & * \\ 0 & 1 & * \\ 0 & 0 & 1 \end{pmatrix}, \quad k \in \bar{D}_7, 
	\\ \label{C178}
&C_8(k) = \begin{pmatrix} 1 & * & * \\ 0 & 1 & 0 \\ 0 & 0 & 1 \end{pmatrix}, \quad k \in \bar{D}_8,
\end{align}
where $*$ denotes an entry which is bounded and continuous except for possible singularities at the points $\varkappa_j$, $K_j$, $0$, and $\infty$. 
For $n = 1,7,8$, the matrices $(\gamma^n)_{ij} = \gamma_{ij}^n$ are given by (see (\ref{gammaijnudef}))
\begin{align*}
& \gamma^1 = \begin{pmatrix} 
\gamma_3 &\gamma_2 &\gamma_2 \\
\gamma_3 &\gamma_3 &\gamma_2 \\
\gamma_3 &\gamma_3 &\gamma_3 
 \end{pmatrix}, \quad
\gamma^7 =   \begin{pmatrix} 
\gamma_3 &\gamma_1 &\gamma_2 \\
\gamma_3 &\gamma_3 &\gamma_2 \\
\gamma_3 &\gamma_3 &\gamma_3 
 \end{pmatrix}, \quad
\gamma^8 =   \begin{pmatrix} 
\gamma_3 &\gamma_2 &\gamma_2 \\
\gamma_3 &\gamma_3 &\gamma_1 \\
\gamma_3 &\gamma_3 &\gamma_3
 \end{pmatrix}.
\end{align*}
Hence evaluation of (\ref{Cndef}) as $(y,t) \to (\infty,0)$ yields
\begin{align*}
C_n(k) = \lim_{y \to \infty} e^{-\nu_0\mathcal{L}} e^{-\mathcal{L} x} \begin{pmatrix} 1 & * & * \\ 0 & 1 & * \\ 0 & 0 & 1 \end{pmatrix}^{-1}\begin{pmatrix} 1 & * & * \\ 0 & 1 & * \\ 0 & 0 & 1 \end{pmatrix}  e^{\mathcal{L} y}
= \begin{pmatrix} 1 & * & * \\ 0 & 1 & * \\ 0 & 0 & 1 \end{pmatrix}, 
\end{align*}
for $k \in D_n$ and $n = 1,7,8$. Moreover, thanks to the assumed decay of the Dirichlet and Neumann values as $t \to \infty$, the functions $\tilde{\Psi}_n$ and $\Psi_n$ are bounded as $(y,t) \to (0, \infty)$ for each $k \in D_n$. Consequently, using that $\re z_2 < \re z_1 < \re z_3$ in $D_7$ and $\re z_1 < \re z_3 < \re z_2$ in $D_8$, evaluation of (\ref{Cndef}) as $(y,t) \to (0,\infty)$ yields $(C_7(k))_{12} = 0$ for $k \in D_7$ and $(C_8(k))_{23} = 0$ for $k \in D_8$. This proves (\ref{C178}).
\end{proof}

\section{A nonlinear steepest descent theorem}\label{steepsec}
We prove a nonlinear steepest descent theorem suitable for determining the asymptotics of (\ref{DP}) in the similarity region.

For $r > 0$, let $X^r = X_1^r \cup \cdots \cup X_4^r$ denote the cross $X = X_1 \cup \cdots \cup X_4$ defined in (\ref{Xdef}) restricted to the disk of radius $r$ centered at the origin, i.e. $X^r = X \cap \{|z| < r\}$. The spaces $\dot{E}^p$ and $\dot{L}^p$ are defined in Appendix \ref{RHapp}.
Let $\mathcal{I} \subset \R$ be a (possibly infinite) interval. Let $\rho, \epsilon:\mathcal{I} \to (0,\infty)$ be bounded strictly positive functions. Let $k_0:\mathcal{I} \to [1/2, 1)$ be a function such that $k_0(\zeta) + \epsilon(\zeta) <1$ and $\epsilon(\zeta) < k_0(\zeta)/2$ for each $\zeta \in \mathcal{I}$. 
We henceforth drop the $\zeta$ dependence of these functions and write simply $\rho$, $\epsilon$, $k_0$ for $\rho(\zeta)$, $\epsilon(\zeta)$, $k_0(\zeta)$, respectively.

Before stating Theorem \ref{steepestdescentth} we list the necessary assumptions.

\begin{figure}
\begin{center}
\bigskip\bigskip
 \begin{overpic}[width=.7\textwidth]{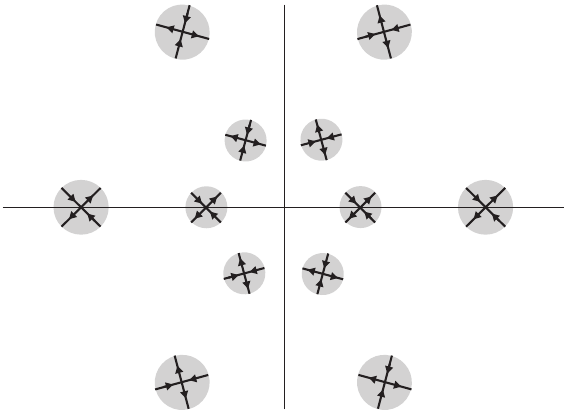}
      \put(101,36){\small $\re k$}
      \put(47,75){\small $\im k$}
      \put(62.6,32){\small $k_0$}
      \put(82.5,30){\small $k_0^{-1}$}
      \put(32.5,32){\small $-k_0$}
      \put(8.5,30){\small $-k_0^{-1}$}
   \end{overpic}
   \begin{figuretext}\label{VGammaX.pdf}
      The contour $\Gamma_X$ and the open set $\mathcal{V}$ (shaded).
             \end{figuretext}
   \end{center}
\end{figure}

\begin{assumptions}[Contour Assumptions]\label{contourassumptions}
%Let $\{p_n\}_1^{12}$ denote the set of `critical points', i.e.,
%$$\{p_n\}_1^{12} = \{\pm k_0, \pm k_0^{-1}, \pm \omega k_0, \pm \omega k_0^{-1}, \pm \omega^2 k_0, \pm \omega^2 k_0^{-1}\}.$$ 
Assume $\Gamma = \Gamma(\zeta)$ is a family of Carleson jump contours parametrized by $\zeta \in \mathcal{I}$ such that:
\begin{itemize}
\item[($\Gamma$1)] For each $\zeta \in \mathcal{I}$, $\Gamma$ contains the small crosses $\pm k_0 + X^{\epsilon}$ as a subset.

\item[($\Gamma$2)] For each $\zeta \in \mathcal{I}$, $\Gamma$ is invariant as a set under the maps 
\begin{align}\label{symmmaps}
 k \mapsto \omega k, \qquad k \mapsto 1/k.
\end{align} 
Moreover, the orientation of $\Gamma$ is such that if $k$ traverses $\Gamma$ in the positive direction, then $\omega k$ and $1/k$ also traverse $\Gamma$ in the positive direction.

\item[($\Gamma$3)]
Let $\mathcal{V}$ denote the union of the two disks $\{|k \pm k_0| < \epsilon\}$ and the sets obtained by letting the symmetries in (\ref{symmmaps}) act repeatedly on these disks, see Figure \ref{VGammaX.pdf}. 
Let $\hat{\Gamma} = \Gamma \cup \partial \mathcal{V}$ and assume that the boundary of each of the $12$ components of $\mathcal{V}$ is oriented counterclockwise. Then, after reversing the orientation on a subcontour if necessary, $\hat{\Gamma}$ is a Carleson jump contour for each $\zeta \in \mathcal{I}$.

\item[($\Gamma$4)] The contour remains a bounded distance away from the point $K_1 := e^{\frac{i\pi}{6}}$ for all $\zeta \in \mathcal{I}$:
\begin{align}\label{distK1}
\inf_{\zeta \in \mathcal{I}} \text{\upshape dist}(K_1, \hat{\Gamma}) > 0.
\end{align}
\end{itemize}
Assume also that the Cauchy singular operator $\mathcal{S}_{\hat{\Gamma}}$ defined by
\begin{align*}
(\mathcal{S}_{\hat{\Gamma}} h)(z) = \lim_{r \to 0} \frac{1}{\pi i} \int_{\hat{\Gamma} \setminus \{|z' - z| < r\}} \frac{h(z')}{z' - z} dz',
\end{align*}
is uniformly\footnote{For any fixed $\zeta \in \mathcal{I}$, $\mathcal{S}_{\hat{\Gamma}}$ is bounded on $L^2(\hat{\Gamma})$ as a consequence of ($\Gamma$3).} bounded on $L^2(\hat{\Gamma})$, i.e. 
\begin{align}\label{cauchysingularbound}
\sup_{\zeta \in \mathcal{I}} \|\mathcal{S}_{\hat{\Gamma}}\|_{\mathcal{B}(L^2(\hat{\Gamma}))} < \infty.
\end{align}
\end{assumptions}

We consider the following family of $L^3$-RH problems parametrized by the two parameters $\zeta \in \mathcal{I}$ and $t > 0$:
\begin{align}\label{RHm}
\begin{cases} m(\zeta, t, \cdot) \in I + \dot{E}^3(\hat{\C} \setminus \Gamma), \\
m_+(\zeta, t, k) = m_-(\zeta, t, k) v(\zeta, t, k) \quad \text{for a.e.} \ k \in \Gamma,
\end{cases} 
\end{align}
where the jump matrix $v$ is assumed to fulfill the following conditions.

\begin{assumptions}[Jump Assumptions]\label{jumpassumptions}
Assume the jump matrix $v$ obeys the symmetries (\ref{symmetriesa}) and (\ref{symmetriesb}) and satisfies
\begin{align}\label{winL1L2Linf}
w(\zeta, t,\cdot) := v(\zeta, t,\cdot) - I \in \dot{L}^1(\Gamma) \cap L^\infty(\Gamma), \qquad \zeta \in \mathcal{I}, \ t > 0.
\end{align}
Let $\tau := t \rho^2$. Let $\Gamma_X$ denote the union of the two small crosses $\pm k_0 + X^{\epsilon}$ and the sets obtained by letting the symmetries in (\ref{symmmaps}) act repeatedly on these crosses. Let $\Gamma' = \Gamma \setminus \Gamma_X$ and suppose 
\begin{subequations}\label{wL12infty}
\begin{align}\label{wL12inftya}
&\|w(\zeta,t,\cdot)\|_{\dot{L}^p(\Gamma')} = O(\epsilon^{\frac{1}{p}} \tau^{-1}), \qquad \tau \to \infty, \ \zeta \in \mathcal{I}, \ p \in \{1, \tfrac{3}{2}, 3\},
	\\ \label{wL12inftyb}
&\|w(\zeta,t,\cdot)\|_{L^\infty(\Gamma')} = O(\tau^{-1}),  \qquad \tau \to \infty, \ \zeta \in \mathcal{I},
\end{align}
\end{subequations}
uniformly with respect to $\zeta \in \mathcal{I}$. 
Moreover, let $\mathcal{C} = \diag(1, -1, 1)$ and suppose that the normalized jump matrices 
\begin{align}\label{vjdef}
\begin{cases}
 v_0(\zeta,t,z) = \mathcal{C}v\Big(\zeta, t, k_0 - \frac{\epsilon z}{\rho}\Big)\mathcal{C}, 
	\vspace{.1cm}\\
\check{v}_0(\zeta,t,z) = \mathcal{C}\overline{v\Big(\zeta, t, \overline{-k_0 + \frac{\epsilon z}{\rho}}\Big)} \mathcal{C}, 
\end{cases}
\qquad z \in X^{\rho}, \quad \zeta \in \mathcal{I},
\end{align}
have the form
\begin{subequations}\label{smallcrossjump}
\begin{align}
v_0(\zeta, t, z) = \begin{cases}
\begin{pmatrix} 1 & 0	 & 0 \\
  R_1(\zeta, t, z)z^{-2i\nu(\zeta)} e^{t\phi(\zeta, z)}	& 1 & 0 \\ 0 & 0 & 1 \end{pmatrix}, &  z \in X_1^{\rho}, \\
\begin{pmatrix} 1 & -R_2(\zeta, t, z)z^{2i\nu(\zeta)}e^{-t\phi(\zeta, z)} & 0	\\
0 & 1 & 0 \\ 0 & 0 & 1  \end{pmatrix}, &  z \in X_2^{\rho}, \\
\begin{pmatrix} 1 &0 & 0  \\
 -R_3(\zeta, t, z)z^{-2i\nu(\zeta)} e^{t\phi(\zeta, z)}	& 1 & 0 \\ 0 & 0 & 1 \end{pmatrix}, &  z \in X_3^{\rho},
 	\\
 \begin{pmatrix} 1  & R_4(\zeta, t, z)z^{2i\nu(\zeta)}e^{-t\phi(\zeta, z)} & 0 \\ 0	& 1 & 0 \\ 0 & 0 & 1 \end{pmatrix}, & z \in X_4^{\rho}, 
\end{cases}
\end{align}
and
\begin{align}
\check{v}_0(\zeta, t, z) =  \begin{cases}
\begin{pmatrix} 1 & 0	 & 0 \\
  \check{R}_1(\zeta, t, z)z^{-2i\check{\nu}(\zeta)} e^{t\phi(\zeta, z)}	& 1 & 0 \\ 0 & 0 & 1 \end{pmatrix}, &  z \in X_1^{\rho}, \\
\begin{pmatrix} 1 & -\check{R}_2(\zeta, t, z)z^{2i\check{\nu}(\zeta)}e^{-t\phi(\zeta, z)} & 0	\\
0 & 1 & 0 \\ 0 & 0 & 1  \end{pmatrix}, &  z \in X_2^{\rho}, \\
\begin{pmatrix} 1 &0 & 0  \\
 -\check{R}_3(\zeta, t, z)z^{-2i\check{\nu}(\zeta)} e^{t\phi(\zeta, z)}	& 1 & 0 \\ 0 & 0 & 1 \end{pmatrix}, &  z \in X_3^{\rho},
 	\\
 \begin{pmatrix} 1	& \check{R}_4(\zeta, t, z)z^{2i\check{\nu}(\zeta)}e^{-t\phi(\zeta, z)} & 0 \\ 0	& 1 & 0 \\ 0 & 0 & 1 \end{pmatrix}, & z \in X_4^{\rho}, 
\end{cases}
\end{align}
\end{subequations}
where:
\begin{itemize}
\item The phase $\phi(\zeta, z)$ is a smooth function of $(\zeta, z) \in \mathcal{I} \times \C$ 
such that 
\begin{align}\label{phiassumptions}
\phi(\zeta, 0) \in i\R, \qquad \frac{\partial \phi}{\partial z}(\zeta, 0) = 0, \qquad \frac{\partial^2 \phi}{\partial z^2}(\zeta, 0) = i, \qquad \zeta \in \mathcal{I},
\end{align}
and
\begin{subequations}\label{rephiestimates}
\begin{align} \label{rephiestimatea}
 & \re \phi(\zeta,z) \leq -\frac{|z|^2}{4}, &&  z \in X_1^{\rho} \cup X_3^{\rho}, \ \zeta \in \mathcal{I},
  	\\ \label{rephiestimateb}
 & \re \phi(\zeta,z) \geq \frac{|z|^2}{4}, &&  z \in X_2^{\rho} \cup X_4^{\rho},  \ \zeta \in \mathcal{I},
  	\\ \label{Phiz3estimate}
 & \biggl|\phi(\zeta, z) - \phi(\zeta,0) - \frac{iz^2}{2}\biggr| \leq C \frac{|z|^3}{\rho}, && z \in X^{\rho}, \ \zeta \in \mathcal{I},
\end{align}
\end{subequations}
where $C > 0$ is a constant.

\item There exist smooth functions $q, \check{q}:\mathcal{I} \to \C$ and constants $(\alpha, L) \in [\frac{1}{3},1) \times (0,\infty)$ such that 
$$\sup_{\zeta \in \mathcal{I}} |q(\zeta)| < 1, \qquad \sup_{\zeta \in \mathcal{I}} |\check{q}(\zeta)| < 1,$$
and
\begin{subequations}\label{Lipschitzconditions}
\begin{align} 
& \begin{cases}
   |R_1(\zeta, t, z) - q(\zeta)| \leq  L  \bigl| \frac{z}{\rho}\bigr|^\alpha e^{\frac{t|z|^2}{6}}, & z \in X_1^{\rho},  \\
 |R_2(\zeta, t, z) - \frac{\overline{q(\zeta)}}{1 - |q(\zeta)|^2} | \leq  L \bigl| \frac{z}{\rho}\bigr|^\alpha e^{\frac{t|z|^2}{6}}, \qquad& z \in X_2^{\rho}, \vspace{.1cm}\\ 
  |R_3(\zeta, t, z) - \frac{q(\zeta)}{1 - |q(\zeta)|^2}\big| \leq  L \bigl| \frac{z}{\rho}\bigr|^\alpha e^{\frac{t|z|^2}{6}}, & z \in X_3^{\rho}, 
  	\\
  |R_4(\zeta, t, z) - \overline{q(\zeta)}| \leq L  \bigl| \frac{z}{\rho}\bigr|^\alpha e^{\frac{t|z|^2}{6}}, & z \in X_4^{\rho}, 
\end{cases} 
	\\
& \begin{cases}
   |\check{R}_1(\zeta, t, z) - \check{q}(\zeta)| \leq  L  \bigl| \frac{z}{\rho}\bigr|^\alpha e^{\frac{t|z|^2}{6}}, & z \in X_1^{\rho},  \\
 |\check{R}_2(\zeta, t, z) - \frac{\overline{\check{q}(\zeta)}}{1 - |\check{q}(\zeta)|^2} | \leq  L \bigl| \frac{z}{\rho}\bigr|^\alpha e^{\frac{t|z|^2}{6}}, \qquad& z \in X_2^{\rho}, \vspace{.1cm}\\ 
  |\check{R}_3(\zeta, t, z) - \frac{\check{q}(\zeta)}{1 - |\check{q}(\zeta)|^2}\big| \leq  L \bigl| \frac{z}{\rho}\bigr|^\alpha e^{\frac{t|z|^2}{6}}, & z \in X_3^{\rho}, 
  	\\
  |\check{R}_4(\zeta, t, z) - \overline{\check{q}(\zeta)}| \leq L  \bigl| \frac{z}{\rho}\bigr|^\alpha e^{\frac{t|z|^2}{6}}, & z \in X_4^{\rho}, 
\end{cases} 
\end{align}
\end{subequations}
for $\zeta \in \mathcal{I}$ and $t > 0$.
\item The functions $\nu(\zeta)$ and $\check{\nu}(\zeta)$ are defined by 
\begin{align}\label{nuchecknudef}
\nu(\zeta) = -\frac{1}{2\pi} \log(1 - |q(\zeta)|^2), \qquad 
\check{\nu}(\zeta) = -\frac{1}{2\pi} \log(1 - |\check{q}(\zeta)|^2).
\end{align}
\end{itemize}
\end{assumptions}

\begin{theo}[Nonlinear steepest descent]\label{steepestdescentth}
If the contour $\Gamma$ and the jump matrix $v(\zeta, t, k)$ satisfy Assumptions \ref{contourassumptions} and Assumptions \ref{jumpassumptions}, then the $L^3$-RH problem (\ref{RHm}) has a unique solution for all sufficiently large $\tau = t\rho^2$ and this solution satisfies
\begin{align}\nonumber
& (1, 1,1) m(\zeta,t,K_1)
=  (1,1,1) 
	\\ \label{limlm12}
& + \frac{2\epsilon}{k_0 \sqrt{\tau}} \re\big(\mathcal{F}_1 \beta - \bar{\mathcal{F}}_2\check{\beta}, 
\mathcal{F}_3 \beta - \bar{\mathcal{F}}_3\check{\beta}, \mathcal{F}_2 \beta - \bar{\mathcal{F}}_1\check{\beta}\big) + O\bigl(\epsilon\tau^{-\frac{1+\alpha}{2}} \bigr)
\end{align}
for $\zeta \in \mathcal{I}$  as $\tau \to \infty$, where the error term is uniform with respect to $\zeta \in \mathcal{I}$, the functions $\mathcal{F}_j = \mathcal{F}_j(\zeta)$, $j = 1,2,3,$ are defined by
\begin{align}\nonumber
 & \mathcal{F}_1(\zeta) = \frac{1 - k_0^2\omega}{(i + k_0)(1- k_0 K_1)},  \qquad
  \mathcal{F}_2(\zeta) = \frac{1 - k_0^2\omega}{(i - k_0)(1 + k_0 K_1)}, 
  	\\\label{Xjdef}
&  \mathcal{F}_3(\zeta) = -\frac{i(1 - k_0^2\omega)}{1 + k_0^2 \omega},
\end{align}
and the functions $\beta = \beta(\zeta, t)$ and $\check{\beta} = \check{\beta}(\zeta, t)$ are defined by
\begin{subequations}\label{betadef}
\begin{align}
&  \beta(\zeta, t) = \sqrt{\nu(\zeta)} e^{i\left(\frac{\pi}{4} - \arg q(\zeta) + \arg \Gamma(i\nu(\zeta))\right)} e^{-t\phi(\zeta, 0)} t^{-i\nu(\zeta)},
	\\
& \check{\beta}(\zeta, t) = \sqrt{\check{\nu}(\zeta)} e^{i\left(\frac{\pi}{4} - \arg \check{q}(\zeta) + \arg \Gamma(i\check{\nu}(\zeta))\right)} e^{-t\phi(\zeta, 0)} t^{-i\check{\nu}(\zeta)}. 
\end{align}  
\end{subequations}
\end{theo}
\begin{proof}
Since $\det v = 1$ and we are considering an $L^3$-RH problem for a $3 \times 3$-matrix valued function, uniqueness follows from Lemma \ref{uniquelemma}.

Let $m^X$ be the solution of Theorem \ref{crossth} and let
\begin{align*}
& D(\zeta, t) = \diag\big(e^{-\frac{t\phi(\zeta, 0)}{2}}t^{-\frac{i\nu(\zeta)}{2}}, e^{\frac{t\phi(\zeta, 0)}{2}}t^{\frac{i\nu(\zeta)}{2}}, 1\big), 
	\\
& \check{D}(\zeta, t) = \diag\big(e^{-\frac{t\phi(\zeta, 0)}{2}}t^{-\frac{i\check{\nu}(\zeta)}{2}}, e^{\frac{t\phi(\zeta, 0)}{2}}t^{\frac{i\check{\nu}(\zeta)}{2}}, 1\big).
\end{align*}
Define $m_0(\zeta,t,k)$ in neighborhoods of $k = k_0$ and $k = -k_0$ by
\begin{align*}
& m_0(\zeta, t, k) = \begin{cases}
\mathcal{C} D(\zeta, t) m^X\Bigl(q(\zeta),  -\frac{\sqrt{\tau}}{\epsilon} (k-k_0)\Bigr) D(\zeta, t)^{-1}\mathcal{C}, & |k - k_0| \leq \epsilon, 
	\vspace{.1cm} \\
\mathcal{C}\check{D}(\zeta, t)^{-1} \overline{m^X\Bigl(\check{q}(\zeta), \frac{\sqrt{\tau}}{\epsilon} \overline{(k+k_0)}\Bigr)} \check{D}(\zeta, t)\mathcal{C}, & |k + k_0| \leq \epsilon,	
\end{cases}
\end{align*}
and extend it to all of $\mathcal{V}$ in such a way that $m_0$ obeys the symmetries (\ref{symmetriesa}) and (\ref{symmetriesb}).

Lemma \ref{deformationlemma} implies that $m$ satisfies the $L^3$-RH problem (\ref{RHm}) if and only if the function $\hat{m}(\zeta, t, k)$ defined by
\begin{align*}
\hat{m}(\zeta,t,k) = \begin{cases}
m(\zeta, t, k)m_0(\zeta,t,k)^{-1}, & k \in \mathcal{V},\\
m(\zeta, t, k), & \text{otherwise},
\end{cases}
\end{align*}
satisfies the $L^3$-RH problem
\begin{align}\label{RHmhat}
\begin{cases}
\hat{m}(\zeta,t,\cdot) \in I + \dot{E}^3(\hat{\C} \setminus \hat{\Gamma}), \\
\hat{m}_+(\zeta,t,k) = \hat{m}_-(\zeta, t, k) \hat{v}(x, t, k) \quad \text{for a.e.} \ k \in \hat{\Gamma}, 
\end{cases}
\end{align}
where the jump matrix $\hat{v}$ is given by
\begin{align*}
\hat{v}(\zeta, t, k) 
=  \begin{cases}
 m_{0-}(\zeta, t, k) v(\zeta, t, k) m_{0+}(\zeta,t,k)^{-1}, & k \in \mathcal{V}, \\
m_0(\zeta, t, k)^{-1}, & k \in \partial \mathcal{V}, \\
v(\zeta, t, k),  & \text{otherwise}.
\end{cases}
\end{align*}
By construction, $\hat{w} = \hat{v} - I$ satisfies the symmetries (\ref{symmetriesa}) and (\ref{symmetriesb}).

The proofs of the following five claims can be found in Appendix \ref{proofapp}.

\begin{claim}\label{claim1}
The function $\hat{w} = \hat{v} - I$ satisfies
\begin{align}\label{hatwestimate}
\hat{w}(\zeta,t,k) = O\bigl(\tau^{-\frac{\alpha}{2}} e^{-\frac{\tau}{24\epsilon^2}|k \mp k_0|^2}\bigr), \qquad \tau \to \infty, \  \zeta \in \mathcal{I}, \  k \in \pm k_0 + X^\epsilon,
\end{align}
where the error term is uniform with respect to $(\zeta, k)$ in the given ranges.
\end{claim}

\begin{claim}\label{claim2}
We have
\begin{subequations}
\begin{align} \label{hatwestimatea}
&\|\hat{w}(\zeta, t, \cdot)\|_{\dot{L}^3(\hat{\Gamma})} = O(\epsilon^{\frac{1}{3}} \tau^{-\frac{\alpha}{2}}), \qquad \tau \to \infty, 
\  \zeta \in \mathcal{I}, 
	\\ \label{hatwestimateb}
&\|\hat{w}(\zeta, t, \cdot)\|_{L^\infty(\hat{\Gamma})} = O(\tau^{-\frac{\alpha}{2}}), \qquad \tau \to \infty, 
\  \zeta \in \mathcal{I},
\end{align}
\end{subequations}
and, for any $p \in [1, \infty)$,
\begin{align}\label{hatwestimatec}
& \|\hat{w}(\zeta, t, \cdot)\|_{\dot{L}^p(\pm k_0 + X^\epsilon)} = O(\epsilon^{\frac{1}{p}} \tau^{-\frac{1}{2p} - \frac{\alpha}{2}}), \qquad \tau \to \infty, \  \zeta \in \mathcal{I},
	\\ \label{star}
& \|m_0(\zeta, t, k)^{-1} - I \|_{L^p(|k - k_0| = \epsilon)} = O(\epsilon^{\frac{1}{p}} \tau^{-\frac{1}{2}}), \qquad \tau \to \infty, \  \zeta \in \mathcal{I},
\end{align}
where the error terms are uniform with respect to $\zeta \in \mathcal{I}$.
\end{claim}

Let $\hat{\mathcal{C}}$ denote the Cauchy operator associated with $\hat{\Gamma}$:
$$(\hat{\mathcal{C}} f)(z) = \frac{1}{2\pi i} \int_{\hat{\Gamma}} \frac{f(s)}{s - z} ds, \qquad z \in \C \setminus \hat{\Gamma}.$$
The operator $\hat{\mathcal{C}}_{\hat{w}}: \dot{L}^3(\hat{\Gamma}) + L^\infty(\hat{\Gamma}) \to \dot{L}^p(\hat{\Gamma})$ is defined by $\hat{\mathcal{C}}_{\hat{w}}(h) = \hat{\mathcal{C}}_-(h \hat{w})$, where $\hat{\mathcal{C}}_-f$ denotes the nontangential boundary value of $\hat{\mathcal{C}}f$ from the right side of $\hat{\Gamma}$.

\begin{claim}\label{claim3}
There exists a $T > 0$ such that $I - \hat{\mathcal{C}}_{\hat{w}(\zeta, t, \cdot)} \in \mathcal{B}(\dot{L}^3(\hat{\Gamma}))$ is invertible for all $(\zeta,t) \in \mathcal{I} \times (0, \infty)$ with $\tau > T$.
\end{claim}

In view of Claim \ref{claim3}, we may define the $3 \times 3$-matrix valued function $\hat{\mu}(\zeta, t, z)$ whenever $\tau > T$ by
\begin{align}\label{hatmudef}
\hat{\mu} = I + (I - \hat{\mathcal{C}}_{\hat{w}})^{-1}\hat{\mathcal{C}}_{\hat{w}}I  \in I + \dot{L}^3(\hat{\Gamma}).
\end{align}

\begin{claim}\label{claim4}
The function $\hat{\mu}(\zeta, t, k)$ satisfies
\begin{align}\label{muhatestimate}
\|\hat{\mu}(\zeta,t,\cdot) - I\|_{\dot{L}^3(\hat{\Gamma})} = O\big(\epsilon^{\frac{1}{3}}\tau^{-\frac{\alpha}{2}}\big), \qquad \tau \to \infty, \  \zeta \in \mathcal{I},
\end{align}
where the error term is uniform with respect to $\zeta \in \mathcal{I}$.
\end{claim}

\begin{claim}\label{claim5}
There exists a unique solution $\hat{m} \in I + \dot{E}^3(\hat{\C} \setminus \hat{\Gamma})$ of the $L^3$-RH problem (\ref{RHmhat}) whenever $\tau > T$. This solution is given by
\begin{align}\label{hatmrepresentation}
\hat{m}(\zeta, t, k) = I + \hat{\mathcal{C}}_{\hat{w}}\hat{\mu} = I + \frac{1}{2\pi i}\int_{\hat{\Gamma}} \hat{\mu}(\zeta, t, s) \hat{w}(\zeta, t, s) \frac{ds}{s - k}.
\end{align}
\end{claim}

Using the above claims we can complete the proof of Theorem \ref{steepestdescentth} as follows.

Let $C(\zeta)$ denote the union of the two circles $|k - k_0|=\epsilon$ and $|k+  k_0|=\epsilon$ oriented counterclockwise. Let $C(\zeta)^{-1}$ denote the image of $C(\zeta)$ under the map  $k \mapsto k^{-1}$. The symmetry properties of $v$ imply that  $\mathcal{A} \hat{m}(\zeta, t, \omega k) \mathcal{A}^{-1} \in I + \dot{E}^3(\hat{\C} \setminus \hat{\Gamma})$ and $\hat{m}(\zeta, t, k)$ both satisfy the $L^3$-RH problem (\ref{RHmhat}); by uniqueness they are equal, i.e.,
$$\hat{m}(\zeta, t, k) = \mathcal{A} \hat{m}(\zeta, t, \omega k) \mathcal{A}^{-1}, \qquad k \in \hat{\C} \setminus \hat{\Gamma}.$$
Using this symmetry in (\ref{hatmrepresentation}), we obtain
\begin{align}\nonumber
 m(\zeta,t,K_1) = \hat{m}(\zeta,t,K_1) 
 = &\; I + \frac{1}{2\pi i} \sum_{n=0}^2 \mathcal{A}^{-n} \big[F_n(\zeta, t) + G_n(\zeta, t) \big]\mathcal{A}^n
	\\ \label{matK1}
& + \frac{1}{2\pi i}\int_{\Gamma} \hat{\mu}(\zeta,t,k) \hat{w}(\zeta,t,k) \frac{dk}{k - K_1},
\end{align}
where
\begin{align*}
& F_n(\zeta, t) = \int_{C(\zeta)} \frac{\hat{\mu}(\zeta,t,k) \hat{w}(\zeta,t,k) dk}{k - \omega^{-n}  K_1}, 
	\\
& G_n(\zeta, t) = \int_{C(\zeta)^{-1}} \frac{\hat{\mu}(\zeta,t,k) \hat{w}(\zeta,t,k) dk}{k - \omega^{-n}  K_1}.
\end{align*}
By (\ref{mcasymptotics}), we have, as $\tau \to \infty$ with $\zeta \in \mathcal{I}$, 
\begin{subequations}\label{mjinvasymptotics}
\begin{align} \nonumber
 m_0(\zeta, t, k)^{-1} =&\;   \mathcal{C}D(\zeta, t) m^X\biggl(q(\zeta), - \frac{\sqrt{\tau}}{\epsilon}(k - k_0)\biggr)^{-1} D(\zeta, t)^{-1}\mathcal{C}
  	\\ 
= &\; I - \frac{B(\zeta, t)}{\sqrt{\tau}(k - k_0)} + O(\tau^{-1}), \qquad |k - k_0| = \epsilon,
	\\ \nonumber
  m_0(\zeta, t, k)^{-1} =&\; \mathcal{C}\check{D}(\zeta, t)^{-1} \overline{m^X\biggl(\check{q}(\zeta), \frac{\sqrt{\tau}}{\epsilon}\overline{(k + k_0)}\biggr)^{-1}} \check{D}(\zeta, t)\mathcal{C}
  	\\ 
= &\; I + \frac{\overline{\check{B}(\zeta, t)}}{\sqrt{\tau}(k +  k_0)} + O(\tau^{-1}), \qquad \  |k + k_0| = \epsilon,
\end{align}
\end{subequations}
where $B(\zeta, t)$ and $\check{B}(\zeta, t)$ are defined by
\begin{align*}
& B(\zeta, t) = i \epsilon \begin{pmatrix} 0 & -\beta(\zeta, t) & 0 \\ 
\overline{\beta(\zeta, t)} & 0 & 0 \\
0 & 0 & 0 \end{pmatrix},
	\\ 
& \check{B}(\zeta, t) = i \epsilon \begin{pmatrix} 0 & -\check{\beta}(\zeta, t) & 0 \\ 
 \overline{\check{\beta}(\zeta, t)} & 0 & 0 \\
0 & 0 & 0 \end{pmatrix},
\end{align*}
with 
$$\beta(\zeta, t) = \beta^X(q(\zeta))e^{-t\phi(\zeta,0)}t^{-i\nu(\zeta)}, \quad
\check{\beta}(\zeta, t) =  \beta^X(\check{q}(\zeta))e^{-t\phi(\zeta,0)}t^{-i\check{\nu}(\zeta)}.$$

Using (\ref{star}), (\ref{muhatestimate}), and (\ref{mjinvasymptotics}) we find
\begin{align}\nonumber
F_n(\zeta, t)
= &\;  \int_{C(\zeta)} \frac{\hat{\mu}(\zeta,t,k) (m_0(\zeta,t,k)^{-1} - I) dk}{k - \omega^{-n}  K_1}
	\\ \nonumber
= &\; \int_{C(\zeta)} \frac{(m_0(\zeta,t,k)^{-1} - I) dk}{k - \omega^{-n}  K_1}	
	\\\nonumber
& + \int_{C(\zeta)} \frac{(\hat{\mu}(\zeta,t,k) - I) (m_0(\zeta,t,k)^{-1} - I) dk}{k - \omega^{-n}  K_1}
	\\\nonumber
= & - \int_{|k - k_0| = \epsilon} \frac{B(\zeta, t)}{\sqrt{\tau}(k - k_0)} \frac{ dk}{k - \omega^{-n}  K_1}
	\\\nonumber
& + \int_{|k + k_0| = \epsilon} \frac{\overline{\check{B}(\zeta, t)}}{\sqrt{\tau}(k + k_0)} \frac{ dk}{k - \omega^{-n}  K_1}
	\\\nonumber
& + O(\epsilon\tau^{-1} )
+ O\bigl(\|\hat{\mu} - I\|_{\dot{L}^3(\hat{\Gamma})} \|m_0^{-1} - I\|_{\dot{L}^{3/2}(C(\zeta))} \bigr)
	\\\label{Fn}
= & - \frac{2\pi i}{\sqrt{\tau}}\bigg(\frac{B(\zeta, t)}{k_0 - \omega^{-n}  K_1}
+ \frac{\overline{\check{B}(\zeta, t)}}{k_0 + \omega^{-n}  K_1}\bigg)
 + O(\epsilon\tau^{-\frac{1+\alpha}{2}})
\end{align}
uniformly with respect to $\zeta \in \mathcal{I}$ as $\tau \to \infty$. 

In order to compute the contribution from $G_n$ we note that (\ref{mjinvasymptotics}) implies
\begin{align*} \nonumber
m_0(\zeta, t, k)^{-1} & = 
\mathcal{B}m_0(\zeta, t, k^{-1})^{-1}  \mathcal{B}
 	\\
& = 
\begin{cases}
I - \frac{ \mathcal{B}B(\zeta, t) \mathcal{B}}{\sqrt{\tau}(k^{-1} - k_0)} + O(\tau^{-1}), & |k^{-1} - k_0| = \epsilon, \\
I + \frac{ \mathcal{B} \overline{\check{B}(\zeta, t)} \mathcal{B}}{\sqrt{\tau}(k^{-1} +  k_0)} + O(\tau^{-1}),
&  |k^{-1} + k_0| = \epsilon.
\end{cases}
\end{align*}
Hence, proceeding as in (\ref{Fn}), we find
\begin{align*}\nonumber
 G_n(\zeta, t)
= &\; \int_{C(\zeta)^{-1}} \frac{\hat{\mu}(\zeta,t,k) (m_0(\zeta,t,k)^{-1} - I) dk}{k - \omega^{-n}  K_1}
	\\ \nonumber
= &\; \int_{C(\zeta)^{-1}} \frac{(m_0(\zeta,t,k)^{-1} - I) dk}{k - \omega^{-n}  K_1}	
	\\\nonumber
& + \int_{C(\zeta)^{-1}} \frac{(\hat{\mu}(\zeta,t,k) - I) (m_0(\zeta,t,k)^{-1} - I) dk}{k - \omega^{-n}  K_1}.
\end{align*}
That is,
\begin{align}\nonumber
 G_n(\zeta, t) = & - \int_{|k^{-1} - k_0| = \epsilon} \frac{\mathcal{B}B(\zeta, t)\mathcal{B}}{\sqrt{\tau}(k^{-1} - k_0)} \frac{ dk}{k - \omega^{-n}  K_1}
	\\\nonumber
& + \int_{|k^{-1} + k_0| = \epsilon} \frac{\mathcal{B}\overline{\check{B}(\zeta, t)}\mathcal{B}}{\sqrt{\tau}(k^{-1} + k_0)} \frac{ dk}{k - \omega^{-n}  K_1}
	\\\nonumber
& + O(\epsilon\tau^{-1}) + O\bigl(\|\hat{\mu} - I\|_{\dot{L}^3(\hat{\Gamma})} \|m_0^{-1} - I\|_{\dot{L}^{3/2}(C(\zeta)^{-1})} \bigr)
	\\\label{Gn}
= &\; \frac{2\pi i}{\sqrt{\tau}}\bigg(\frac{\mathcal{B} B(\zeta, t)\mathcal{B}}{k_0(1 - k_0 \omega^{-n}  K_1)}
 + \frac{\mathcal{B}\overline{\check{B}(\zeta, t)}\mathcal{B}}{k_0(1 + k_0 \omega^{-n}  K_1)}\bigg)
 + O(\epsilon\tau^{-\frac{1+\alpha}{2}})
\end{align}
uniformly with respect to $\zeta \in \mathcal{I}$. 
On the other hand, using (\ref{distK1}),
\begin{align}\nonumber
\biggl|\int_{\Gamma} \frac{\hat{\mu}(\zeta,t,k) \hat{w}(\zeta,t,k)dk}{k - K_1}\biggr|
 & = \biggl|\int_{\Gamma} \frac{(\hat{\mu} -I) \hat{w}dk}{k - K_1} + \int_{\Gamma} \frac{\hat{w} dk}{k - K_1}\biggr|
	\\ \nonumber
& \leq C \|\hat{\mu} - I\|_{\dot{L}^3(\Gamma)}  \|\hat{w}\|_{\dot{L}^{\frac{3}{2}}(\Gamma)} + C \|\hat{w}\|_{\dot{L}^1(\Gamma)}.
\end{align}
The $\dot{L}^1$-norm of $\hat{w}$  is $O(\epsilon \tau^{-1})$ on $\Gamma'$ by (\ref{wL12inftya}) and is $O(\epsilon\tau^{-\frac{1+\alpha}{2}})$ on $\{\pm k_0 + X^\epsilon\}$ by (\ref{hatwestimatec}). Hence $\|\hat{w}\|_{\dot{L}^1(\Gamma)} = O(\epsilon \tau^{-\frac{1+\alpha}{2}})$.
Similarly, 
$$\|\hat{w}\|_{\dot{L}^{\frac{3}{2}}(\Gamma)} = O\big(\epsilon^{\frac{2}{3}} \tau^{-1} + \epsilon^{\frac{2}{3}} \tau^{-\frac{1}{3} -\frac{\alpha}{2}}\big)$$ 
by (\ref{wL12inftya}) and (\ref{hatwestimatec}). Since $\|\hat{\mu} - I\|_{\dot{L}^3(\Gamma)} = O(\epsilon^{1/3}\tau^{-\frac{\alpha}{2}})$ by (\ref{muhatestimate}) and $1/3 \leq \alpha < 1$, we infer that
\begin{align}\label{intSigmahatmuhatw}
\biggl|\int_{\Gamma} \frac{\hat{\mu}(\zeta,t,k) \hat{w}(\zeta,t,k) dk}{k - K_1}\biggr| = O(\epsilon \tau^{-\frac{1+\alpha}{2}}), \qquad \tau \to \infty, \  \zeta \in \mathcal{I},
\end{align}
uniformly with respect to $\zeta \in \mathcal{I}$.
Equations (\ref{matK1}), (\ref{Fn}), (\ref{Gn}), and (\ref{intSigmahatmuhatw}) yield
\begin{align}\nonumber
m(\zeta, t, K_1) 
= &\; I-\frac{1}{\sqrt{\tau}}\sum_{n=0}^2 \mathcal{A}^{-n}\bigg(
\frac{B(\zeta, t)}{k_0 - \omega^{-n}K_1}
+ \frac{\overline{\check{B}(\zeta, t)}}{k_0 + \omega^{-n}K_1}
	\\ \label{mzetatK1}
& - \frac{\mathcal{B} B(\zeta, t)\mathcal{B}}{k_0(1  - k_0 \omega^{-n}K_1)}
- \frac{\mathcal{B}\overline{\check{B}(\zeta, t)}\mathcal{B}}{k_0(1 + k_0\omega^{-n}K_1)}\bigg)\mathcal{A}^n
+ O(\epsilon \tau^{-\frac{1 + \alpha}{2}}).
\end{align}
Recalling that $K_1 = e^{\frac{i\pi}{6}}$, a tedious but straightforward computation gives (\ref{limlm12}).
This completes the proof of the theorem.
%Repeating the above arguments with $K_1$  replaced with $K_1^{-1}$, we find that $m(\zeta, t, K_1^{-1})$ is given by the right-hand side of (\ref{mzetatK1}) but with $K_1$ replaced by $K_1^{-1}$. Using these expressions in the right-hand side of 
%$$m(\zeta, t, K_1)  = \frac{1}{2}\big(m(\zeta, t, K_1) + \mathcal{B} m(\zeta, t, K_1^{-1}) \mathcal{B}\big)$$
%and recalling that $K_1 = e^{\frac{i\pi}{6}}$, a straightforward computation gives (\ref{limlm12}).
\end{proof}

\section{Long-time asymptotics in the similarity sector}\label{similaritysec}
Theorem \ref{mainth} gives explicit formulas in terms of $r(k)$ for the leading order asymptotics of the solution of the DP equation on the half-line in the asymptotic region
$$0 < c < \xi < 3, \qquad (3-\xi)^{\frac{3}{2}} t \to \infty,$$
where $\xi = x/t$.

\begin{theo}[Long-time asymptotics in the similarity region]\label{mainth}
Let $u_0$ and $\{g_j\}_0^2$ be functions in the Schwartz class $\mathcal{S}(\R_+)$ that satisfy the assumptions (\ref{qassumptions}). 
Suppose there exists a unique solution $u(x,t)$ of equation (\ref{DP}) with $\kappa = 1$ in the half-line domain $\Omega = \{x \geq 0, t \geq 0\}$ such that
\begin{itemize}
\item $u$ is a smooth function of $(x,t) \in \Omega$,
\item $u$ satisfies the initial and boundary conditions (\ref{boundaryconditions}),
\item $u(\cdot, t) \in \mathcal{S}(\R_+)$ for each $t \geq 0$.
\end{itemize}
Let $q(x,t) = (u(x,t) - u_{xx}(x,t) + 1)^{\frac{1}{3}}$.  
Define $\tilde{\Psi}_n(x,t,k)$ for $(x,t)$ in the set 
$$\{x\geq 0, t = 0\} \cup \{x=0, t \geq 0\}$$ 
by the linear integral equations (\ref{tildePhiFredholm}). Define $r(k)$ for $k \in (\frac{\sqrt{5}-1}{2}, \frac{\sqrt{5}+1}{2})$ by 
$$r(k) = (\tilde{\Psi}_{18}(0,0,k)^{-1}\tilde{\Psi}_7(0,0,k))_{21} e^{i(k - k^{-1})\int_0^\infty(q(x,0) - 1) dx}.$$
Suppose the set $\{k_j\}$ defined in Section \ref{prelsec} is empty and that
$$\sup_{\frac{\sqrt{5}-1}{2} < k < \frac{\sqrt{5}+1}{2}} |r(k)| < 1.$$

Then, for any $\alpha \in [\frac{1}{3},1)$ and $0 < c < 3$, the following asymptotic formulas are valid:
\begin{align}\label{qfinalxi}
q(x,t) = &\; 1 + \frac{b_1(\xi)}{\sqrt{t}}
\cos \big(b_2(\xi) t - \nu(\xi) \log(t) + b_3(\xi)\big) 
	\\ \nonumber
& + O\big((3 - \xi)^{-\frac{1 + 3\alpha}{4}} t^{-\frac{1+\alpha}{2}}\big), \qquad (3-\xi)^{\frac{3}{2}} t \to \infty, \  c \leq \xi < 3,
	\\\label{ufinalxi}
u(x,t) = &\; \frac{3 b_1(\xi)}{(1 + 4 \tilde{k}_0^2(\xi) )\sqrt{t}}
\cos\big(b_2(\xi) t - \nu(\xi) \log t + b_3(\xi)\big) 
	\\ \nonumber
&+ O\big((3 - \xi)^{-\frac{1 + 3\alpha}{4}} t^{-\frac{1+\alpha}{2}}\big),	
\qquad (3-\xi)^{\frac{3}{2}} t \to \infty, \  c \leq \xi < 3,
\end{align}
where the error terms are uniform with respect to $\xi$ in the given ranges and the functions $\{b_j(\xi)\}_1^3$, $\nu(\xi)$, and $\tilde{k}_0(\xi)$ are defined by
\begin{align*}
 b_1(\xi) = & \; \frac{1 - k_0^2(\xi)}{1 + k_0^2(\xi)} \sqrt{\frac{(3 + 4\tilde{k}_0^2(\xi))(1 + 4 \tilde{k}_0^2(\xi))\nu(\xi)}{3\tilde{k}_0 (\xi)- 4\tilde{k}_0^3(\xi)}},
	\\
 b_2(\xi) = &\; \frac{48 \tilde{k}_0^3(\xi)}{(1 + 4 \tilde{k}_0^2(\xi))^2},
\end{align*}
\begin{align*}
b_3(\xi) = & \;\frac{\pi}{4} - \chi_0(\xi) + \nu(\xi) \log\bigg(\frac{(4 \tilde{k}_0^2(\xi)+1)^2 (4 \tilde{k}_0^2(\xi)+3)}{576 \tilde{k}_0^3(\xi)(3 - 4 \tilde{k}_0^2(\xi))}\bigg)  
	\\
&+ \arg \Gamma(i\nu(\xi)) - \arg r(k_0(\xi))  - \arctan\bigg(\sqrt{3}\frac{1 + k_0^2(\xi)}{1 - k_0^2(\xi)}\bigg) 
	\\
& + \frac{3\tilde{k}_0(\xi)}{\pi} \int_{k_0(\xi)}^{\frac{1}{k_0(\xi)}} \log(1 - |r(s)|^2) \frac{1 + s^4}{1 + s^6} ds,	
	\\
\nu(\xi) = & -\frac{1}{2\pi} \log(1 - |r(k_0(\xi))|^2), \qquad \tilde{k}_0(\xi) = \sqrt{\frac{-2 \xi -3+\sqrt{24 \xi +9}}{8 \xi }},
\end{align*}
with
\begin{align} \nonumber
k_0(\xi) = & - \tilde{k}_0(\xi) + \sqrt{1 + \tilde{k}_0^2(\xi)}, 
	\\\nonumber
\chi_0(\xi) %= \im(2\chi(\xi, k_0(\xi)) - \chi(\xi, \omega k_0(\xi)) - \chi(\xi, \omega^2 k_0(\xi)))
= & -\frac{3}{2\pi} \int_{k_0(\xi)}^{\frac{1}{k_0(\xi)}}  \log\bigg(\frac{1 - |r(s)|^2}{1 - |r(k_0(\xi))|^2}\bigg) 
	\\ \label{chi0def}
&\times\frac{(k_0(\xi)-k_0^{11}(\xi)) (s^6+s^4)+(k_0^5(\xi)-k_0^7(\xi)) (s^{10}+1)}{k_0^{12}(\xi) s^6-k_0^6(\xi) (s^{12}+1)+s^6} ds.
\end{align}
\end{theo}

\subsection{Proof of Theorem \ref{mainth}}
The basic idea of the proof consists of deforming the contour of the RH problem (\ref{RHnu}) so that the jump matrix is exponentially close to the identity everywhere except near a certain set of critical points. Near these critical points the problem is solved locally and the asymptotic expansion of the solution is extracted from the local solution.

\begin{figure}
\begin{center}
 \begin{overpic}[width=.7\textwidth]{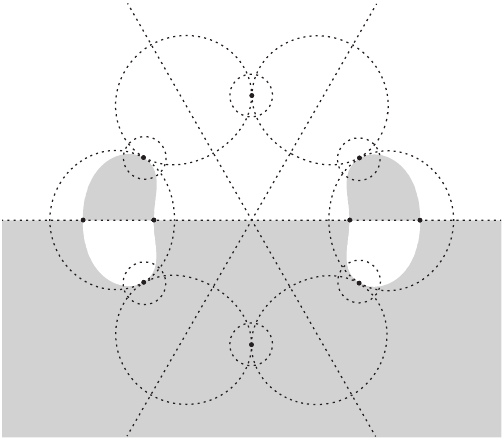}
  \put(44,4){\small $\re \Phi < 0$}
 \put(44,82){\small $\re \Phi > 0$}
       \put(69,57.3){\tiny $K_1$}
       \put(46,67.5){\tiny $K_2$}
       \put(28,57){\tiny $K_3$}
       \put(28,28.5){\tiny $K_4$}
       \put(46,18){\tiny $K_5$}
       \put(68.2,29){\tiny $K_6$}
              \put(68,40){\small $k_0$}
      \put(81,40){\small $k_0^{-1}$}
      \put(27,40){\small $-k_0$}
      \put(12,40){\small $-k_0^{-1}$}
   \end{overpic}
   \begin{figuretext}\label{rePhi.pdf}
      The regions in the complex $k$-plane where $\re \Phi(\zeta, k) < 0$ (shaded) and $\re \Phi(\zeta, k) > 0$ (white). The contour separating the sets $E_n$ is dashed.
             \end{figuretext}
   \end{center}
\end{figure}

Suppose $\zeta = y/t \in (0,3)$. In order to find the critical points, we note that, for fixed $\zeta$, the function $\Phi(\zeta, k)$ defined in (\ref{Phizetakdef}) can be written as $\Phi(\zeta, k) = \tilde{\Phi}(\zeta, \tilde{k}(k))$ where
$$\tilde{\Phi}(\zeta, \tilde{k}) = 2i\zeta \tilde{k} - \frac{6i\tilde{k}}{1 + 4 \tilde{k}^2}, \qquad \tilde{k}(k) = \frac{1}{2}\bigg(\frac{1}{k} - k\bigg).$$
The equation $\partial_{\tilde{k}}\tilde{\Phi}(\zeta, \tilde{k}) = 0$ has two real solutions $\pm \tilde{k}_0$, where $\tilde{k}_0 = \tilde{k}_0(\zeta)$ is defined by
$$\tilde{k}_0 = \sqrt{\frac{-2 \zeta -3 +\sqrt{24 \zeta +9}}{8\zeta }}.$$
Consequently, there are four real points at which $\partial \Phi/\partial k = 0$; these are given by $\pm k_0, \pm k_0^{-1}$, where $k_0 = k_0(\zeta)$ is defined by
$$k_0 = - \tilde{k}_0 + \sqrt{1 + \tilde{k}_0^2}.$$
It follows that for $\zeta \in (0,3)$ the RH problem (\ref{RHnu}) has twelve critical points associated with it (the four points $\pm k_0, \pm k_0^{-1}$ as well as the eight points obtained by multiplying these four points by $\omega$ and $\omega^2$). 
Note that $0 < \tilde{k}_0 < 1/2$ and $\frac{\sqrt{5}-1}{2} < k_0 < 1$ for $0 < \zeta < 3$. We henceforth assume that $\zeta \in \mathcal{I}$, where $\mathcal{I} = [c, 3)$ and $c > 0$ is a small constant.  
The signature table of $\re \Phi$ for $\zeta$ in this range is displayed in Figure \ref{rePhi.pdf}. The form of the signature table can be deduced from the explicit formula
\begin{align*}
& \re \Phi(\zeta, k_1 + ik_2) = k_2 \left(k_1^2+k_2^2+1\right) \bigg(\frac{\zeta}{k_1^2+k_2^2} + \frac{P(k_1,k_2)}{Q(k_1,k_2)}\bigg),
\end{align*}
where
\begin{align*}
P(k_1,k_2) = &\; 3 \left(k_1^4+k_1^2 \left(2  k_2^2-3\right)+k_2^4+k_2^2+1\right), 
	\\
Q(k_1,k_2) = &\; k_1^8+k_1^6 \left(4 k_2^2-2\right)+k_1^4 \left(6 k_2^4-2 k_2^2+3\right)
	\\
& +2 k_1^2 \left(2 k_2^6+k_2^4-5k_2^2-1\right) +\left(k_2^4+k_2^2+1\right)^2.
\end{align*}
The equation $P = 0$ implies $\re(z_2 - z_1) = 0$ and defines the curves $(\bar{D}_1 \cap \bar{D}_7) \cup (\bar{D}_6 \cap \bar{D}_{18})$ and $(\bar{D}_3 \cap \bar{D}_{12}) \cup (\bar{D}_4 \cap \bar{D}_{13})$.
The remainder of the proof proceeds through seven steps.

\begin{figure}
\begin{center}
 \begin{overpic}[width=.9\textwidth]{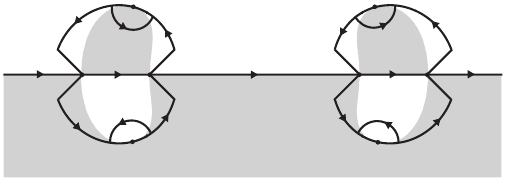}
  \put(44,12){$\re \Phi < 0$}
 \put(44,29){$\re \Phi > 0$}
       \put(72,35.6){\small $K_1$}
       \put(25,35.6){\small $K_3$}
       \put(25,5){\small $K_4$}
       \put(73,4.8){\small $K_6$}
       \put(71,18.5){\small $k_0$}
      \put(80,18.2){\small $k_0^{-1}$}
      \put(25.5,18.5){\small $-k_0$}
      \put(15.4,18.2){\small $-k_0^{-1}$}
   \end{overpic}
   \begin{figuretext}\label{contour2.pdf}
      The jump contour $\hat{\Gamma}$ for  $k$  near $\R$ together with the regions where $\re \Phi < 0$ (shaded) and $\re \Phi > 0$ (white).
             \end{figuretext}
   \end{center}
\end{figure}

\begin{figure}
\begin{center}
 \begin{overpic}[width=.75\textwidth]{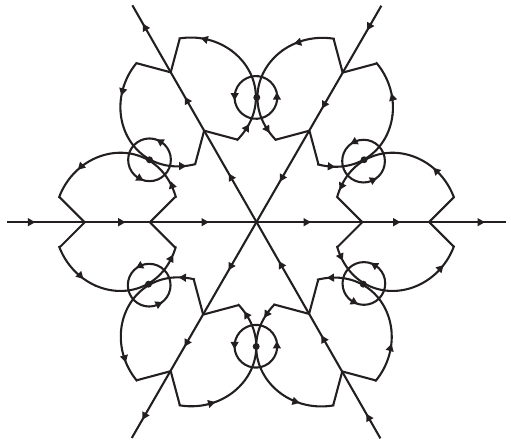}
       \put(71,53.3){\tiny $K_1$}
      \put(46.4,66.6){\tiny $K_2$}
      \put(28,57){\tiny $K_3$}
      \put(28.8,29){\tiny $K_4$}
      \put(50.4,17.8){\tiny $K_5$}
      \put(70.6,32.2){\tiny $K_6$}
      \put(70.2,40.2){\scriptsize $k_0$}
      \put(79.8,40.2){\scriptsize $k_0^{-1}$}
      \put(25,40.2){\scriptsize $-k_0$}
      \put(13.4,40.2){\scriptsize $-k_0^{-1}$}
     \put(93,48){\small $F_1$}
      \put(58,48){\small $F_1$}
     \put(4,48){\small $F_3$}
      \put(39,48){\small $F_3$}
     \put(4,37){\small $F_4$}
      \put(39,37){\small $F_4$}
       \put(93,37){\small $F_6$}
      \put(58,37){\small $F_6$}
      \put(75,48){\small $F_7$}
       \put(20,37){\small $F_{12}$}
      \put(20,48){\small $F_{13}$}
        \put(75,37){\small $F_{18}$}
   \end{overpic}
   \begin{figuretext}\label{contour2all.pdf}
      The jump contour $\hat{\Gamma}$ for $\hat{M}$. 
       \end{figuretext}
   \end{center}
\end{figure}

\bigskip
{\bf Step 1: Deform contour.}
We begin by deforming the contour so that it passes through the twelve critical points, see Figures \ref{contour2.pdf} and \ref{contour2all.pdf}. 
For $k$ near $k_0$ and $k_0^{-1}$, the contour deformation is achieved by defining $\hat{M} = M J_{1,7}^{-1}$ for $k \in F_1 \cap E_7$ and $\hat{M} = M J_{6,18}^{-1}$ for $k \in F_6 \cap E_{18}$, where the sets $\{F_n\}$ are as in Figure \ref{contour2all.pdf}. We define $\hat{M}$ analogously near the other critical points and set $\hat{M} = M$ otherwise. 

Let $\hat{\Gamma}$ denote the contour displayed in Figure \ref{contour2all.pdf}. 
The matrix $J_{6,18}$ involves the factor $h(k)e^{t\Phi}$, which is bounded and analytic in $F_6 \cap E_{18}$. Similarly, $J_{1,7}$ involves the factor $\overline{h(\bar{k})}e^{-t\Phi}$, which is bounded and analytic in $F_1 \cap E_7$. 
Thus, except for possible singularities at the points $\{\varkappa_j\}_1^6$, $\hat{M}(y,t,k)$ is a bounded and analytic function of $k \in \hat{\C} \setminus \hat{\Gamma}$.
Using Lemma \ref{deformationlemma}, we conclude that $\hat{N} = (1,1,1) \hat{M}$ satisfies the $L^3$-RH problem
\begin{align}\label{RHNhat}
\begin{cases} \hat{N}(y,t, \cdot) \in (1,1,1) + \dot{E}^3(\hat{\C} \setminus \hat{\Gamma}), 
	\\
\hat{N}_+(y,t,k) =  \hat{N}_-(y,t,k) \hat{J}(y,t,k) \quad \text{for a.e.} \ k \in \hat{\Gamma},
\end{cases} 
\end{align}
where the expression for $\hat{J} = \hat{J}_{m,n}$ for $k \in \bar{F}_m \cap \bar{F}_n$ coincides with the expression for $J$ on $\bar{E}_m \cap \bar{E}_n$, see equations (\ref{Jformula1}), (\ref{Jformula2}), and (\ref{Jformula3}).

\bigskip
{\bf Step 2: Conjugate.}
On the circles where the $E_n$'s and $E_{n+18}$'s meet, the jump matrix $\hat{J}(y,t,k)$ has exponential decay as $t \to \infty$. Indeed, the circles where the $E_n$'s and $E_{n+18}$'s meet are the small circles centered at $K_j$, $j = 1, \dots, 6$, see Figure \ref{Econtour.pdf}. By symmetry, it is enough to consider the small circle centered at $K_1$. On this circle, the jump matrix is given by the three formulas (\ref{Jformula3}) in Lemma \ref{Jlemma}, see Figure 4.
The $(12)$, $(13)$, and $(23)$ entries of the matrices in (\ref{Jformula3}) involve the exponentials $e^{-t\Phi(\zeta, k)}$, $e^{t\Phi(\zeta, \omega^2 k)}$, and $e^{-t\Phi(\zeta, \omega k)}$, respectively; the decay now follows from (\ref{Jformula3}) and the signature table of $\re \Phi$, see Figure \ref{rePhi.pdf}.

In order to arrive at a jump matrix with the appropriate decay properties also on the remaining part of the contour, we need to perform a triangular factorization of $\hat{J}$. Such a factorization can be achieved by conjugating the RH problem as follows.
Let
$$\Delta(\zeta, k) = \begin{pmatrix} \delta(\zeta, k) \delta(\zeta, \omega^2k)^{-1} & 0 & 0 \\ 0 & \delta(\zeta, k)^{-1}\delta(\zeta, \omega k) & 0 \\ 0 & 0 & \delta(\zeta, \omega^2k) \delta(\zeta, \omega k)^{-1}   \end{pmatrix},$$
where, for $k \in \C \setminus \R$,
\begin{align}\nonumber
 \delta(\zeta, k) = \exp\bigg\{&\frac{1}{4\pi i} \int_{-\frac{1}{k_0}}^{-k_0}\log(1 - |\check{r}(s)|^2)\bigg(\frac{1}{s - k} - \frac{1}{s - \frac{1}{k}}\bigg)ds
	\\ \label{deltadef}
&+ \frac{1}{4\pi i} \int_{k_0}^{\frac{1}{k_0}} \log(1 - |r(s)|^2)\bigg(\frac{1}{s - k} - \frac{1}{s - \frac{1}{k}}\bigg) ds\bigg\}.
\end{align}
The identities
$$\delta(\zeta, k) = \frac{1}{\overline{\delta(\zeta, \bar{k})}} = \frac{1}{\delta(\zeta, k^{-1})},$$
imply that $\Delta(\zeta, k)$ obeys the three symmetries in (\ref{symmetries}).
The function $\delta$ satisfies 
\begin{align*}
 \delta_+(\zeta, k) = \begin{cases} 
\delta_-(\zeta, k) (1 - |r(k)|^2), & k \in (k_0, k_0^{-1}), \\
\delta_-(\zeta, k) (1 - |\check{r}(k)|^2), & k \in (-k_0^{-1}, -k_0), \\
\delta_-(\zeta, k), & \text{otherwise},
\end{cases}
\end{align*}
and
\begin{align*}
& \delta(\zeta, k) = e^{i\varphi} + O(k^{-1}), \qquad k \to \infty, \  k \in \C,
\end{align*}
where the constant $\varphi \in \R$ is given by
$$\varphi = \frac{1}{4\pi} \int_{-\frac{1}{k_0}}^{-k_0}\log(1 - |\check{r}(s)|^2) \frac{ds}{s}
+ \frac{1}{4\pi} \int_{k_0}^{\frac{1}{k_0}} \log(1 - |r(s)|^2)  \frac{ds}{s}.$$
Moreover, the representation 
$$\delta(\zeta,k) = \bigg(\frac{(k_0^{-1} - k)(k_0 - k^{-1})}{(k_0 - k)(k_0^{-1} - k^{-1})}\bigg)^{\frac{i\nu}{2}} \bigg(\frac{(k_0 + k)(k_0^{-1} + k^{-1})}{(k_0^{-1} + k)(k_0 + k^{-1})}\bigg)^{\frac{i\check{\nu}}{2}}  e^{\chi(\zeta, k)},$$
where
\begin{align*}
\chi(\zeta, k) = & \;\frac{1}{4\pi i}\int_{-\frac{1}{k_0}}^{-k_0} \log\bigg(\frac{1 - |\check{r}(s)|^2}{1 - |\check{r}(-k_0)|^2}\bigg) \bigg(\frac{1}{s-k} - \frac{1}{s - \frac{1}{k}}\bigg) ds
	\\
& + \frac{1}{4\pi i}\int_{k_0}^{\frac{1}{k_0}} \log\bigg(\frac{1 - |r(s)|^2}{1 - |r(k_0)|^2}\bigg)  \bigg(\frac{1}{s-k} - \frac{1}{s - \frac{1}{k}}\bigg) ds,
\end{align*}
shows that
$$\delta(\zeta, \cdot), \delta(\zeta, \cdot)^{-1} \in E^\infty(\C \setminus \R).$$
We conclude that
$$\Delta(\zeta, \cdot), \Delta(\zeta, \cdot)^{-1} \in I + \dot{E}^3(\C \setminus \hat{\Gamma}) \cap E^\infty(\C \setminus \hat{\Gamma}).$$

The function $\tilde{M}$ defined by
$$\tilde{M}(y,t,k) = \hat{M}(y,t,k)\Delta(\zeta, k)$$
satisfies the jump condition $\tilde{M}_+ = \tilde{M}_- \tilde{J}$ on $\hat{\Gamma}$ with $\tilde{J} = \Delta_-^{-1} \hat{J} \Delta_+$. Define $r_2(k)$ by
$$r_2(k) = \frac{r(k)}{1 - r(k)\overline{r(\bar{k})}}, \qquad k \in \bar{E}_7 \cap \bar{E}_{18}.$$
We find from (\ref{Jformula1}) that
\begin{align*}
  \tilde{J} = 
 \begin{cases} 
 b_l^{-1}b_u, & k \in \bar{F}_1 \cap \bar{F}_6, 
	\\
B_u^{-1} B_l, & k \in \bar{F}_7 \cap \bar{F}_{18}, 
	\\
\begin{pmatrix}
 1 & \frac{\delta(\zeta, \omega k) \delta(\zeta, \omega^2k)}{\delta(\zeta, k)^2} \overline{h(\bar{k})} e^{-t\Phi} & 0 \\
0 & 1 & 0 \\
 0 & 0 & 1 
\end{pmatrix}, & k \in \bar{F}_1 \cap \bar{F}_7, \\
\begin{pmatrix}
 1 & 0 & 0 \\
\frac{\delta(\zeta, k)^2}{\delta(\zeta, \omega k) \delta(\zeta, \omega^2k)} h(k) e^{t\Phi} & 1 & 0 \\
 0 & 0 & 1 
\end{pmatrix}, & k \in \bar{F}_6 \cap \bar{F}_{18}, 
	\\
\check{b}_l^{-1}\check{b}_u, & k \in \bar{F}_3 \cap \bar{F}_4, 
	\\
\check{B}_u^{-1} \check{B}_l, & k \in \bar{F}_{12} \cap \bar{F}_{13}, 
	\\
\begin{pmatrix}
 1 & \frac{\delta(\zeta, \omega k) \delta(\zeta, \omega^2k)}{\delta(\zeta, k)^2}  \overline{\check{h}(\bar{k})} e^{-t\Phi}  & 0 \\
0 & 1 & 0 \\
 0 & 0 & 1 
\end{pmatrix}, & k \in \bar{F}_3 \cap \bar{F}_{12}, 
	\\
\begin{pmatrix}
 1 & 0 & 0 \\
\frac{\delta(\zeta, k)^2}{\delta(\zeta, \omega k) \delta(\zeta, \omega^2k)} \check{h}(k) e^{t\Phi} & 1 & 0 \\
 0 & 0 & 1 
\end{pmatrix}, & k \in \bar{F}_4 \cap \bar{F}_{13},
\end{cases}
\end{align*}
where
\begin{align*}\nonumber
& b_l =  \begin{pmatrix} 1 & 0 & 0 \\
-\frac{\delta(\zeta, k)^2}{\delta(\zeta, \omega k) \delta(\zeta, \omega^2k)} r_1(k) e^{t\Phi(\zeta,k)}	& 1 & 0 \\ 0 & 0 & 1 \end{pmatrix}, 
	 \\
& B_u = \begin{pmatrix} 1 & \frac{\delta(\zeta, \omega k) \delta(\zeta, \omega^2 k)}{\delta_-(\zeta, k)^2} \overline{r_2(\bar{k})} e^{-t\Phi(\zeta,k)} & 0	\\
0 & 1 & 0 \\ 0 & 0 & 1 \end{pmatrix}, 
	\\ 
&B_l = \begin{pmatrix} 1 &0 & 0 \\
 \frac{\delta_+(\zeta, k)^2}{\delta(\zeta, \omega k) \delta(\zeta, \omega^2 k)} r_2(k) e^{t\Phi(\zeta,k)} & 1 & 0 \\ 0 & 0 & 1 \end{pmatrix},
	\\
& b_u =  \begin{pmatrix} 1 & -\frac{\delta(\zeta, \omega k) \delta(\zeta, \omega^2k)}{\delta(\zeta, k)^2} \overline{r_1(\bar{k})} e^{-t\Phi(\zeta,k)}& 0	\\
0	& 1 & 0 \\ 0 & 0 & 1 \end{pmatrix},
\end{align*}
and $\check{b}_l$, $\check{B}_u$, $\check{B}_l$, $\check{b}_u$ are given by analogous expressions but with $r$ and $h$ replaced with  $\check{r}$ and $\check{h}$, respectively.

\medskip
{\bf Step 3: Introduce analytic approximations.}
The next step consists of splitting each of the functions $r_j$, $j = 1,2$, into an analytic part $r_{j,a}$ and a small remainder $r_{j,r}$. 

\begin{figure}
\begin{center}
\bigskip
 \begin{overpic}[width=.98\textwidth]{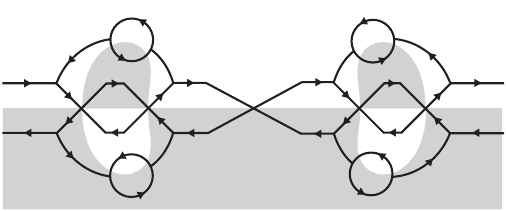}
 \put(92,23){\small $U_1$}
 \put(82,27){\small $U_2$}
 \put(76,23){\small $U_3$}
 \put(76,18){\small $U_4$}
 \put(82,14){\small $U_5$}
 \put(92,18){\small $U_6$}
 \put(61,23){\small $U_1$}
 \put(61,18){\small $U_6$}
 \put(36,23){\small $U_7$}
 \put(28,27){\small $U_8$}
 \put(21,23){\small $U_9$}
 \put(21,18){\small $U_{10}$}
 \put(27,14){\small $U_{11}$}
 \put(36,18){\small $U_{12}$}
 \put(4,23){\small $U_7$}
 \put(4,18){\small $U_{12}$}
 \put(92, 28){\small $V_1$}
 \put(61, 28){\small $V_1$}
  \put(36, 28){\small $V_3$}
 \put(4, 28){\small $V_3$}
  \put(4,13){\small $V_4$}
 \put(36,13){\small $V_4$}
  \put(61,13){\small $V_6$}
 \put(92,13){\small $V_6$}
   \put(45,9){\small $\re \Phi < 0$}
 \put(45,33){\small $\re \Phi > 0$}
   \end{overpic}
   \begin{figuretext}\label{contour3.pdf}
      The open sets $\{U_j\}_1^{12}$ and the jump contour $\Gamma$ for $k$  near $\R$. 
       \end{figuretext}
   \end{center}
\end{figure}

Define open sets $U_j := U_j(\zeta)$, $j = 1, \dots, 12$, as in Figure \ref{contour3.pdf} so that $\re \Phi(\zeta, k) > 0$ in $U_1 \cup U_4 \cup U_7 \cup U_{10}$ and $\re \Phi(\zeta, k) < 0$ in $U_3 \cup U_6 \cup U_9 \cup U_{12}$.
Write $U_6 = U_6^+ \cup U_6^-$, where $U_6^+ = U_6 \cap \{\re k > k_0^{-1}\}$ and $U_6^- = U_6 \cap \{\re k < k_0\}$. 

\begin{lemma}[Analytic approximations of $r_1(k)$ and $r_2(k)$]\label{decompositionlemma}
There exist decompositions
\begin{align*}
& r_1(k) = r_{1,a}(y, t, k) + r_{1,r}(y, t, k), \qquad k \in \bar{U}_1 \cap \bar{U}_6, 
	\\
& r_2(k) = r_{2,a}(y, t, k) + r_{2,r}(y, t, k), \qquad k \in \bar{U}_3 \cap \bar{U}_4, 
\end{align*}
where the functions $\{r_{j,a}, r_{j,r}\}_{j=1}^2$ have the following properties:
\begin{enumerate}[$(a)$]
\item For each $\zeta \in \mathcal{I}$ and each $t > 0$, $r_{1,a}(y, t, k)$ is continuous for $k \in \bar{U}_6$ and analytic for $k \in U_6$, while $r_{2,a}(y, t, k)$ is continuous for $k \in \bar{U}_3$ and analytic for $k \in U_3$. 

\item There exists a constant $C$ independent of $\zeta, t, k$ such that
\begin{align}\nonumber
& \begin{cases}
|r_{1, a}(y, t, k) - r_1(k_0)| \leq 
C |k - k_0| e^{\frac{t}{4}|\re \Phi(\zeta,k)|},  &k \in \bar{U}_6^-, \\
|r_{1, a}(y, t, k) - r_1(k_0^{-1})| \leq 
C |k - k_0^{-1}| e^{\frac{t}{4}|\re \Phi(\zeta,k)|},  &k \in \bar{U}_6^+, \\
|r_{1, a}(y, t, k)| \leq \frac{C}{1 + |k|} e^{\frac{t}{4}|\re \Phi(\zeta,k)|}, & k \in \bar{U}_6, 
\end{cases}
\end{align}
and
\begin{align}\label{rjaestimateb}
&\begin{cases} |r_{2, a}(y, t, k) - r_2(k_0)| \leq 
C |k - k_0| e^{\frac{t}{4}|\re \Phi(\zeta,k)|}, \\
 |r_{2, a}(y, t, k) - r_2(k_0^{-1})| \leq 
C |k - k_0^{-1}| e^{\frac{t}{4}|\re \Phi(\zeta,k)|}, 
\end{cases} \  k \in \bar{U}_3,
\end{align}
for $\zeta \in \mathcal{I}$ and $t > 0$.

\item The $L^1, L^2$, and $L^\infty$ norms on $\bar{U}_1 \cap \bar{U}_6$ of the function $r_{1,r}(y, t, \cdot)$ are $O(t^{-3/2})$ as $t \to \infty$ uniformly with respect to $\zeta \in \mathcal{I}$.

\item The $L^1, L^2$, and $L^\infty$ norms on $\bar{U}_3 \cap \bar{U}_4$ of the function $r_{2,r}(y, t, \cdot)$ are $O(t^{-3/2})$ as $t \to \infty$ uniformly with respect to $\zeta \in \mathcal{I}$.

\item The following symmetries are valid:
\begin{align}\label{rsymmetries}
r_{j,a}(y, t, k) = \overline{r_{j,a}(y, t, \bar{k}^{-1})}, \quad
r_{j,r}(y, t, k) = \overline{r_{j,r}(y, t, \bar{k}^{-1})}, \qquad j = 1, 2.
\end{align}
\end{enumerate}
\end{lemma}
\begin{proof}
Analytic approximations of this type were introduced in \cite{DZ1993}. The proof here follows the presentation of \cite{Lmkdvhalfline, Lnonlinearsteepest} (see in particular Lemma 4.8 of \cite{Lnonlinearsteepest}). 
We will derive the decomposition of $r_1$ in $U_6^+$. This decomposition can easily be extended to $U_6^-$ by means of the symmetry $r_1(k) = \overline{r_1(\bar{k}^{-1})}$. The decomposition of $r_2(k)$ can be derived in a similar way.

Our assumption that $u_0$ and $\{g_j\}_0^2$ belong to the Schwartz class $\mathcal{S}(\R_+)$ implies that $r_1(k)$ has the following properties: 
\begin{itemize}
\item $r_1(k)$ is a smooth function of $k \in (0, \infty)$.
 
\item There are functions $\{p_j(\zeta)\}_0^7$ such that
\begin{align*}
r_1^{(n)}(k) = \frac{d^n}{dk^n}\bigg(\sum_{j=0}^7 p_j(\zeta) (k - k_0^{-1})^j\bigg) + O((k-k_0^{-1})^{8-n}) 
\end{align*}
as $k \in \R$ approaches $k_0^{-1}$ for $n = 0,1,2$ (in fact, $p_j(\zeta) := r_1^{(j)}(k_0^{-1})/j!$).

\item There are constants $\{r_{1,j}\}_1^2 \subset \C$ such that
\begin{align*}
r_1^{(n)}(k) = \frac{d^n}{dk^n}\bigg(\frac{r_{1,1}}{k} + \frac{r_{1,2}}{k^2} \bigg) + O(k^{-3}), \qquad k \to \infty, \ n = 0,1,2.
\end{align*}
\end{itemize}

We let
$$f_0(\zeta, k) = \sum_{j=1}^{10} \frac{a_j(\zeta)}{k^j},$$
where the coefficients $\{a_j(\zeta)\}_1^{10}$ are determined by the conditions
\begin{align*}
& f_0(\zeta, k) = \sum_{j=0}^7 p_j(\zeta) (k-k_0^{-1})^j + O((k-k_0^{-1})^8), \qquad k \to k_0^{-1}, \ \zeta \in \mathcal{I},
	\\
& f_0(\zeta, k) = \frac{r_{1,1}}{k} + \frac{r_{1,2}}{k^2} + O(k^{-3}), \qquad k \to \infty, \ \zeta \in \mathcal{I}.
\end{align*}
These ten linear conditions determine the coefficients $a_j(\zeta)$ uniquely and we have $\sup_{\zeta \in \mathcal{I}} |a_j(\zeta)| < \infty$ for each $j$. For each $\zeta \in \mathcal{I}$, $f_0(\zeta, k)$ is a rational function of $k \in \C$ with no poles in $\bar{U}_6^+$, which coincides with $r_1(k)$ to seventh order at $k_0$ and to second order at $\infty$. In other words, the function $f(\zeta,k)$ defined by $f(\zeta,k) = r_1(k) - f_0(\zeta,k)$ satisfies
\begin{align}\label{fcoincide2}
 \frac{\partial^n f}{\partial k^n} (\zeta, k) =
\begin{cases}
 O((k-k_0^{-1})^{8 - n}), & k \to k_0^{-1}, 
	\\
O(k^{-3}), & k \to \infty, 
 \end{cases}
 \   k \in \R, \ \zeta \in \mathcal{I}, \ n = 0,1,2,
\end{align}
where the error terms are uniform with respect to $\zeta \in \mathcal{I}$.

\begin{figure}
\begin{center}
 \begin{overpic}[width=.6\textwidth]{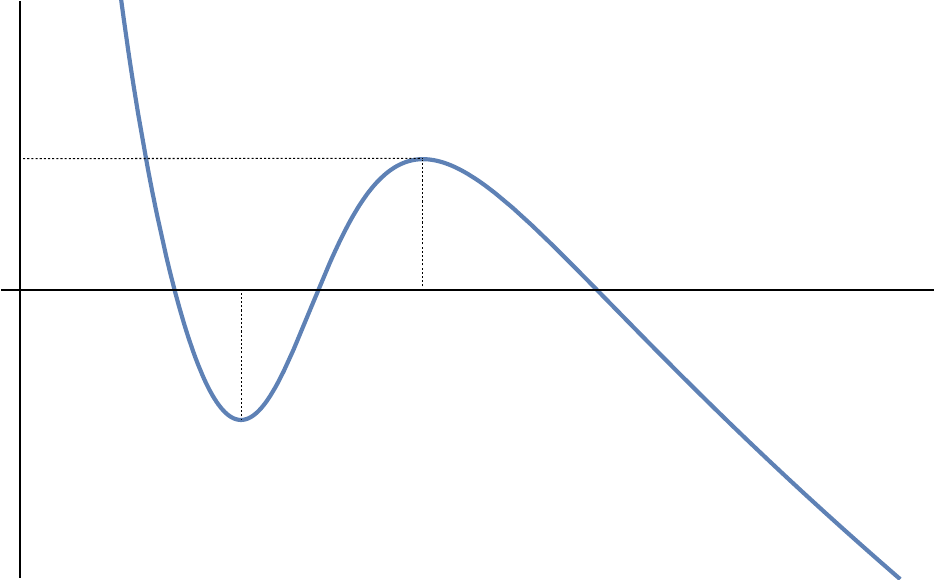}
      \put(102,30.2){$k$}
      \put(.7,65){$\phi$}
      \put(23.8,33.2){$k_0$}
      \put(42,26){$k_0^{-1}$}
      \put(-4,44.5){$\phi_0$}
       \end{overpic}
     \begin{figuretext}\label{phiofk.pdf}
       Graph of the function $k \mapsto \phi(\zeta, k)$ defined in (\ref{phidef}).
     \end{figuretext}
     \end{center}
\end{figure}

The decomposition of $r_1(k)$ can now be derived as follows.
Define the new variable $\phi$ by
\begin{align}\label{phidef}
  \phi = -i\Phi(\zeta, k) = 2\tilde{k} \zeta - \frac{6\tilde{k}}{1 + 4\tilde{k}^2}, \qquad \tilde{k} = \frac{k^{-1} - k}{2},
\end{align}  
and let $\phi_0 := \phi_0(\zeta)$ denote the value of $\phi$ at $k_0^{-1}$, i.e., $\phi_0(\zeta) = \phi(\zeta, k_0^{-1})$. For each  $\zeta \in \mathcal{I}$, the map $k \mapsto \phi = \phi(\zeta, k)$ is a decreasing bijection  $[k_0^{-1}, \infty) \to (-\infty, \phi_0]$, see Figure \ref{phiofk.pdf}.

Hence we may define a function $F(\zeta, \phi)$ by
\begin{align}\label{Fdef2}
F(\zeta, \phi) = \begin{cases} \frac{k^3}{k-k_0^{-1}} f(\zeta, k), &  \phi < \phi_0, \\
0, & \phi \geq \phi_0,
\end{cases} \qquad \zeta \in \mathcal{I}, \ \phi \in \R.
\end{align}
For each $\zeta \in \mathcal{I}$, the function $F(\zeta, \cdot)$ is smooth on $\R \setminus \{\phi_0\}$ and 
\begin{align}\label{dnFdphin}
\frac{\partial^n F}{\partial \phi^n}(\zeta, \phi) = \bigg(\frac{1}{\partial \phi/\partial k} \frac{\partial }{\partial k}\bigg)^n 
\bigg[\frac{k^3}{k-k_0^{-1}} f(\zeta, k)\bigg], \qquad \zeta \in \mathcal{I}, \ \phi < \phi_0.
\end{align}
We have
\begin{align}\label{dphidkatinfty}
\frac{\partial \phi}{\partial k}(\zeta, k) = -\zeta + O(k^{-2}), \qquad k \to \infty, \ \zeta \in \mathcal{I}. 
\end{align}
On the other hand, a computation shows that
\begin{align*}
&  \frac{\partial \phi}{\partial k}(\zeta, k) = d_1(\zeta) (k-k_0^{-1}) + d_2(\zeta) (k-k_0^{-1})^2 + O((k-k_0^{-1})^3), \qquad \zeta \in \mathcal{I},
\end{align*}
as $k \to k_0^{-1}$, where the coefficients $\{d_j(\zeta)\}_1^2$ are bounded for $\zeta \in \mathcal{I}$ and satisfy $d_1(\zeta) < 0$ for $\zeta \in (0,3)$, $d_1(3) = 0$, and $d_2(3) \neq 0$. In particular, there exist constants $c,C > 0$ independent of $k, \zeta$ such that
\begin{align}\label{dphidkatk0m1}
& \bigg|\frac{\partial \phi}{\partial k}(\zeta, k)\bigg| \geq c |k-k_0^{-1}|^2 \quad \text{and} \quad \bigg|\frac{\partial^2 \phi}{\partial k^2}(\zeta, k)\bigg| < C
\end{align}
for all $k > k_0^{-1}$ and $\zeta \in \mathcal{I}$.
Equations (\ref{fcoincide2}), (\ref{dnFdphin}), (\ref{dphidkatinfty}), and (\ref{dphidkatk0m1}) show that $F(\zeta, \cdot) \in C^1(\R)$ for each $\zeta$ and that
\begin{align*}
\bigg| \frac{\partial^n F}{\partial \phi^n}(\zeta, \phi)\bigg| \leq 
 \frac{C}{1 + |\phi|}, \qquad \phi < \phi_0, \ \zeta \in \mathcal{I}, \ n = 0,1,2,
\end{align*}
where $C > 0$ is independent of $\zeta, \phi,n$.
Hence
\begin{align}\label{supdFdphi}
\sup_{\zeta \in \mathcal{I}} \bigg\|\frac{\partial^n F}{\partial \phi^n}(\zeta, \cdot)\bigg\|_{L^2(\R)} < \infty, \qquad n = 0,1,2.
\end{align}
In particular, $F(\zeta, \cdot)$ belongs to the Sobolev space $H^2(\R)$ for each $\zeta \in \mathcal{I}$.
We conclude that the Fourier transform $\hat{F}(\zeta, s)$ defined by
\begin{align}\label{Fhatdef}
\hat{F}(\zeta, s) = \frac{1}{2\pi} \int_{\R} F(\zeta, \phi) e^{-i\phi s} d\phi
\end{align}
satisfies
\begin{align}\label{FFhat}
F(\zeta, \phi) =  \int_{\R} \hat{F}(\zeta, s) e^{i\phi s} ds
\end{align}
and
\begin{align}\label{x2Fhat}
\sup_{\zeta \in \mathcal{I}} \|s^2 \hat{F}(\zeta, s)\|_{L^2(\R)} < \infty.
\end{align}
In view of (\ref{Fdef2}) and (\ref{FFhat}), we find
$$ \frac{k-k_0^{-1}}{k^3}\int_{\R} \hat{F}(\zeta, s) e^{s\Phi(\zeta,k)} ds 
= \begin{cases} f(\zeta, k), &  k > k_0^{-1}, \\
0, & k \leq k_0^{-1}, 
 \end{cases} \qquad \zeta \in \mathcal{I}.$$
Let us write
$$f(\zeta, k) = f_a(y, t, k) + f_r(y, t, k), \qquad \zeta \in \mathcal{I}, \ t > 0, \ k \geq k_0,$$
where the functions $f_a$ and $f_r$ are defined by
\begin{align*}
& f_a(y,t,k) = \frac{k-k_0^{-1}}{k^3}\int_{-\frac{t}{4}}^\infty \hat{F}(\zeta,s) e^{s\Phi(\zeta,k)} ds, \qquad \zeta \in \mathcal{I}, \ t > 0, \ k \in \bar{U}_6^+,  
\end{align*}
and
\begin{align*}
& f_r(y,t,k) = \frac{k-k_0^{-1}}{k^3}\int_{-\infty}^{-\frac{t}{4}} \hat{F}(\zeta,s) e^{s\Phi(\zeta,k)} ds,\qquad \zeta \in \mathcal{I}, \ t > 0, \  k \geq k_0^{-1}.
\end{align*}
The function $f_a(y, t, \cdot)$ is continuous in $\bar{U}_6^+$ and analytic in $U_6^+$, because $\re \Phi \leq 0$ in $\bar{U}_6$. 
Furthermore,
\begin{align}\label{faest}
 |f_a(y, t, k)| 
&\leq \frac{|k-k_0^{-1}|}{|k|^3} \|\hat{F}(\zeta,\cdot)\|_{L^1(\R)}  \sup_{s \geq -\frac{t}{4}} e^{s \re \Phi(\zeta,k)}
	\\ \nonumber
&\leq \frac{C|k-k_0^{-1}|}{|k|^3}  e^{\frac{t}{4} |\re \Phi(\zeta,k)|}, \qquad \zeta \in \mathcal{I}, \ t > 0, \ k \in \bar{U}_6^+,
\end{align}
and
\begin{align}\label{frest}
|f_r(y, t, k)| & \leq \frac{|k-k_0^{-1}|}{|k|^3}  \int_{-\infty}^{-\frac{t}{4}} s^2 |\hat{F}(\zeta,s)| s^{-2} ds
	\\ \nonumber
& \leq \frac{C}{1 + |k|^2}  \| s^2 \hat{F}(\zeta,s)\|_{L^2(\R)} \sqrt{\int_{-\infty}^{-\frac{t}{4}} s^{-4} ds}  
 	\\ \nonumber
&  \leq \frac{C}{1 + |k|^2} t^{-3/2}, \qquad \zeta \in \mathcal{I}, \ t > 0, \ k \geq k_0^{-1}.
\end{align}
Hence the $L^1$, $L^2$, and $L^\infty$ norms of $f_r$ on $(k_0^{-1}, \infty)$ are $O(t^{-3/2})$ uniformly with respect to $\zeta \in \mathcal{I}$. Letting
\begin{align*}
& r_{1,a}(y, t, k) = f_0(\zeta, k) + f_a(y, t, k), \qquad k \in \bar{U}_6^+,
	\\
& r_{1,r}(y, t, k) = f_r(y, t, k), \qquad k \geq k_0^{-1}.
\end{align*}
we find a decomposition of $r_1$ for $k > k_0^{-1}$ with the properties listed in the statement of the lemma. 
\end{proof}

Using the decompositions of $r_1(k)$ and $r_2(k)$ established in Lemma \ref{decompositionlemma}, we can factorize the matrices $\{b_u, b_l, B_u, B_l\}$ as follows:
$$b_u = b_{u,a}b_{u,r}, \qquad
b_l = b_{l,a}b_{l,r}, \qquad
B_u = B_{u,a}B_{u,r}, \qquad
B_l = B_{l,a}B_{l,r},$$
where $\{b_{u,a}, b_{l,a}, B_{u,a}, B_{l,a}\}$ and $\{b_{u,r}, b_{l,r}, B_{u,r}, B_{l,r}\}$ are defined by the same formulas as $\{b_u, b_l, B_u, B_l\}$ except that the functions $\{r_j(k)\}_1^2$ are replaced by $\{r_{j,a}(k)\}_1^2$ and $\{r_{j,r}(k)\}_1^2$, respectively.

In the same way, we introduce decompositions
\begin{align*}
& \check{r}_1(k) = \check{r}_{1,a}(y, t, k) + \check{r}_{1,r}(y, t, k), \qquad k \in \bar{U}_7 \cap \bar{U}_{12}, 
	\\
& \check{r}_2(k) = \check{r}_{2,a}(y, t, k) + \check{r}_{2,r}(y, t, k), \qquad k \in \bar{U}_9 \cap \bar{U}_{10}, 
\end{align*}
which lead to the factorizations
$$\check{b}_u = \check{b}_{u,a}\check{b}_{u,r}, \qquad
\check{b}_l = \check{b}_{l,a}\check{b}_{l,r}, \qquad
\check{B}_u = \check{B}_{u,a}\check{B}_{u,r}, \qquad
\check{B}_l = \check{B}_{l,a}\check{B}_{l,r}.$$

{\bf Step 4: Deform again.}
Let the sets $\{U_j\}$ and $\{V_j\}$ be as in Figure \ref{contour3.pdf} and define $\mathfrak{m}(y,t,k)$ for $k$ near $\R$ by
\begin{align}\label{betweendef}
\mathfrak{m}(y,t,k) = \tilde{M}(y,t,k) G(y,t,k),
\end{align}
where
$$G(y,t,k) = \begin{cases}  
 b_{u,a}(y,t,k)^{-1}, & k \in U_1, \\
B_{l,a}(y,t,k)^{-1}, & k \in U_3, \\
B_{u,a}(y,t,k)^{-1}, & k \in U_4, \\
 b_{l,a}(y,t,k)^{-1}, & k \in U_6, \\
\check{b}_{u,a}(y,t,k)^{-1}, & k \in U_7, \\
\check{B}_{l,a}(y,t,k)^{-1}, & k \in U_9, \\
 \check{B}_{u,a}(y,t,k)^{-1}, & k \in U_{10}, \\
 \check{b}_{l,a}(y,t,k)^{-1}, & k \in U_{12}.
\end{cases}
$$
We define $\mathfrak{m}$ analogously near the lines $\omega \R$  and $\omega^2 \R$ and set $\mathfrak{m} = \tilde{M}$ elsewhere.
Let $\Gamma$ denote the jump contour for $\mathfrak{m}$; the part of $\Gamma$ near $\R$ is displayed in Figure \ref{contour3.pdf}. 
The function $\mathfrak{m}$ satisfies the jump condition $\mathfrak{m}_+ = \mathfrak{m}_- v$ on $\Gamma$, where, for $k$ near $\R_+$, 
\begin{subequations}\label{vdef}
\begin{align}
v(y,t,k) = \begin{cases}
b_{u,a}, & k \in \bar{U}_1 \cap \bar{V}_1, 
	\\
B_{l,a}, & k \in \bar{U}_2 \cap \bar{U}_3, 
	\\
B_{u,a}, & k \in \bar{U}_4 \cap \bar{U}_5,
	\\
b_{l,a}, & k \in \bar{U}_6 \cap \bar{V}_6,
 	\\
b_{u,a} \tilde{J}_{1,7}, & k \in \bar{U}_1 \cap \bar{U}_2,		
	\\
b_{l,a}  \tilde{J}_{6,18}, & k \in \bar{U}_5 \cap \bar{U}_6,
 	\\
b_{l,r}^{-1}b_{u,r}, & k \in \bar{U}_1 \cap \bar{U}_6,
 	\\
B_{u,r}^{-1}B_{l,r}, & k \in \bar{U}_3 \cap \bar{U}_4,
\end{cases}
\end{align}
and, for $k$ near $\R_-$,
\begin{align}
v(y,t,k) = \begin{cases}
\check{B}_{l,a}, & k \in \bar{U}_8 \cap \bar{U}_9, 
	\\
\check{b}_{u,a}, & k \in \bar{U}_7 \cap \bar{V}_3, 
	\\
\check{b}_{l,a}, & k \in \bar{U}_{12} \cap \bar{V}_4,
 	\\
\check{B}_{u,a}, & k \in \bar{U}_{10} \cap \bar{U}_{11},
	\\
\check{b}_{u,a} \tilde{J}_{3,12}, & k \in \bar{U}_7 \cap \bar{U}_8,		
	\\
\check{b}_{l,a}  \tilde{J}_{4,13}, & k \in \bar{U}_{11} \cap \bar{U}_{12},
	\\
\check{b}_{l,r}^{-1}\check{b}_{u,r}, & k \in \bar{U}_7 \cap \bar{U}_{12},
 	\\
\check{B}_{u,r}^{-1}\check{B}_{l,r}, & k \in \bar{U}_9 \cap \bar{U}_{10},
\end{cases}
\end{align}
\end{subequations}
with
\begin{align*}
& b_{u,a} \tilde{J}_{1,7} = 
\begin{pmatrix}
 1 & -\frac{\delta(\zeta, \omega k) \delta(\zeta, \omega^2k)}{\delta(\zeta, k)^2} (\overline{r_{1,a}(\bar{k})} - \overline{h(\bar{k})}) e^{-t\Phi} & 0 \\
0 & 1 & 0 \\
 0 & 0 & 1 
\end{pmatrix},
	\\
& b_{l,a}  \tilde{J}_{6,18}
= \begin{pmatrix} 1 & 0 & 0 \\
- \frac{\delta(\zeta, k)^2}{\delta(\zeta, \omega k) \delta(\zeta, \omega^2k)} (r_{1,a}(k) - h(k)) e^{t\Phi} & 1 & 0 \\
 0 & 0 & 1 \end{pmatrix},
	\\
& \check{b}_{u,a} \tilde{J}_{3,12} = 
\begin{pmatrix}
 1 & - \frac{\delta(\zeta, \omega k) \delta(\zeta, \omega^2k)}{\delta(\zeta, k)^2} (\overline{\check{r}_{1,a}(\bar{k})} - \overline{\check{h}(\bar{k})})  e^{-t\Phi} & 0 \\
0 & 1 & 0 \\
 0 & 0 & 1 
\end{pmatrix},
	\\
& \check{b}_{l,a}  \tilde{J}_{4,13}
= \begin{pmatrix}
 1 & 0 & 0 \\
- \frac{\delta(\zeta, k)^2}{\delta(\zeta, \omega k) \delta(\zeta, \omega^2k)} (\check{r}_{1,a}(k) - \check{h}(k)) e^{t\Phi} & 1 & 0 \\
 0 & 0 & 1 
\end{pmatrix}.
\end{align*}
The jump matrix $v$ obeys the symmetries (\ref{symmetriesa}) and (\ref{symmetriesb}). 
Lemma \ref{decompositionlemma} implies that 
$$G(y,t,\cdot) \in I + \dot{E}^3(\mathcal{U}) \cap E^\infty(\mathcal{U}),$$
where $\mathcal{U} = U_1 \cup U_3 \cup U_4 \cup U_6 \cup U_7 \cup U_9 \cup U_{10} \cup U_{12}$.
It follows from Lemma \ref{deformationlemma} that $\hat{N}(y,t,k)$ is a row-vector solution of the $L^3$-RH problem (\ref{RHNhat}) if and only if the function $n(y,t, k)$ defined by
\begin{equation}\label{ndef}
n(y, t, k) = (1,1,1) \mathfrak{m}(y,t,k), \qquad k \in \hat{\C} \setminus \Gamma
\end{equation}
is a row-vector solution of the $L^3$-RH problem
\begin{align}\label{RHn}
\begin{cases} n(y,t, \cdot) \in (1,1,1) + \dot{E}^3(\hat{\C} \setminus \Gamma), 
	\\
n_+(y,t,k) =  n_-(y,t,k) v(y,t,k) \quad \text{for a.e.} \ k \in \Gamma.
\end{cases} 
\end{align}

\bigskip
{\bf Step 5: Apply Theorem \ref{steepestdescentth}.}
Let  $\alpha \in [\frac{1}{3},1)$ and $c \in (0,3)$. Let $\epsilon = \frac{1-k_0}{2}$ (then $0<\epsilon <1$ and $k_0+\epsilon<1$). We claim that Theorem \ref{steepestdescentth} can be applied to the contour $\Gamma$ and jump matrix $v$ with 
\begin{align} \label{phizetazdef}
\mathcal{I} & = [c, 3), \quad \rho = \epsilon \sqrt{-i  F''(\tilde{k}_0) \tilde{k}'(k_0)^2}
=  \epsilon \frac{1+k_0^{-2}}{2} \sqrt{\frac{48 \tilde{k}_0 (3-4
   \tilde{k}_0^2)}{ (4 \tilde{k}_0^2+1)^3}}, 
   	\\ \nonumber
q(\zeta) & = \bigg(\frac{2}{k_0^{-1} + k_0}\bigg)^{2i\check{\nu}(\zeta)} 
\frac{e^{2\chi(\zeta, k_0)}r(k_0) }{\delta(\zeta, \omega k_0) \delta(\zeta, \omega^2k_0)} \bigg(\frac{\rho}{\epsilon}(k_0^{-1} - k_0)k_0\bigg)^{2i\nu(\zeta)},
\end{align}
\begin{align*}
 \check{q}(\zeta) &= \bigg(\frac{2}{k_0^{-1} + k_0}\bigg)^{2i\nu(\zeta)} 
\frac{\delta(\zeta, - \omega k_0) \delta(\zeta, -\omega^2 k_0)}{e^{2\chi(\zeta, -k_0)}}
\overline{\check{r}(-k_0)} 
		\\\nonumber
&\quad \times \bigg(\frac{\rho}{\epsilon}(k_0^{-1} - k_0)k_0\bigg)^{2i\check{\nu}(\zeta)},
	\\\nonumber
\nu(\zeta) & = -\frac{1}{2\pi} \log(1 - |r(k_0)|^2), \qquad
\check{\nu}(\zeta) = -\frac{1}{2\pi} \log(1 - |\check{r}(-k_0)|^2), 
	\\\nonumber
\phi(\zeta, z) & = \Phi\biggl(\zeta, k_0 - \frac{\epsilon z}{\rho}\biggr) = \overline{\Phi\biggl(\zeta, \overline{-k_0 + \frac{\epsilon z}{\rho}}\biggr)} 
	\\ \nonumber
&= -\frac{48i \tilde{k}_0^3}{(1 + 4 \tilde{k}_0^2)^2} + \frac{i}{2}z^2 + O(z^3), \qquad z \to 0.
\end{align*}
Indeed, by adding a number of arcs on which $v = I$, we can ensure that $\Gamma$ is a Carleson jump contour which satisfies $(\Gamma1)$-$(\Gamma4)$. 

\begin{lemma}
The $3 \times 3$-matrix valued function $w = v - I$ satisfies (\ref{winL1L2Linf}) and (\ref{wL12infty}). 
\end{lemma}
\begin{proof}
This follows from the decay properties of $e^{\pm t \Phi}$. Indeed, for $k \in \Gamma'$ near $\R$, the decay is a consequence of the expressions for the jump matrix $v$ given in (\ref{vdef}) and the estimates in Lemma \ref{decompositionlemma}. By symmetry, it follows that $w$ has decay also near $\omega \R$ and $\omega^2 \R$. On the other hand, we already noted that $w$ has exponential decay on the circles $\cup_{j=1}^6 \partial B_j$, where $B_j$ denotes the disk centered at $K_j$ along which the $E_n$'s and $E_{n+18}$'s meet (see Step 2 above). Thus, it only remains to verify that $w$ is small on $\Gamma' \cap (\cup_{j=1}^6 B_j)$. We will show that $w$ is small on $\bar{E}_{19} \cap \bar{E}_{25}$; the other parts of $\Gamma' \cap (\cup_{j=1}^6 B_j)$ can be handled in a similar way. 

In view of (\ref{Jformula2}), we have
$$v(y,t,k) =  \begin{pmatrix}
1 & \frac{\delta(\zeta, \omega k) \delta(\zeta, \omega^2 k)}{\delta(\zeta, k)^2} f_1(k) 
e^{-t\Phi(\zeta,k)} & 0 \\
0 & 1 & 0 \\
0 & 0 & 1 
\end{pmatrix},\quad  k \in \bar{E}_{19} \cap \bar{E}_{25},$$
where the functions $\frac{\delta(\zeta, \omega k) \delta(\zeta, \omega^2 k)}{\delta(\zeta, k)^2}$ and $f_1(k)$ are bounded on $\bar{E}_{19} \cap \bar{E}_{25}$.
In particular,
$$|w(y,t,k)| \leq C e^{-t \re \Phi(\zeta,k)}, \qquad \zeta \in \mathcal{I}, \ k \in \bar{E}_{19} \cap \bar{E}_{25}.$$
It is therefore enough to show that there exists a constant $c_1 > 0$ independent of $\zeta, k$ such that
\begin{align}\label{rePhiE19E25}
\re \Phi(\zeta,k) > c_1 > 0, \qquad \zeta \in \mathcal{I}, \ k \in \bar{E}_{19} \cap \bar{E}_{25}.
\end{align}
In order to prove (\ref{rePhiE19E25}), we note that the curve $\bar{E}_{19} \cap \bar{E}_{25} \subset \bar{D}_1 \cap \bar{D}_7$ is determined by the equation $\re(z_1 - z_2) = 0$. Hence the definition (\ref{Phizetakdef}) of $\Phi(\zeta,k)$ yields
\begin{align*}
\re \Phi(\zeta, k) & = \re (l_2 - l_1)\zeta = \frac{k_2 (k_1^2+k_2^2+1) \zeta}{k_1^2+k_2^2}, \qquad k \in \bar{D}_1 \cap \bar{D}_7.
\end{align*}
Since $k_2 = \im k$ is positive and bounded away from $0$ on $\bar{E}_{19} \cap \bar{E}_{25}$ and $\mathcal{I} = [c, 3)$ with $c > 0$, the inequality (\ref{rePhiE19E25}) follows.
\end{proof}

The definition (\ref{vdef}) of $v$ implies that (\ref{vjdef}) and (\ref{smallcrossjump}) are satisfied with
\begin{align*}
\begin{cases}
R_1(\zeta, t, z) = \frac{\delta(\zeta,k)^2}{\delta(\zeta,\omega k) \delta(\zeta,\omega^2k)} (r_{1,a}(k) - h(k))  z^{2i\nu(\zeta)},
	\\
R_2(\zeta, t, z) = \frac{\delta(\zeta,\omega k) \delta(\zeta,\omega^2k)}{\delta(\zeta,k)^2} \overline{r_{2,a}(\bar{k})} z^{-2i\nu(\zeta)},
	\\
R_3(\zeta, t, z) =  \frac{\delta(\zeta,k)^2}{\delta(\zeta,\omega k) \delta(\zeta,\omega^2k)} r_{2,a}(k) z^{2i\nu(\zeta)},
	\\
R_4(\zeta, t, z) = \frac{\delta(\zeta,\omega k) \delta(\zeta,\omega^2k)}{\delta(\zeta,k)^2}  (\overline{r_{1,a}(\bar{k})} - \overline{h(\bar{k})}) z^{-2i\nu(\zeta)},
\end{cases}
\end{align*}
where $k$ and $z$ are related by $k = k_0 - \frac{\epsilon z}{\rho}$, and
\begin{align*}
\begin{cases}
\check{R}_1(\zeta, t, z) =  \frac{\delta(\zeta,\omega k) \delta(\zeta,\omega^2k)}{\delta(\zeta,k)^2} (\overline{\check{r}_{1,a}(\bar{k})} - \overline{\check{h}(\bar{k})})  z^{2i\check{\nu}(\zeta)},
	\\
\check{R}_2(\zeta, t, z) = \frac{\delta(\zeta,k)^2}{\delta(\zeta,\omega k) \delta(\zeta,\omega^2k)} \check{r}_{2,a}(k) z^{-2i\check{\nu}(\zeta)},
	\\
\check{R}_3(\zeta, t, z) = \frac{\delta(\zeta,\omega k) \delta(\zeta,\omega^2k)}{\delta(\zeta,k)^2} \overline{\check{r}_{2,a}(\bar{k})}  z^{2i\check{\nu}(\zeta)},
	\\
\check{R}_4(\zeta, t, z) = \frac{\delta(\zeta,k)^2}{\delta(\zeta,\omega k) \delta(\zeta,\omega^2k)} (\check{r}_{1,a}(k) - \check{h}(k)) z^{-2i\check{\nu}(\zeta)},
\end{cases}
\end{align*}
where $k$ and $z$ are related by $k = -k_0 + \frac{\epsilon z}{\rho}$. 
The definition (\ref{phizetazdef}) of $\phi(\zeta, z)$ shows that (\ref{phiassumptions}) and (\ref{rephiestimates}) hold. 
The symmetry $\delta(\zeta,k) = 1/\overline{\delta(\zeta,\bar{k})}$ implies that $|\delta(\zeta,\omega k_0) \delta(\zeta,\omega^2 k_0)| = 1$. Hence $|q(\zeta)| = |r(k_0)|$ and $|\check{q}(\zeta)| = |\check{r}(-k_0)|$; this yields (\ref{nuchecknudef}).
To establish (\ref{Lipschitzconditions}), we note that if $k = k_0 - \frac{\epsilon z}{\rho}$, then
\begin{align*}
 R_1(\zeta, t, z)  = &\; \big((k_0^{-1} - k)(k^{-1} - k_0)kk_0\big)^{i\nu} \bigg(\frac{(k_0 + k)(k_0^{-1} + k^{-1})}{(k_0^{-1} + k)(k_0 + k^{-1})}\bigg)^{i\check{\nu}} 
	\\
& \times \frac{e^{2\chi(\zeta, k)}}{\delta(\zeta,\omega k) \delta(\zeta,\omega^2k)} (r_{1,a}(k) - h(k)) \bigg(\frac{\rho}{\epsilon}\bigg)^{2i\nu(\zeta)}.
\end{align*}
Since $r_{1,a}(k_0) - h(k_0) = r(k_0)$, this shows that $q(\zeta) = R_1(\zeta, t, 0)$. Similarly, if $k = -k_0 + \frac{\epsilon z}{\rho}$, then
\begin{align*}
\check{R}_1(\zeta, t, z) 
%= \bigg(\frac{(k_0^{-1} - k)(k_0 - k^{-1})}{(k_0 - k)(k_0^{-1} - k^{-1})}\bigg)^{-i\nu} \bigg(\frac{(k_0 + k)(k_0^{-1} + k^{-1})}{(k_0^{-1} + k)(k_0 + k^{-1})}\bigg)^{-i\check{\nu}}  \frac{\delta(\zeta,\omega k) \delta(\zeta,\omega^2k)}{e^{2\chi(\zeta, k)}} \overline{\check{r}(\bar{k})} \bigg(\frac{\rho}{\epsilon}(k_0 + k)\bigg)^{2i\check{\nu}(\zeta)}
%	\\
= &\; \bigg(\frac{(k_0^{-1} - k)(k_0 - k^{-1})}{(k_0 - k)(k_0^{-1} - k^{-1})}\bigg)^{-i\nu} 
\bigg(\frac{1}{(k_0^{-1} + k)(k_0 + k^{-1})kk_0}\bigg)^{-i\check{\nu}}  
	\\
& \times \frac{\delta(\zeta,\omega k) \delta(\zeta,\omega^2k)}{e^{2\chi(\zeta, k)}} (\overline{\check{r}_{1,a}(\bar{k})} - \overline{\check{h}(\bar{k})})\bigg(\frac{\rho}{\epsilon}\bigg)^{2i\check{\nu}(\zeta)}.
\end{align*}
Since $(\overline{\check{r}_{1,a}(-k_0)} - \overline{\check{h}(-k_0)}) = \overline{\check{r}(-k_0)}$, this shows that $\check{q}(\zeta) = \check{R}_1(\zeta, t, 0)$. The inequalities in (\ref{Lipschitzconditions}) for $R_1$ and $\check{R}_1$ are now a consequence of standard estimates, cf. \cite{DZ1993}. The inequalities in (\ref{Lipschitzconditions}) for $\{R_j, \check{R}_j\}_2^4$ are proved in a similar way.
This shows that the conditions of Theorem \ref{steepestdescentth} are satisfied. 

The conclusion (\ref{limlm12}) of Theorem \ref{steepestdescentth} implies that the solution $n(y,t,k)$ of the $L^3$-RH problem (\ref{RHn}) satisfies
 \begin{align}\nonumber
n(y,t,K_1)
= &\; (1,1,1) + \frac{2\epsilon}{k_0 \sqrt{\tau}} \re\big(\mathcal{F}_1 \beta - \bar{\mathcal{F}}_2\check{\beta}, 
\mathcal{F}_3 \beta - \bar{\mathcal{F}}_3\check{\beta}, \mathcal{F}_2 \beta - \bar{\mathcal{F}}_1\check{\beta}\big) 
	\\ \label{limn}
& + O\bigl(\epsilon\tau^{-\frac{1+\alpha}{2}} \bigr), \qquad \tau = t\rho^2 \to \infty, \  \zeta \in \mathcal{I},
\end{align}
where the error term is uniform with respect to $\zeta \in \mathcal{I}$ and the functions $\{\mathcal{F}_j\}_1^3$, $\beta$, $\check{\beta}$ are defined in (\ref{Xjdef}) and (\ref{betadef}).

\begin{rema}
In general, the solution $m$ of the $L^3$-RH problem (\ref{RHm}) featured in Theorem \ref{steepestdescentth} is different than the function $\mathfrak{m}$ used in this section; the former is regular at the points $\{\varkappa_j\}_1^6$ whereas the latter, in general, is singular at these points. However, by Lemma \ref{rowvectorlemma}, this discrepancy disappears when premultiplying by $(1,1,1)$; hence the row vector solution $n$ of the $L^3$-RH problem (\ref{RHn}) satisfies (\ref{limn}).
\end{rema}

\bigskip
{\bf Step 6: Find $q(x,t)$ and $u(x,t)$.}
For $k \in E_{25}$ near $K_1$ we have 
$$\mathfrak{m} = M\Delta =
P(k)^{-1}\mathcal{D}(x,t)^{-1} P(k) \Psi_7 e^{(x-y + \nu_0) \mathcal{L}} \Delta.$$
Using the identity
$$(1,1,1)P(k)^{-1} \mathcal{D}(x,t)^{-1} P(k) = q(x,t)(1,1,1),$$ 
this gives
\begin{align}\nonumber
n(y, t, k)
 = &\; q(x,t)(1,1,1) \Psi_7(x, t, k) e^{(x-y + \nu_0) \mathcal{L}} \Delta(\zeta, k).
\end{align}
Since $\Psi_n(x,t,K_1) = I$, evaluation of this equation at $k = K_1$ yields
\begin{align*}
n(y, t, K_1)
& = q(x,t) (1,1,1) e^{(x-y + \nu_0)\mathcal{L}(K_1)} \Delta(\zeta, K_1)
 	\\
& = q(x,t)  (e^{y-x -  \nu_0}\Delta_{11}(\zeta, K_1) , \Delta_{22}(\zeta, K_1), e^{x-y+ \nu_0}\Delta_{33}(\zeta, K_1)).
\end{align*}
Hence, by (\ref{limn}),\footnote{All error terms of the form $O(\cdot)$ in the remainder of the proof are uniform with respect to $\zeta \in \mathcal{I}$.}
$$\frac{\Delta_{11}(\zeta, K_1)\Delta_{33}(\zeta, K_1)}{\Delta_{22}^2(\zeta, K_1)} = \frac{n_1(\zeta, t, K_1) n_3(\zeta, t, K_1)}{n_2^2(\zeta, t, K_1)} = 1 + O(\epsilon \tau^{-1/2})$$
as $\tau \to \infty$. Fixing $\zeta \in \mathcal{I}$ on the left-hand side of this equation and letting $t \to \infty$, we deduce that
$\Delta_{11}(\zeta, K_1)\Delta_{33}(\zeta, K_1) = \Delta_{22}^2(\zeta, K_1)$ for $\zeta \in \mathcal{I}$. Proceeding as in the proof of Proposition 4.2 of \cite{BS2013}, we infer that $\Delta_{22}(\zeta, K_1) = 1$ and $|r(k_0)| = |\check{r}(-k_0)|$ for all  $\zeta \in \mathcal{I}$. It follows that $\nu(\zeta) = \check{\nu}(\zeta)$ for $\zeta \in \mathcal{I}$. The equation
$$\delta(\zeta, k) = e^{\frac{1}{4\pi i} \int_{k_0}^{k_0^{-1}} \log(1 - |r(s)|^2)\big(\frac{1}{s-k} - \frac{1}{s + k} - \frac{1}{s-k^{-1}} + \frac{1}{s + k^{-1}}\big)ds},$$
now shows that $\delta(\zeta, k) = \delta(\zeta, -k)^{-1}$ and $\chi(\zeta, k) = -\chi(\zeta, -k)$. 
In particular, $q(\zeta) = \check{q}(\zeta) e^{i(\arg \check{r}(-k_0) + \arg r(k_0))}$. 

A computation shows that
\begin{align}\label{argq}
 \arg q(\zeta) = &\; 2 \nu \log\bigg(\frac{2}{k_0 + k_0^{-1}}\bigg) + \chi_0(\zeta) + \arg r(k_0) 
 	\\\nonumber
 &+ 2 \nu \log\bigg(\frac{\rho}{\epsilon}(k_0^{-1} - k_0)k_0\bigg)
	\\\nonumber
& - \frac{\nu}{2} \log\bigg|\frac{(k_0^{-1} - k)(k_0 - k^{-1})}{(k_0 - k)(k_0^{-1} - k^{-1})}\frac{(k_0 + k)(k_0^{-1} + k^{-1})}{(k_0^{-1} + k)(k_0 + k^{-1})}\bigg|_{k = \omega k_0}
	\\\nonumber
 &  - \frac{\nu}{2} \log\bigg|\frac{(k_0^{-1} - k)(k_0 - k^{-1})}{(k_0 - k)(k_0^{-1} - k^{-1})}\frac{(k_0 + k)(k_0^{-1} + k^{-1})}{(k_0^{-1} + k)(k_0 + k^{-1})}\bigg|_{k = \omega^2 k_0}
	\\ \nonumber
= &\; \chi_0(\zeta) + \arg r(k_0) - \nu \log Y,
\end{align}
where 
$$\chi_0(\zeta) = \im(2\chi(\zeta, k_0) - \chi(\zeta, \omega k_0) - \chi(\zeta, \omega^2 k_0))$$
and the function $Y = Y(\zeta)$ is defined by
$$Y(\zeta) = \frac{(4 \tilde{k}_0^2+1)^2 (4
   \tilde{k}_0^2+3)}{576 \tilde{k}_0^3 (3 - 4\tilde{k}_0^2)}.$$
Equations (\ref{limn}) and (\ref{argq}) yield
\begin{align} \nonumber
n_1(y, t, K_1) = &\;1 + \frac{d_1}{\sqrt{t}}
\re \Big(  \big(\mathcal{F}_1 e^{-i\arg r(k_0)} 
-  \bar{\mathcal{F}}_2  e^{i \arg \check{r}(-k_0)}\big)e^{i(d_2 t - \nu \log t + d_3)}\Big) 
	\\\label{n1expression}
& + O(\epsilon \tau^{-\frac{1+\alpha}{2}}), \qquad \tau \to \infty, \  \zeta \in \mathcal{I},
\end{align}
and
\begin{align}\nonumber
q(x,t) = &\; n_2(y, t, K_1) 
	\\\nonumber
=&\;  1  + \frac{d_1}{\sqrt{t} }\re\Big( \big(\mathcal{F}_3 e^{-i\arg r(k_0)}
- \bar{\mathcal{F}}_3e^{i \arg \check{r}(-k_0)}\big)e^{i(d_2 t - \nu \log t + d_3)}\Big)
	\\ \label{qformula1}
&  + O(\epsilon \tau^{-\frac{1+\alpha}{2}}), \qquad \tau \to \infty, \  \zeta \in \mathcal{I},
 \end{align}
where the functions $d_j = d_j(\zeta)$, $j = 1, 2, 3$, are defined by
\begin{align*}
&d_1 = \frac{2\epsilon\sqrt{\nu}}{k_0 \rho}, \qquad d_2 = \frac{48 \tilde{k}_0^3}{(1 + 4 \tilde{k}_0^2)^2}, \qquad d_3 = \frac{\pi}{4} - \chi_0 + \nu \log Y + \arg \Gamma(i\nu). 
\end{align*}

It can be seen from the proof of Theorem \ref{mainth} that thanks to the uniform decay and smooth dependence of the jump matrix $v$ on $t$, the asymptotic formula (\ref{limlm12}) can be differentiated in time without affecting the error term. Hence equations (\ref{n1expression}) and (\ref{qformula1}) together with the fact that 
\begin{align}\label{yminusx}
y - x = \log\bigg(\frac{n_1(y, t, K_1)}{n_2(y, t, K_1)}\bigg) - \log \Delta_{11}(\zeta, K_1) + \nu_0,
\end{align}
yield
\begin{align} \nonumber
 u(x,t) = &\; \frac{\partial}{\partial t}\bigg|_{\text{$y$ fixed}} (x - y)
  = \frac{\partial}{\partial t}\bigg|_{\text{$y$ fixed}} \frac{d_1}{\sqrt{t}}
\re\Big( f e^{i(d_2 t - \nu \log t + d_3)}\Big) 
 	\\ \label{uform1}
& + \frac{\partial}{\partial t}\bigg|_{\text{$y$ fixed}} \log \Delta_{11}(\zeta, K_1)
+ O(\epsilon \tau^{-\frac{1+\alpha}{2}}),
\end{align}
where 
$$f(\zeta) = (\mathcal{F}_3 - \mathcal{F}_1) e^{-i\arg r(k_0)} - ( \bar{\mathcal{F}}_3 -  \bar{\mathcal{F}}_2)e^{i \arg \check{r}(-k_0)}.$$
As $\zeta \to 3^-$, we have the expansions
\begin{align*}
& \epsilon = \frac{1}{12}(3 - \zeta)^{\frac{1}{2}} + O(3-\zeta), &&
k_0 = 1 - \frac{(3 - \zeta)^{\frac{1}{2}}}{6}  + O(3 - \zeta),
	\\
& \tilde{k}_0 =  \frac{(3 - \zeta)^{\frac{1}{2}}}{6}  + O((3 - \zeta)^{\frac{3}{2}}),
&&
 \rho =  \frac{(3-\zeta)^{\frac{3}{4}}}{\sqrt{6}}  + O((3-\zeta)^{\frac{5}{4}}),
	\\
& d_1 = \frac{\sqrt{\nu}(3-\zeta)^{-\frac{1}{4}}}{\sqrt{6}}  + O((3-\zeta)^{\frac{3}{4}}), &&
 d_2 = \frac{2(3-\zeta)^{\frac{3}{2}}}{9}  + O((3-\zeta)^{\frac{5}{2}}), 
 	\\
& \mathcal{F}_1 = \frac{3 + \sqrt{3}}{2} + O((3 - \zeta)^{\frac{1}{2}}), &&
 \mathcal{F}_2 = \frac{-3 + \sqrt{3}}{2} + O((3 - \zeta)^{\frac{1}{2}}),  
	\\
&\mathcal{F}_3 = -\sqrt{3} + O((3 - \zeta)^{\frac{1}{2}}), &&
Y = \frac{3(3 - \zeta)^{-\frac{3}{2}}}{8}+ O((3 - \zeta)^{-\frac{1}{2}}).
\end{align*}
We deduce that there exist constants $c_1, c_2> 0$ such that
\begin{align}\nonumber
& c_1 \epsilon^{3/2} < \rho(\zeta) < c_2 \epsilon^{3/2},
\qquad \nu'(\zeta) < c_2 \epsilon^{-1}, \qquad
f'(\zeta) < c_2 \epsilon^{-1}, 
	\\ \label{epsilonestimates}
& d_1'(\zeta) < c_2 \epsilon^{-\frac{5}{2}}, \qquad
d_3'(\zeta) < c_2 \epsilon^{-2},
\end{align}
for all $\zeta \in \mathcal{I}$. Moreover, since
\begin{align}\label{Delta11}
\Delta_{11}(\zeta, K_1) = e^{\frac{3}{2\pi} \int_{k_0}^{k_0^{-1}} \log(1 - |r(s)|^2) \frac{1 + s^4}{1 + s^6} ds},
\end{align}
we obtain the estimates
\begin{align*}
& c_1 \epsilon < |1- \Delta_{11}(\zeta, K_1)| < c_2 \epsilon, \qquad
c_1 \epsilon^{-1}  < \big|\partial_\zeta \Delta_{11}(\zeta, K_1)\big| < c_2 \epsilon^{-1},
\end{align*}
which show that
\begin{align}\label{logDelta11}
\bigg|\frac{\partial}{\partial t}\bigg|_{\text{$y$ fixed}} \log \Delta_{11}(\zeta, K_1)\bigg| 
= \bigg|\frac{\zeta \partial_\zeta\Delta_{11}(\zeta, K_1)}{t \Delta_{11}(\zeta, K_1)}\bigg|
< C\epsilon^{-1}t^{-1} = O( \epsilon^2 \tau^{-1})
\end{align}
as $\tau \to \infty$. The above and the following estimates are uniformly valid for $\zeta \in \mathcal{I}$. Using the estimates (\ref{epsilonestimates}) and (\ref{logDelta11}) together with the identities 
$$\frac{\partial \zeta}{\partial t}\bigg|_{\text{$y$ fixed}} = -\frac{\zeta}{t} \quad \text{and} \quad d_2 + t \frac{\partial d_2}{\partial t} = \frac{6\tilde{k}_0}{1 + 4 \tilde{k}_0^2},$$
equation (\ref{uform1}) yields, as $\tau \to \infty$,
\begin{align} \label{uformula1}
 u(x,t) = \frac{d_1}{\sqrt{t}}\frac{6\tilde{k}_0}{1 + 4 \tilde{k}_0^2}
\re\Big(i f e^{i(d_2 t - \nu \log t + d_3)}\Big) 
+ O(\epsilon \tau^{-\frac{1+\alpha}{2}}).
\end{align}
Since $\frac{d}{dx} = q \frac{d}{dy}$ we find $u_{xx} =  u_{yy} + O(\epsilon \tau^{-\frac{1+\alpha}{2}})$. Using that $\frac{1}{1 + 4 \tilde{k}_0^2}(1 + (\frac{dd_2}{d\zeta})^2) 
= 1$, we conclude that
\begin{align}\label{uuxx1}
u - u_{xx} + 1= & \;1 + \frac{d_1}{\sqrt{t}} 6\tilde{k}_0
\re\Big( i f e^{i(d_2 t - \nu \log t + d_3)}\Big) + O(\epsilon \tau^{-\frac{1+\alpha}{2}}).
\end{align}
Substituting (\ref{qformula1}) and (\ref{uuxx1}) into the relation $q^3 = u - u_{xx} + 1$, the terms of order $O(t^{-1/2})$ yield
$$\big(6\tilde{k}_0 i   (\mathcal{F}_3 - \mathcal{F}_1)  - 3 \mathcal{F}_3\big) - \big(6\tilde{k}_0 i ( \bar{\mathcal{F}}_3 -  \bar{\mathcal{F}}_2) - 3 \bar{\mathcal{F}}_3\big)e^{i \arg \check{r}(-k_0) + i\arg r(k_0)} = 0,$$
that is,
\begin{align}\label{eiarg}
e^{i \arg \check{r}(-k_0) + i\arg r(k_0)}= \frac{6\tilde{k}_0 i   (\mathcal{F}_3 - \mathcal{F}_1)  - 3 \mathcal{F}_3}{6\tilde{k}_0 i ( \bar{\mathcal{F}}_3 -  \bar{\mathcal{F}}_2) - 3 \bar{\mathcal{F}}_3} = \frac{1 - \omega k_0^2}{k_0^2 - \omega}.
\end{align}
Using (\ref{eiarg}) and the identity
$$\mathcal{F}_3 - \bar{\mathcal{F}}_3 \frac{1 - \omega k_0^2}{k_0^2 - \omega}
= \frac{2 \sqrt{3} \tilde{k}_0 \sqrt{4 \tilde{k}_0^2 + 3}}{4 \tilde{k}_0^2 + 1} e^{-i \arctan\big(\sqrt{3}\frac{1 + k_0^2}{1 - k_0^2}\big)}$$
in (\ref{qformula1}), we find, as $\tau \to \infty$,
\begin{align}\label{qfinalzeta}
q(x,t) = \; & 1 + \frac{c_1(\zeta)}{\sqrt{t}}
\cos \big(c_2(\zeta) t - \nu(\zeta) \log(t) + c_3(\zeta)\big) + O(\epsilon \tau^{-\frac{1+\alpha}{2}}),
\end{align}
where 
\begin{align*}
c_1(\zeta) &= \frac{1 - k_0^2}{1 + k_0^2} \sqrt{\frac{(3 + 4\tilde{k}_0^2)(1 + 4 \tilde{k}_0^2)\nu}{3\tilde{k}_0 - 4\tilde{k}_0^3}}, \qquad c_2(\zeta) = \frac{48 \tilde{k}_0^3}{(1 + 4 \tilde{k}_0^2)^2},
	\\
 c_3(\zeta) &= d_3 - \arg r(k_0) - \arctan\bigg(\sqrt{3}\frac{1 + k_0^2}{1 - k_0^2}\bigg)
	\\
& =  \frac{\pi}{4} - \chi_0(\zeta) + \nu \log Y 
+ \arg \Gamma(i\nu) - \arg r(k_0) - \arctan\bigg(\sqrt{3}\frac{1 + k_0^2}{1 - k_0^2}\bigg).
\end{align*}
	
Similarly, using (\ref{eiarg}) in (\ref{uformula1}), we find
\begin{align*}
u(x,t) = &\; \frac{6\tilde{k}_0 d_1}{(1 + 4 \tilde{k}_0^2) \sqrt{t}}
\re\bigg\{i \bigg(\mathcal{F}_3 - \mathcal{F}_1 - ( \bar{\mathcal{F}}_3 -  \bar{\mathcal{F}}_2)\frac{1 - \omega k_0^2}{k_0^2 - \omega}\bigg)
		\\
&\times e^{i(d_2 t - \nu \ln t + d_3 -\arg r(k_0))}\bigg\} + O(\epsilon \tau^{-\frac{1+\alpha}{2}}).
\end{align*}
In view of the identity
$$\mathcal{F}_3 - \mathcal{F}_1 - ( \bar{\mathcal{F}}_3 -  \bar{\mathcal{F}}_2)\frac{1 - \omega k_0^2}{k_0^2 - \omega}
= -\frac{i\sqrt{3}\sqrt{4 \tilde{k}_0^2 + 3}}{4 \tilde{k}_0^2 + 1} e^{-i\arctan\big(\sqrt{3}\frac{ 1 + k_0^2}{1 - k_0^2}\big)},$$
this yields
\begin{align}\label{ufinalzeta}
u(x,t) =&\;  \frac{3 c_1(\zeta)}{(1 + 4 \tilde{k}_0^2(\zeta)) \sqrt{t}}
\cos\big(c_2(\zeta) t - \nu(\zeta) \log t + c_3(\zeta)\big) 
	\\\nonumber
&+ O(\epsilon \tau^{-\frac{1+\alpha}{2}})
\end{align}
as $\tau \to \infty$ uniformly with respect to $\zeta \in \mathcal{I}$.

\bigskip
{\bf Step 7: Replace $\zeta$ with $\xi$.}	
In the last step of the proof, we show that, up to a phase shift, $\zeta = y/t$ can be replaced with $\xi = x/t$ in the asymptotic formulas (\ref{qfinalzeta}) and (\ref{ufinalzeta}) without affecting the error term. For clarity, we reinsert the dependence on $\zeta$ of the functions $k_0(\zeta)$, $\tilde{k}_0(\zeta)$, and $\epsilon(\zeta)$. 

By (\ref{n1expression}), (\ref{qformula1}), and (\ref{yminusx}), we have 
$$y-x= - \log \Delta_{11}(\zeta, K_1) + \nu_0 + O(\epsilon(\zeta) \tau^{-1/2}).$$ 
Hence, as $\tau \to \infty$,
$$\zeta - \xi = \frac{c(\zeta)}{t} + O(\epsilon(\zeta) t^{-1} \tau^{-\frac{1}{2}}) = O((3-\zeta)^{\frac{1}{2}} t^{-1})
= O((3-\xi)^{\frac{1}{2}} t^{-1}).$$
The asymptotic sector $\{\zeta \in \mathcal{I}, \tau \to \infty\}$ is equivalent to $\{\zeta \in \mathcal{I}, t(3 - \zeta)^{\frac{3}{2}}  \to \infty\}$ and hence also to $\{\xi \in \mathcal{I}, t(3 - \xi)^{\frac{3}{2}} \to \infty\}$.
If $\alpha\in \R$ and $g(\zeta)$ is a smooth function such that $|g'(\zeta)| < C (3 - \zeta)^{\alpha}$ for all $\zeta \in \mathcal{I}$, then
\begin{align*}
 |g(\xi) - g(\zeta)| & = \bigg|\int_\zeta^{\xi} g'(\eta) d\eta\bigg|
< C|(3 - \xi)^{\alpha + 1} - (3 - \zeta)^{\alpha + 1}|
	\\
& < C |\xi - \zeta| \max\big\{(3 - \zeta)^\alpha, (3 - \xi)^\alpha\big\}
< C |\xi - \zeta| (3 - \xi)^\alpha
	\\
& = O\big((3 - \xi)^{\alpha + \frac{1}{2}} t^{-1}\big), \qquad \tau \to \infty, \  \zeta \in \mathcal{I}.
\end{align*}
The estimates 
\begin{align*}
& |c_1'(\zeta)| < C (3-\zeta)^{-\frac{3}{4}}, \quad
|c_2'(\zeta)| < C (3-\zeta)^{\frac{1}{2}}, \quad
|c_3'(\zeta)| < C (3-\zeta)^{-1},
	\\
& |\nu'(\zeta)| < C (3-\zeta)^{-\frac{1}{2}},\quad
|\epsilon'(\zeta)| < C (3-\zeta)^{-\frac{1}{2}},\quad
|\tilde{k}_0'(\zeta)| < C (3-\zeta)^{-\frac{1}{2}},
\end{align*}
therefore imply
\begin{align}\nonumber
& |c_1(\xi) - c_1(\zeta)| = O\big((3 - \xi)^{-\frac{1}{4}} t^{-1}\big), \qquad
|c_2(\xi) - c_2(\zeta)| = O\big((3 - \xi) t^{-1}\big), 
	\\\nonumber
& |c_3(\xi) - c_3(\zeta)| = O\big((3 - \xi)^{-\frac{1}{2}} t^{-1}\big), \qquad
|\nu(\xi) - \nu(\zeta)| = O\big(t^{-1}\big), 
	\\\label{xizetaestimates}
& |\epsilon(\xi) - \epsilon(\zeta)| = O\big(t^{-1}\big), \qquad
|\tilde{k}_0(\xi) - \tilde{k}_0(\zeta)| = O\big(t^{-1}\big).
\end{align}
On the other hand, the identity
$$c_2'(\zeta) = \frac{dc_2}{d\tilde{k}_0}\bigg/ \frac{d\zeta}{d\tilde{k}_0} = -2\tilde{k}_0(\zeta), \qquad \zeta \in \mathcal{I},$$
and the estimate 
$$|c_2''(\zeta)| < C (3-\zeta)^{-\frac{1}{2}}, \qquad \zeta \in \mathcal{I},$$
imply
\begin{align}\nonumber
c_2(\zeta) - c_2(\xi) & = -2\tilde{k}_0(\xi) (\zeta - \xi) + O((3 - \xi)^{-\frac{1}{2}} (\zeta - \xi)^2) 
	\\ \label{c2zetaxi}
& = 2\tilde{k}_0(\xi) t^{-1}  \log \Delta_{11}(\zeta, K_1) + O((3-\xi)^{\frac{1}{4}} t^{-\frac{3}{2}}).
\end{align}
Since
\begin{align*}
\chi(\zeta, k) = &\; \frac{1}{4\pi i} \int_{k_0}^{\frac{1}{k_0}}  \log\bigg(\frac{1 - |r(s)|^2}{1 - |r(k_0)|^2}\bigg) 	
	\\
&\times \bigg(\frac{1}{s-k} - \frac{1}{s + k} - \frac{1}{s-k^{-1}} + \frac{1}{s + k^{-1}}\bigg) ds,
\end{align*}
we see that $\chi_0(\xi)$ can be expressed as in (\ref{chi0def}).
Employing equations (\ref{Delta11}), (\ref{xizetaestimates}), and (\ref{c2zetaxi}), the asymptotic formulas (\ref{qfinalxi}) and (\ref{ufinalxi}) follow from (\ref{qfinalzeta}) and (\ref{ufinalzeta}), respectively.

\begin{rema}
Substituting the asymptotic formula (\ref{ufinalxi}) for $u(x,t)$ into (\ref{DP}), we can verify explicitly that the DP equation is satisfied to leading order in the similarity region. Indeed, by (\ref{ufinalxi}), the nonlinear terms in (\ref{DP}) are easily seen to be of order $O(\epsilon \tau^{-\frac{1+\alpha}{2}})$ as $\tau \to \infty$, and the linear terms satisfy
$$(u - u_{xx})_t + 3u_x = O(\epsilon \tau^{-\frac{1+\alpha}{2}})$$
as a consequence of the identity
$$\bigg(b_2 - \xi \frac{\partial b_2}{\partial \xi}\bigg)\bigg(1 + \bigg(\frac{\partial b_2}{\partial \xi}\bigg)^2\bigg) + 3 \frac{db_2}{d\xi} = 0.$$
\end{rema}

\begin{rema}
The main contributions to the asymptotic formula (\ref{ufinalxi}) come from the critical points $\omega^j k_0^{\pm 1}$ located on the lines $\omega^j \R$, $j = 0,1,2$. On the other hand, the same is true for the asymptotics of the solution of the whole line problem \cite{BS2013}. Therefore, the structure of the asymptotics for the whole line and half-line problems is the same, the only difference being in the determination of $r(k)$ which, in turn, determines $\nu$ ($h_0$ in the notation of \cite{BS2013}). 

Following \cite{BS2013}, these contributions can be determined by parametrizing the neighborhood of $k = k_0$ using the rescaled spectral parameter $z = \frac{\rho_{\rm new}}{\epsilon}(\tilde{k} - \tilde{k}_0)$,
where (recall that $\Phi(\zeta, k) = \tilde{\Phi}(\zeta, \tilde{k}(k))$)
$$\phi_{\rm new}(\zeta, z) = \tilde{\Phi}(\zeta, \tilde{k}_0 + \frac{\epsilon}{\rho_{\rm new}}z)$$
$$\rho_{\rm new} = \sqrt{-i\epsilon^2 \partial_{\tilde k}^2\tilde\Phi(\zeta,\tilde k_0)}
= \epsilon\sqrt{ \frac{48 \tilde{k}_0 (3 - 4\tilde{k}_0^2)}{(1 + 4\tilde{k}_0^2)^3}},$$
and we denote quantities defined using this rescaled spectral parameter by the subscript/superscript $new$.
It follows that (\ref{vjdef}) and (\ref{smallcrossjump}) are satisfied with
\begin{align*}
\begin{cases}
R_1^{\rm new}(\zeta, t, z) = \frac{\delta(\zeta, k)^2}{\delta(\zeta, \omega k) \delta(\zeta, \omega^2k)} r(k)  z^{2i\nu(\zeta)},
	\\
R_2^{\rm new}(\zeta, t, z) = \frac{\delta(\zeta, \omega k) \delta(\zeta, \omega^2k)}{\delta(\zeta, k)^2} \frac{\overline{r(\bar{k})}}{1 - r(k)\overline{r(\bar{k})}} z^{-2i\nu(\zeta)},
	\\
R_3^{\rm new}(\zeta, t, z) =  \frac{\delta(\zeta, k)^2}{\delta(\zeta, \omega k) \delta(\zeta, \omega^2k)} \frac{r(k)}{1 - r(k)\overline{r(\bar{k})}} z^{2i\nu(\zeta)},
	\\
R_4^{\rm new}(\zeta, t, z) = \frac{\delta(\zeta, \omega k) \delta(\zeta, \omega^2k)}{\delta(\zeta, k)^2}  \overline{r(\bar{k})} z^{-2i\nu(\zeta)},
\end{cases}
\end{align*}
where $k$ and $z$ are related by $\tilde{k} = \tilde{k}_0 + \frac{\epsilon z}{\rho_{\rm new}}$. 
Hence
$$R_1^{\rm new} = R_1 \bigg(\frac{k_0 + k^{-1}}{k_0 + k_0^{-1}}\bigg)^{2i\nu},$$
and so $q^{\rm new} = q$.
The proof of Theorem \ref{steepestdescentth} proceeds in the same way as before except that equation (\ref{mjinvasymptotics}) is replaced with 
\begin{align} \nonumber
 m_0(\zeta, t, k)^{-1} & =  \mathcal{C}D(\zeta, t) m^X\biggl(q(\zeta), \frac{\sqrt{\tau_{\rm new}}}{\epsilon}(\tilde{k} - \tilde{k}_0)\biggr)^{-1} D(\zeta, t)^{-1}\mathcal{C}
  	\\  \nonumber
& = I + \frac{B(\zeta, t)}{\sqrt{\tau_{\rm new}}(\tilde{k} - \tilde{k}_0)} + O(\tau_{\rm new}^{-1}), \qquad \  \zeta \in \mathcal{I}, \  |k - k_0| = \epsilon,
\end{align}
as $\tau_{\rm new} := t\rho_{\rm new}^2 \to \infty$. Hence, the contribution from the critical point at $k_0$ to $m(\zeta, t, K_1)$ is 
\begin{align}\nonumber
&\; \frac{1}{2\pi i} \int_{|k - k_0|=\epsilon } \frac{\hat{\mu}(\zeta,t,k) (m_0(\zeta,t,k)^{-1} - I) dk}{k - K_1}
	\\ \nonumber
= &\; \frac{1}{2\pi i} \int_{|k - k_0|=\epsilon } \frac{(m_0(\zeta,t,k)^{-1} - I) dk}{k - K_1}	
	\\\nonumber
& + \frac{1}{2\pi i}\int_{|k - k_0|=\epsilon } \frac{(\hat{\mu}(\zeta,t,k) - I) (m_0(\zeta,t,k)^{-1} - I) dk}{k -  K_1}
	\\\nonumber
= &\; \frac{1}{2\pi i}\int_{|k - k_0| = \epsilon} \frac{B(\zeta, t)}{\sqrt{\tau_{\rm new}}(\tilde{k} - \tilde{k}_0)} \frac{ dk}{k - K_1}
	\\\nonumber
& + O(\epsilon\tau_{\rm new}^{-1} )
+ O\bigl(\|\hat{\mu} - I\|_{\dot{L}^3(\hat{\Gamma})} \|m_0^{-1} - I\|_{L^{3/2}(|k - k_0|=\epsilon)} \bigr)
	\\\nonumber
= &\; \frac{1}{2\pi i}\frac{B(\zeta, t)}{\sqrt{\tau_{\rm new}}} \int_{|k - k_0| = \epsilon} \frac{2 dk}{(\frac{1}{k} - \frac{1}{k_0} - (k - k_0))(k - K_1)}
  + O(\epsilon\tau_{\rm new}^{-\frac{1+\alpha}{2}})
	\\ \label{newk0contribution}
= & - \frac{2}{\sqrt{\tau_{\rm new}}(\frac{1}{k_0^2} + 1)}\frac{B(\zeta, t)}{k_0 - K_1} + O(\epsilon\tau_{\rm new}^{-\frac{1+\alpha}{2}}).
\end{align}
This leads to the same formula for the asymptotics of $u(x,t)$ as above because $\frac{1}{2}(1 + k_0^{-2})\rho_{\rm new} = \rho$. 

Taking into account the correspondence of notations ($\nu$, $k_0$, and $\tilde{k}_0$ in the current paper correspond to $h_0$, $\kappa_0$, and $p_0$ in \cite{BS2013}, respectively), formulas (\ref{ufinalxi}) and (\ref{chi0def}) actually correct the coefficients $c_1$ and $c_4$ in the corresponding asymptotic formula (4.1) in \cite{BS2013}, where the contribution of the critical points to $m(K_1)$ was treated incorrectly. 
\end{rema}

\begin{rema}
We have assumed in Theorem \ref{mainth} that the set of singularities $\{k_j\}$ related to the presence of solitons is empty. If there are a finite number of points $\{k_j\}$ present at which the solution has simple poles, then the RH problem can be supplemented with residue conditions at these points, see Section 5 of \cite{L2013}. The supplemented RH problem can be mapped to a regular one coupled with a system of algebraic equations, see Proposition 2.4 of \cite{FI1996}. 
\end{rema}

\appendix
\section{$L^p$-Riemann--Hilbert problems} \label{RHapp}
\renewcommand{\theequation}{A.\arabic{equation}}
Since the jump contour for the RH problem associated with equation (\ref{DP}) on the half-line has nontransversal intersection points (see Figure \ref{Dns.pdf}), special care has to be taken when defining the notion of an $L^p$-RH problem. We will follow \cite{LCarleson} where a theory of $L^p$-RH problems with jumps across Carleson contours is developed using generalized Smirnoff classes.

Let $\mathcal{J}$ denote the collection of all subsets $\Gamma$ of the Riemann sphere $\hat{\C} = \C \cup \{\infty\}$ such that $\Gamma$ is homeomorphic to the unit circle and
\begin{align}\label{carlesondef}
 \sup_{z \in \Gamma \cap \C} \sup_{r > 0} \frac{|\Gamma \cap D(z, r)|}{r} < \infty,
\end{align}
where $D(z, r)$ denotes the disk of radius $r$ centered at $z$. Curves satisfying (\ref{carlesondef}) are called {\it Carleson curves}; all contours considered in this paper are Carleson.
Let $p \in [1,\infty)$. If $D$ is the bounded component of $\hat{\C} \setminus \Gamma$ where $\Gamma \in \mathcal{J}$ and $\infty \notin \Gamma$, then a function $f$ analytic in $D$ belongs to the {\it Smirnoff class $E^p(D)$} if there exists a sequence of rectifiable Jordan curves $\{C_n\}_1^\infty$ in $D$, tending to the boundary in the sense that $C_n$ eventually surrounds each compact subdomain of $D$, such that
\begin{align}\label{Epsup}
\sup_{n \geq 1} \int_{C_n} |f(z)|^p |dz| < \infty.
\end{align}
If $D$ is a subset of $\hat{\C}$ bounded by an arbitrary curve in $\mathcal{J}$, $E^p(D)$ is defined as the set of functions $f$ analytic in $D$ for which $f \circ \varphi^{-1} \in E^p(\varphi(D))$, where $\varphi(z) = \frac{1}{z - z_0}$ and $z_0$ is any point in $\C \setminus \bar{D}$. The subspace of $E^p(D)$ consisting of all functions $f \in E^p(D)$ such that $z f(z) \in E^p(D)$ is denoted by $\dot{E}^p(D)$. If $D = D_1 \cup \cdots \cup D_n$ is the union of a finite number of disjoint subsets of $\hat{\C}$ each of which is bounded by a curve in $\mathcal{J}$, then $E^p(D)$ and $\dot{E}^p(D)$ denote the set of functions $f$ analytic in $D$ such that $f|_{D_j} \in E^p(D_j)$ and $f|_{D_j} \in \dot{E}^p(D_j)$ for each $j$, respectively. We define $E^\infty(D)$ as the space of bounded analytic functions on $D$.

A {\it Carleson jump contour} is a connected subset $\Gamma$ of $\hat{\C}$ such that:
\begin{enumerate}[$(a)$]
\item $\Gamma \cap \C$ is the union of finitely many oriented arcs\footnote{A subset $\Gamma \subset \C$ is an {\it arc} if it is homeomorphic to a connected subset of the real line which contains at least two distinct points.} each pair of which have at most endpoints in common.

\item $\hat{\C} \setminus \Gamma$ is the union of two disjoint open sets $D_+$ and $D_-$ each of which has a finite number of simply connected components in $\hat{\C}$.

\item $\Gamma$ is the positively oriented boundary of $D_+$ and the negatively oriented boundary of $D_-$, i.e. $\Gamma = \partial D_+ = -\partial D_-$.

\item If $\{D_j^+\}_1^n$ and $\{D_j^-\}_1^m$ are the components of $D_+$ and $D_-$, then $\partial D_j^+ \in \mathcal{J}$ for $j = 1, \dots, n$, and $\partial D_j^- \in \mathcal{J}$ for $j = 1, \dots, m$.
\end{enumerate}

We henceforth make the following assumptions: (a) $p \in (1, \infty)$ and $n \geq 1$ is an integer, (b) $\Gamma  = \partial D_+ = -\partial D_-$ is a Carleson jump contour, and (c) $v: \Gamma \to GL(n, \C)$ is an $n \times n$-matrix valued function. 
We define $\dot{L}^p(\Gamma)$ as the set of all measurable functions on $\Gamma$ such that $|z - z_0|^{1 - \frac{2}{p}}h(z) \in L^p(\Gamma)$ for some (and hence all) $z_0 \in \C \setminus \Gamma$.
If $f \in \dot{E}^p(D_+)$ or  $f \in \dot{E}^p(D_-)$, the nontangential limits of $f(z)$ as $z$ approaches the boundary exist a.e. on $\Gamma$ and the boundary function belongs to $\dot{L}^p(\Gamma)$. Let $D = D_+ \cup D_-$. A {\it solution of the $L^p$-RH problem determined by $(\Gamma, v)$} is an $n \times n$-matrix valued function $m \in I + \dot{E}^p(D)$ such that the nontangential boundary values $m_\pm$ satisfy $m_+ = m_- v$ a.e. on $\Gamma$. 
 
\begin{lemma}[Uniqueness]\label{uniquelemma}
Suppose $1 \leq n \leq p$ and $\det v = 1$ a.e. on $\Gamma$.
If the solution of the $L^p$-RH problem determined by $(\Gamma, v)$ exists, then it is unique and has unit determinant.
\end{lemma}

If $h \in \dot{L}^p(\Gamma)$, then the Cauchy transform $\mathcal{C}h$ defined by
\begin{align}\label{Cauchytransform}
(\mathcal{C}h)(z) = \frac{1}{2\pi i} \int_\Gamma \frac{h(s)}{s - z} ds, \qquad z \in \C \setminus \Gamma,
\end{align}
satisfies $\mathcal{C}h \in \dot{E}^p(D)$. 
We denote the nontangential boundary values of $\mathcal{C}h$ from the left and right sides of $\Gamma$ by $\mathcal{C}_+ h$ and $\mathcal{C}_-h$ respectively. 
We henceforth fix a point $z_0 \in \C \setminus \Gamma$ and turn $\dot{L}^p(\Gamma)$ into a Banach space with the norm
\begin{align}\label{dotLpnorm}
\|h\|_{\dot{L}^p(\Gamma)} := \||\cdot - z_0|^{1 - \frac{2}{p}} h\|_{L^p(\Gamma)}.
\end{align}
Then $\mathcal{C}_+$ and $\mathcal{C}_-$ are bounded operators on $\dot{L}^p(\Gamma)$ and $\mathcal{C}_+ - \mathcal{C}_- = I$.
Given a function $w \in \dot{L}^p(\Gamma) \cap L^\infty(\Gamma)$, we define $\mathcal{C}_{w}: \dot{L}^p(\Gamma) + L^\infty(\Gamma) \to \dot{L}^p(\Gamma)$ by $\mathcal{C}_{w}(h) = \mathcal{C}_-(h w)$.
Then
\begin{align}\label{Cwnorm}
\|\mathcal{C}_w\|_{\mathcal{B}(\dot{L}^p(\Gamma))} \leq C \|w\|_{L^\infty(\Gamma)}.
\end{align}
where $C := \|\mathcal{C}_-\|_{\mathcal{B}(\dot{L}^p(\Gamma))} < \infty$ and $\mathcal{B}(\dot{L}^p(\Gamma))$ denotes the space of bounded linear operators on $\dot{L}^p(\Gamma)$.

\begin{lemma}\label{mulemma}
Suppose $w := v - I \in \dot{L}^p(\Gamma) \cap L^\infty(\Gamma)$.
If $m \in I + \dot{E}^p(D)$ satisfies the $L^p$-RH problem determined by $(\Gamma, v)$, then $\mu = m_- \in I + \dot{L}^p(\Gamma)$ satisfies 
\begin{align}\label{rhoeq}
\mu - I = \mathcal{C}_w(\mu)  \quad \text{in}\quad \dot{L}^p(\Gamma).
\end{align}
Conversely, if $\mu \in I + \dot{L}^p(\Gamma)$ satisfies (\ref{rhoeq}), then
$m = I + \mathcal{C}(\mu w) \in I + \dot{E}^p(D)$ satisfies the $L^p$-RH problem determined by $(\Gamma, v)$. 
\end{lemma}

\begin{lemma}\label{EpCnlemma}
Let $D$ be an open subset of $\hat{\C}$ bounded by a curve $\Gamma \in \mathcal{J}$ with $\infty \in \Gamma$. Let $z_0 \in \C \setminus \bar{D}$ and let $f:D \to \C$ be an analytic function. 
 Then $f \in \dot{E}^p(D)$ if and only if there exist curves $\{C_n\}_1^\infty \subset \mathcal{J}$ in $D$, tending to $\Gamma$ in the sense that $C_n$ eventually surrounds each compact subset of $D$, such that
\begin{align}\label{Epdotsupz0}
\sup_{n \geq 1} \int_{C_n} |z - z_0|^{p-2} |f(z)|^p |dz| < \infty.
\end{align}
\end{lemma}

 \begin{lemma}[Contour deformation]\label{deformationlemma}
Let $\gamma \in \mathcal{J}$. Suppose that, reversing the orientation on a subcontour if necessary, $\hat{\Gamma} = \Gamma \cup \gamma$ is a Carleson jump contour. 
Let $B_+$ and $B_-$ be the two components of $\hat{\C} \setminus \gamma$. Let $\hat{D}_\pm$ be the open sets such that $\hat{\C} \setminus \hat{\Gamma} = \hat{D}_+ \cup \hat{D}_-$ and $\partial \hat{D}_+ = - \partial \hat{D}_- = \hat{\Gamma}$. 
Let $\hat{D} = \hat{D}_+ \cup \hat{D}_-$.
Let $\gamma_+$ and $\gamma_-$ be the parts of $\gamma$ that belong to the boundaries of $\hat{D}_+ \cap B_+$ and $\hat{D}_- \cap B_+$, respectively. 
Suppose $v: \Gamma \to GL(n, \C)$. Suppose $m_0:\hat{D} \cap B_+ \to GL(n,\C)$ satisfies
\begin{align}\label{m0m0inv}
m_0, m_0^{-1} \in I + \dot{E}^p(\hat{D} \cap B_+) \cap E^\infty(\hat{D} \cap B_+).
\end{align}
Define $\hat{v}:\hat{\Gamma} \to GL(n, \C)$ by
\begin{align*}
\hat{v} 
=  \begin{cases}
 m_{0-} v m_{0+}^{-1} & \text{on} \quad  \Gamma \cap B_+, \\
m_{0+}^{-1} & \text{on} \quad \gamma_+, \\
m_{0-} & \text{on} \quad \gamma_-, \\
v & \text{on} \quad \Gamma \cap B_-.
\end{cases}
\end{align*}
Let $m$ and $\hat{m}$ be related by
\begin{align}\label{hatmdefmm0}
\hat{m} = \begin{cases}
mm_0^{-1} & \text{on} \quad \hat{D} \cap B_+,\\
m & \text{on} \quad \hat{D} \cap B_-.
\end{cases}
\end{align}
Then $m(z)$ satisfies the $L^p$-RH problem determined by $(\Gamma,v)$ if and only if $\hat{m}(z)$ satisfies the $L^p$-RH problem determined by $(\hat{\Gamma}, \hat{v})$.
\end{lemma}

Proofs of the above statements can be found in \cite{LCarleson}. We will also need the following uniqueness result for row vector solutions. 

\begin{lemma}\label{rowvectorlemma}
Suppose $1 \leq n \leq p$ and $\det v = 1$ a.e. on $\Gamma$. Suppose the $L^p$-RH problem determined by $(\Gamma, v)$ has a unique solution $m$. 
If $N$ is a row vector solution of the $L^p$-RH problem determined by $(\Gamma, v)$ in the sense that $N \in (1, 1, \dots, 1) + \dot{E}^p(D)$ and $N_+ = N_- v$ a.e. on $\Gamma$, then $N = (1, 1, \dots, 1)m$.
\end{lemma}
\begin{proof}
Let $\hat{m}$ denote the $n \times n$-matrix valued function obtained from $m$ by replacing the first row with the row vector $N$. Then the  $n \times n$-matrix valued function $\tilde{m} \in I + \dot{E}^p(D)$ defined by
$$\tilde{m} = \begin{pmatrix} 1 & -1 & \cdots & -1 \\ 0 & 1 & \cdots & 0 \\ \vdots & \vdots &  & \vdots \\ 0 & 0 & \cdots & 1 \end{pmatrix} \hat{m}$$
satisfies the $L^p$-RH problem determined by $(\Gamma, v)$. Hence $\tilde{m} = m$ by Lemma \ref{uniquelemma}. Consequently, $N = (1, 1, \dots, 1)\tilde{m} = (1, 1, \dots, 1)m$. 
\end{proof}

\section{The solution on a cross} 
\renewcommand{\theequation}{B.\arabic{equation}}

\begin{figure}
\begin{center}
 \begin{overpic}[width=.35\textwidth]{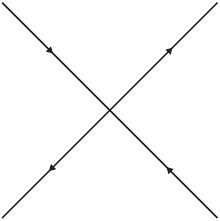}
 \put(67,81){$X_1$}
 \put(24,81){$X_2$}
 \put(24,16){$X_3$}
 \put(67,16){$X_4$}
 \end{overpic}
   \bigskip
   \begin{figuretext}\label{X.pdf}
      The contour $X = X_1 \cup \cdots \cup X_4$.
      \end{figuretext}
   \end{center}
\end{figure}

Consider the cross $X = X_1 \cup \cdots \cup X_4 \subset \hat{\C}$ where
\begin{align} \nonumber
&X_1 = \bigl\{ue^{\frac{i\pi}{4}}\, \big| \, 0 \leq u \leq \infty\bigr\}, && 
X_2 = \bigl\{ue^{\frac{3i\pi}{4}}\, \big| \, 0 \leq u \leq \infty\bigr\},  
	\\ \label{Xdef}
&X_3 = \bigl\{ue^{-\frac{3i\pi}{4}}\, \big| \, 0 \leq u \leq \infty\bigr\}, && 
X_4 = \bigl\{ue^{-\frac{i\pi}{4}}\, \big| \, 0 \leq u \leq \infty\bigr\},
\end{align}
and $X$ is oriented as in Figure \ref{X.pdf}. 
Let $\D \subset \C$ denote the open unit disk and define  the function $\nu:\D \to (0,\infty)$ by 
$\nu(q) = -\frac{1}{2\pi} \log(1 - |q|^2)$.
We consider the following family of $L^3$-RH problems parametrized by $q \in \D$:
\begin{align}\label{RHmc}
\begin{cases} m^X(q, \cdot) \in I + \dot{E}^3(\hat{\C} \setminus X), 
	\\
m_+^X(q, z) =  m_-^X(q, z) v^X(q, z) \quad \text{for a.e.} \ z \in X, 
\end{cases} 
\end{align}
where the jump matrix $v^X(q, z)$ is defined by\footnote{Throughout the paper, complex powers and logarithms are defined using the principal branch: If $z, a \in \C$ and $z \neq 0$, then $\log z := \log|z| + i\Arg{z}$ and $z^a := e^{a\log z}$, where $\Arg z \in (-\pi, \pi]$ denotes the principal value of $\arg z$.}
\begin{align*} 
v^X(q, z) = \begin{cases}
\begin{pmatrix} 1 & 0 & 0	\\
  q z^{-2i\nu(q)} e^{\frac{iz^2}{2}}	& 1 & 0 \\ 0 & 0 & 1 \end{pmatrix}, &   z \in X_1, 
  	\\
\begin{pmatrix} 1 & - \frac{\overline{q}}{1 - |q|^2} z^{2i\nu(q)}e^{-\frac{iz^2}{2}}	& 0\\
0 & 1 & 0 \\ 0 & 0 & 1  \end{pmatrix}, &  z \in X_2, 
	\\
\begin{pmatrix} 1 &0 & 0\\
- \frac{q}{1 - |q|^2}z^{-2i\nu(q)} e^{\frac{iz^2}{2}}	& 1 & 0 \\ 0 & 0 & 1 \end{pmatrix}, &  z \in X_3,
	\\
 \begin{pmatrix} 1	& \overline{q} z^{2i\nu(q)}e^{-\frac{iz^2}{2}} & 0	\\
0	& 1 & 0 \\ 0 & 0 & 1  \end{pmatrix}, &  z \in X_4.
\end{cases}
\end{align*}

The RH problem (\ref{RHmc}) can be solved explicitly in terms of parabolic cylinder functions \cite{I1981} and this leads to the following standard result.

\begin{theo}\label{crossth}
  The $L^3$-RH problem (\ref{RHmc}) has a unique solution $m^X(q, z)$ for each $q \in \D$. This solution satisfies
\begin{align}\label{mcasymptotics}
  m^X(q, z) = I + \frac{i}{z}\begin{pmatrix} 0 & -\beta^X(q) & 0 \\ \overline{\beta^X(q)} & 0 & 0 \\ 0 & 0 & 0 \end{pmatrix} + O\biggl(\frac{1}{z^2}\biggr), \qquad z \to \infty,  \  q \in \D, 
\end{align}  
where the error term is uniform with respect to $\arg z \in [0, 2\pi]$ and $q$ in compact subsets of $\D$, and the function $\beta^X(q)$ is defined by
\begin{align}\label{betacdef}
\beta^X(q) = \sqrt{\nu(q)} e^{i\left(\frac{\pi}{4} - \arg q + \arg \Gamma(i\nu(q)\right)}, \qquad q \in \D.
\end{align}
Moreover, for each compact subset $K$ of $\D$, 
$$\sup_{q \in K} \sup_{z \in \C \setminus X} |m^X(q, z)| < \infty.$$
\end{theo}

\section{Proofs of claims}\label{proofapp}
\renewcommand{\theequation}{C.\arabic{equation}}
This appendix presents the proofs of the five claims used in the proof of Theorem \ref{steepestdescentth}.

\subsection{Proof of Claim \ref{claim1}.}
We first assume $k \in k_0 + X^\epsilon$. Then
\begin{align*}
\hat{w}(\zeta, t,k)
& = m_{0-}(\zeta, t, k) v(\zeta, t, k) m_{0+}(\zeta,t,k)^{-1} - I
	\\
& = m_{0-}(\zeta, t, k)  u(\zeta,t,k) m_{0+}(\zeta,t,k)^{-1},
\end{align*}
where
$$
u(\zeta,t,k) := v(\zeta, t, k) - \mathcal{C} D(\zeta,t) v^X\biggl(q(\zeta), -  \frac{\sqrt{\tau}}{\epsilon}(k-k_0)\biggr) D(\zeta,t)^{-1} \mathcal{C}.$$
The functions $m_{0+}(\zeta, t, k)$ and $m_{0-}(\zeta, t, k)$ are uniformly bounded for $t>0$, $\zeta \in \mathcal{I}$, and $k \in k_0 + X^\epsilon$. Therefore, it is enough to prove that 
\begin{align}\label{Ujbound}
u(\zeta,t,k) = O\bigl(\tau^{-\frac{\alpha}{2}} e^{-\frac{\tau}{24\epsilon^2}|k-k_0|^2}\bigr), \qquad \tau \to \infty, \  \zeta \in \mathcal{I}, \  k \in k_0 + X^\epsilon,
\end{align}
uniformly with respect to $(\zeta, k)$.
Introducing the function $u_0$ by
\begin{align*}
u_0(\zeta,t,z) & =  \mathcal{C}  u\biggl(\zeta,t,k_0 - \frac{\epsilon z}{\rho}\biggr) \mathcal{C} 
	\\
& = v_0(\zeta, t, z) - D(\zeta,t) v^X(q(\zeta), \sqrt{t} z) D(\zeta,t)^{-1},
\end{align*}
we can rewrite the condition (\ref{Ujbound}) as follows:
\begin{align}\label{ujbound}
u_0(\zeta,t,z) = O\bigl(\tau^{-\frac{\alpha}{2}} e^{-\frac{t|z|^2}{24}}\bigr), \qquad \tau \to \infty, \  \zeta \in \mathcal{I}, \  z \in X^{\rho},
\end{align}
uniformly with respect to $(\zeta, z)$ in the given ranges.
Using that
\begin{align*}
& D(\zeta,t) v^X\bigl(q(\zeta), \sqrt{t} z\bigr) D(\zeta,t)^{-1}
	\\
&\qquad  = \begin{cases}
 \begin{pmatrix} 1 & 0 & 0	\\
  q(\zeta)  z^{-2i\nu(\zeta)} e^{\frac{itz^2}{2}}e^{ t\phi(\zeta, 0)}	& 1 & 0 \\ 0 & 0 & 1\end{pmatrix}, &   z \in X_1, 
  	\\
\begin{pmatrix} 1 & - \frac{\overline{q(\zeta)}}{1 - |q(\zeta)|^2} z^{2i\nu(\zeta)}e^{-\frac{itz^2}{2}}e^{-t\phi(\zeta, 0)} & 0	\\
0 & 1 & 0 \\ 0 & 0 & 1 \end{pmatrix}, &  z \in X_2, 
	\\
\begin{pmatrix} 1 &0 & 0\\
- \frac{q(\zeta)}{1 - |q(\zeta)|^2} z^{-2i\nu(\zeta)} e^{\frac{itz^2}{2}}e^{t\phi(\zeta, 0)} & 1 & 0 \\ 0 & 0 & 1 \end{pmatrix}, &  z \in X_3,
	\\
\begin{pmatrix} 1	&  \overline{q(\zeta)} z^{2i\nu(\zeta)}e^{-\frac{itz^2}{2}} e^{-t\phi(\zeta, 0)} & 0 \\
0 & 1 & 0 \\ 0 & 0 & 1 \end{pmatrix}, &  z \in X_4,	
\end{cases}
\end{align*}
equation (\ref{ujbound}) follows from the assumptions (\ref{smallcrossjump})-(\ref{Lipschitzconditions}). Indeed, we will give the details of the proof of (\ref{ujbound}) in the case of $z \in X_1^{\rho}$; the other cases are similar.

Let $z \in X_1^{\rho}$. In this case only the $(21)$ entry of $u_0(\zeta,t,z)$ is nonzero and using that $\arg z = \frac{\pi}{4}$ and $\sup_{\zeta \in \mathcal{I}} |q(\zeta)| < 1$, we find
\begin{align} \nonumber
|(u_0(\zeta, t, z))_{21}| = &\;
\big|R_1(\zeta, t, z)z^{-2i\nu(\zeta)}e^{t\phi(\zeta, z)}
- q(\zeta) z^{-2i\nu(\zeta)}e^{\frac{itz^2}{2}} e^{t\phi(\zeta, 0)}\big|
	\\ \nonumber
= & \; |z^{-2i\nu(\zeta)}| \bigl|R_1(\zeta, t, z)e^{t\hat{\phi}(\zeta, z)} - q(\zeta)\bigr| |e^{t\phi(\zeta, 0)}| e^{-\frac{t|z|^2}{2}} 
  	\\ \nonumber
 \leq &\; e^{\frac{\pi \nu(\zeta)}{2}}\Bigl(\bigl|R_1(\zeta, t, z) - q(\zeta)\bigr|e^{t\re \hat{\phi}(\zeta, z)}
 	\\ \label{Westimate}
& + |q(\zeta)| \bigl|e^{t\hat{\phi}(\zeta, z)} - 1\bigr|\Bigr)
 e^{-\frac{t|z|^2}{2}}, \qquad \zeta \in \mathcal{I}, \  t > 0, \  z \in X_1^{\rho},
\end{align}
where $\hat{\phi}(\zeta, z) = \phi(\zeta,z)- \phi(\zeta, 0) - \frac{iz^2}{2}$.
The simple estimate
$$|e^w -1| = \biggl|\int_0^1 w e^{sw} ds\biggr| \leq |w| \max_{s \in [0,1]}e^{s\re w}, \qquad w \in \C,$$
yields the inequality
\begin{align}\label{ewminus1estimate}  
  |e^w - 1| \leq |w| \max(1, e^{\re w}), \qquad w \in \C.
\end{align}
On the other hand, by (\ref{phiassumptions}) and (\ref{rephiestimatea}),
\begin{align}\label{rehatphi}
\re \hat{\phi}(\zeta, z) = \re \phi(\zeta,z) + \frac{|z|^2}{2} \leq \frac{|z|^2}{4}, \qquad \zeta \in \mathcal{I}, \  z \in X_1^{\rho}.
\end{align}
Using (\ref{ewminus1estimate}), (\ref{rehatphi}), and the fact that $\sup_{\zeta \in \mathcal{I}} |q(\zeta)| < 1$ in (\ref{Westimate}), we find
\begin{align*}
 |(u_0(\zeta, t, z))_{21}| \leq &\;  e^{\frac{\pi \nu(\zeta)}{2}}\Bigl(\bigl|R_1(\zeta, t, z) - q(\zeta)\bigr| 
  + t\bigl|\hat{\phi}(\zeta, z)\bigr|\Bigr) e^{-\frac{t|z|^2}{4}}, 
  	\\
& \hspace{5cm}   \zeta \in \mathcal{I}, \  t > 0,\  z \in X_1^{\rho}. 
\end{align*}  
By (\ref{Phiz3estimate}), (\ref{Lipschitzconditions}), and the fact that $\sup_{\zeta \in \mathcal{I}} |\nu(\zeta)| < \infty$, the right-hand side is of order
\begin{align} \nonumber
& O\biggl(\biggl( \frac{L |z|^\alpha}{\rho^\alpha} + \frac{tC|z|^3}{\rho} \biggr) e^{-\frac{t|z|^2}{4}}  \biggr)
= O\biggl( \biggl(\frac{(t|z|^2)^{\alpha/2}}{\tau^{\alpha/2}} + \frac{(t|z|^2)^{3/2}}{\tau^{1/2}}\biggr) e^{-\frac{t|z|^2}{12}}\biggr) 
 	\\
&  = O\biggl( \biggl(\frac{1}{\tau^{\alpha/2}} + \frac{1}{\tau^{1/2}}\biggr) e^{-\frac{t|z|^2}{24}}\biggr), \qquad \tau \to \infty, \  \zeta \in \mathcal{I}, \  z \in X_1^{\rho},
\end{align}
uniformly with respect to $(\zeta, z)$ in the given ranges. This proves (\ref{ujbound}) in the case of $z \in X_1^{\rho}$.

Now let $k \in -k_0 + X^\epsilon$. Then
\begin{align*}
\hat{w}(\zeta, t,k)
& = m_{0-}(\zeta, t, k) v(\zeta, t, k) m_{0+}(\zeta,t,k)^{-1} - I
	\\
& = m_{0-}(\zeta, t, k)  u(\zeta,t,k) m_{0+}(\zeta,t,k)^{-1},
\end{align*}
where
$$u(\zeta,t,k) := v(\zeta, t, k) - \mathcal{C}\overline{\check{D}(\zeta,t) v^X\biggl(\check{q}(\zeta), \frac{\sqrt{\tau}}{\epsilon}\overline{(k + k_0)}\biggr) \check{D}(\zeta,t)^{-1}}\mathcal{C}.$$
The functions $m_{0+}(\zeta, t, k)$ and $m_{0-}(\zeta, t, k)$ are uniformly bounded for $t>0$, $\zeta \in \mathcal{I}$, and $k \in -k_0 + X^\epsilon$. Therefore, it is enough to prove that 
\begin{align}\label{Ujbound2}
u(\zeta,t,k) = O\bigl(\tau^{-\frac{\alpha}{2}} e^{-\frac{\tau}{8\epsilon^2}|k+k_0|^2}\bigr), \qquad \tau \to \infty, \  \zeta \in \mathcal{I}, \  k \in -k_0 + X^\epsilon,
\end{align}
uniformly with respect to $(\zeta, k)$.
Introducing the function $u_0$ by
\begin{align*}
u_0(\zeta,t,z) & = \mathcal{C}\overline{u\biggl(\zeta,t, \overline{- k_0 + \frac{\epsilon z}{\rho}}\biggr)}\mathcal{C}
	\\
& = \check{v}_0(\zeta, t, z) - \check{D}(\zeta,t) v^X(\check{q}(\zeta), \sqrt{t} z) \check{D}(\zeta,t)^{-1},
\end{align*}
we can rewrite the condition (\ref{Ujbound2}) as follows:
\begin{align*}
u_0(\zeta,t,z) = O\bigl(\tau^{-\frac{\alpha}{2}} e^{-\frac{t|z|^2}{8}}\bigr), \qquad \tau \to \infty, \  \zeta \in \mathcal{I}, \  z \in X^{\rho},
\end{align*}
uniformly with respect to $(\zeta, z)$ in the given ranges. The rest of the proof is as in the case of $k \in k_0 + X^\epsilon$.
\proofend

\subsection{Proof of Claim \ref{claim2}.}
In view of the symmetries (\ref{symmetriesa}) and (\ref{symmetriesb}),
\begin{align}\nonumber
\|\hat{w}(\zeta, t, \cdot)\|_{\dot{L}^3(\hat{\Gamma})} 
=  O\Big(&\|\hat{w}(\zeta, t, \cdot)\|_{\dot{L}^3(\Gamma')} + \|m_0(\zeta, t, \cdot)^{-1} - I \|_{L^3(|k - k_0| = \epsilon)} 
	\\ \nonumber
&+ \|m_0(\zeta, t, \cdot)^{-1} - I \|_{L^3(|k + k_0| = \epsilon)} 
	\\ \label{whatnorm}
&  + \|\hat{w}(\zeta, t, \cdot) \|_{L^3(k_0 + X^\epsilon)}  + \|\hat{w}(\zeta, t, \cdot) \|_{L^3(-k_0 + X^\epsilon)}\Big).
\end{align}
On $\Gamma'$, the matrix $\hat{w}$ is given by either $v - I$ or $m_0(v - I)m_0^{-1}$.
Hence $\|\hat{w}(\zeta, t, \cdot)\|_{\dot{L}^3(\Gamma')} = O(\epsilon^{1/3}\tau^{-1})$ by the assumption (\ref{wL12inftya}).
Moreover, by (\ref{mcasymptotics}), $m^X(q, z) = I + O\bigl(\frac{1}{z}\bigr)$ as $z \to \infty$ uniformly with respect to the argument of $z$ and with respect to $q$ in compact subsets of $\D$. Hence, as the entries of $D(\zeta, t)$ have unit modulus,
\begin{align}\nonumber
&\|m_0(\zeta, t, k)^{-1} - I \|_{L^p(|k - k_0| = \epsilon)} 
 	\\ \nonumber
& = 
\biggl \|\mathcal{C}D(\zeta, t) \biggl[m^X\biggl(q(\zeta),  - \frac{\sqrt{\tau}}{\epsilon}( k - k_0)\biggr)^{-1}  - I \biggr]D(\zeta, t)^{-1}\mathcal{C}\biggr \|_{L^p(|k - k_0| = \epsilon)}  
	\\ \nonumber
& = \begin{cases} O(\epsilon^{1/p} \tau^{-1/2}), & p \in [1,\infty), \\ 
O(\tau^{-1/2}), & p = \infty,
\end{cases}
\end{align}
uniformly with respect to $\zeta \in \mathcal{I}$. This proves (\ref{star}). The third term on the right-hand side of (\ref{whatnorm}) can be estimated in a similar way. The last two terms in (\ref{whatnorm}) can be estimated using (\ref{hatwestimate}). This yields (\ref{hatwestimatea}). The proof of (\ref{hatwestimateb}) uses the assumption (\ref{wL12inftyb}) and is similar.

In order to prove (\ref{hatwestimatec}), we note that (\ref{hatwestimate}) implies
\begin{align} \nonumber
  \|\hat{w}(\zeta, t, \cdot)\|_{L^p(k_0 + X^\epsilon)} 
 & = O\biggl(  \tau^{-\frac{\alpha}{2}} \bigg(\int_{k_0 + X^\epsilon} e^{-\frac{p\tau}{24\epsilon^2} |k-k_0|^2} |dk|\bigg)^{\frac{1}{p}} \biggr)
  	\\ \label{hatwL1estimate}
& = O\biggl(\tau^{-\frac{\alpha}{2}}  \bigg(\int_0^\epsilon e^{-\frac{p\tau}{24\epsilon^2} u^2} du\bigg)^{\frac{1}{p}} \biggr), \qquad \tau \to \infty, \  \zeta \in \mathcal{I}.
\end{align}
Letting $v = \frac{p\tau}{24\epsilon^2} u^2$ we find 
\begin{align}\label{int0epsilonestimate}
\int_0^\epsilon e^{-\frac{p\tau}{24\epsilon^2} u^2} du
\leq \int_0^\infty  e^{-\frac{p\tau}{24\epsilon^2} u^2} du 
= \frac{\epsilon\sqrt{6}}{\sqrt{p\tau}} \int_0^\infty \frac{e^{-v}}{\sqrt{v}} dv = \frac{\epsilon \sqrt{6\pi}}{\sqrt{p\tau}}.
\end{align}
Equations (\ref{hatwL1estimate}) and (\ref{int0epsilonestimate}) yield (\ref{hatwestimatec}).
\proofend

\subsection{Proof of Claim \ref{claim3}.}
By (\ref{Cwnorm}) and (\ref{hatwestimateb}),
\begin{align}\label{Chatwnorm}
\|\hat{\mathcal{C}}_{\hat{w}}\|_{\mathcal{B}(\dot{L}^3(\hat{\Gamma}))} \leq C \|\hat{w}\|_{L^\infty(\hat{\Gamma})}  
= O(\tau^{-\frac{\alpha}{2}}), \qquad \tau \to \infty.
\end{align}
This proves the claim. 
\proofend

\subsection{Proof of Claim \ref{claim4}.}
The Neumann series 
\begin{align}\label{neumannseries}
(I - \hat{\mathcal{C}}_{\hat{w}})^{-1} = \sum_{j=0}^\infty \hat{\mathcal{C}}_{\hat{w}}^j
\end{align}
implies that
\begin{align*}
\|(I - \hat{\mathcal{C}}_{\hat{w}})^{-1}\|_{\mathcal{B}(\dot{L}^3(\hat{\Gamma}))} 
& \leq \sum_{j=0}^\infty \|\hat{\mathcal{C}}_{\hat{w}}\|_{\mathcal{B}(\dot{L}^3(\hat{\Gamma}))}^j
 = \frac{1}{1 - \|\hat{\mathcal{C}}_{\hat{w}}\|_{\mathcal{B}(\dot{L}^3(\hat{\Gamma}))}}.
%	\\
%& \leq \frac{1}{1 - c\|\hat{w}\|_{L^\infty(\hat{\Gamma})}}.
\end{align*}
Now (\ref{cauchysingularbound}) and the Sokhotski-Plemelj formula $\mathcal{C}_- = \frac{1}{2}(-I + \mathcal{S}_\Gamma)$ show that 
$$\sup_{\zeta \in \mathcal{I}} \|\hat{\mathcal{C}}_-\|_{\mathcal{B}(\dot{L}^3(\hat{\Gamma}))} < \infty.$$ 
Thus,
\begin{align*}
\|\hat{\mu} - I\|_{\dot{L}^3(\hat{\Gamma})} & = 
\|(I- \hat{\mathcal{C}}_{\hat{w}})^{-1}\hat{\mathcal{C}}_{\hat{w}}I\|_{\dot{L}^3(\hat{\Gamma})} 
	\\
&  \leq \|(I - \hat{\mathcal{C}}_{\hat{w}})^{-1}\|_{\mathcal{B}(\dot{L}^3(\hat{\Gamma}))}
\|\hat{\mathcal{C}}_-(\hat{w})\|_{\dot{L}^3(\hat{\Gamma})}
\leq \frac{C\|\hat{w}\|_{\dot{L}^3(\hat{\Gamma})}}{1 - \|\hat{\mathcal{C}}_{\hat{w}}\|_{\mathcal{B}(\dot{L}^3(\hat{\Gamma}))}}.
\end{align*}
In view of (\ref{hatwestimatea}) and (\ref{Chatwnorm}), this gives (\ref{muhatestimate}).
\proofend

\subsection{Proof of Claim \ref{claim5}.}
Uniqueness follows from Lemma \ref{uniquelemma} since $\det \hat{v} = 1$. Moreover, equation (\ref{hatmudef}) implies that $\hat{\mu} - I = \hat{\mathcal{C}}_{\hat{w}} \hat{\mu}$. Hence, by Lemma \ref{mulemma}, $\hat{m} = I + \hat{\mathcal{C}}(\hat{\mu} \hat{w})$ satisfies the $L^3$-RH problem (\ref{RHmhat}). 
\proofend

\nocite{*}
\bibliographystyle{cdraifplain}
\bibliography{bibliography}

\end{document}